\documentclass[11pt]{amsart}
\usepackage{amsfonts}
\usepackage{mathrsfs}
\usepackage{amscd}
\usepackage{amsmath}
\usepackage{amssymb}
\usepackage{latexsym}
\usepackage{lscape}
\usepackage{xypic}

\newtheorem{thm}{Theorem}[section]
\newtheorem{prop}[thm]{Proposition}
\newtheorem{lem}[thm]{Lemma}

\newtheorem{lem-def}[thm]{Lemma-Definition}
\newtheorem{cor}[thm]{Corollary}
\theoremstyle{remark}

\newtheorem{rmk}{Remark}[section]
\theoremstyle{definition}
\newtheorem{dfn}{Definition}[section]
\newtheorem{Thm}[]{Theorem}

\numberwithin{equation}{subsection}

\input xy
\xyoption{all}

\newcommand{\quash}[1]{}  
\newcommand{\nc}{\newcommand}
\nc{\on}{\operatorname}

\newcommand{\frakc}{{\mathfrak c}}

\newcommand{\frakl}{{\mathfrak l}}

\newcommand{\fraks}{{\mathfrak s}}

\newcommand{\fraku}{{\mathfrak u}}

\newcommand{\frakD}{{\mathfrak D}}

\newcommand{\frakR}{{\mathfrak R}}

\newcommand{\frakX}{{\mathfrak X}}

\newcommand{\bbA}{{\mathbb A}}
\newcommand{\bbB}{{\mathbb B}}

\newcommand{\bbG}{{\mathbb G}}

\newcommand{\bbP}{{\mathbb P}}
\newcommand{\bbQ}{{\mathbb Q}}
\newcommand{\bbR}{{\mathbb R}}

\newcommand{\bbT}{{\mathbb T}}

\newcommand{\bbZ}{{\mathbb Z}}
\newcommand{\calA}{{\mathcal A}}
\newcommand{\calB}{{\mathcal B}}

\newcommand{\calE}{{\mathcal E}}
\newcommand{\calF}{{\mathcal F}}
\newcommand{\calG}{{\mathcal G}}
\newcommand{\calH}{{\mathcal H}}
\newcommand{\calI}{{\mathcal I}}

\newcommand{\calL}{{\mathcal L}}
\newcommand{\calM}{{\mathcal M}}
\newcommand{\calN}{{\mathcal N}}
\newcommand{\calO}{{\mathcal O}}
\newcommand{\calP}{{\mathcal P}}

\newcommand{\calS}{{\mathcal S}}
\newcommand{\calT}{{\mathcal T}}
\newcommand{\calU}{{\mathcal U}}
\newcommand{\calV}{{\mathcal V}}

\newcommand{\calX}{{\mathcal X}}

\newcommand{\calZ}{{\mathcal Z}}

\nc{\al}{{\alpha}} \nc{\be}{{\beta}}
\newcommand{\ga}{{\gamma}}
\nc{\ve}{{\varepsilon}}
\nc{\Ga}{{\Gamma}}
\newcommand{\la}{{\lambda}}
\nc{\La}{{\Lambda}}

\nc{\ad}{{\on{ad}}}
\newcommand{\Ad}{{\on{Ad}}}
\nc{\Adm}{{\on{Adm}}}
\nc{\aff}{{\on{aff}}}
\nc{\Aff}{{\mathbf{Aff}}}
\newcommand{\Aut}{{\on{Aut}}}
\nc{\Bun}{{\on{Bun}}}

\nc{\der}{{\on{der}}}

\nc{\diag}{{\on{diag}}}
\newcommand{\End}{{\on{End}}}
\nc{\Fl}{{\calF\ell}}
\newcommand{\Gal}{{\on{Gal}}}
\newcommand{\Gr}{{\on{Gr}}}
\newcommand{\Hom}{{\on{Hom}}}
\nc{\IC}{{\on{IC}}}
\newcommand{\id}{{\on{id}}}
\nc{\Id}{{\on{Id}}}

\nc{\Ind}{{\on{Ind}}}
\newcommand{\Lie}{{\on{Lie\ \!\!}}}
\newcommand{\Pic}{{\on{Pic}}}
\newcommand{\pr}{{\on{pr}}}
\newcommand{\Res}{{\on{Res}}}
\nc{\res}{{\on{res}}}
\newcommand{\s}{{\on{sc}}}

\newcommand{\Spec}{{\on{Spec\ \!\!}}}

\nc{\tr}{{\on{tr}}}
\newcommand{\Tr}{{\on{Tr}}}

\newcommand{\GL}{{\on{GL}}}
\nc{\GSp}{{\on{GSp}}}
\nc{\GU}{{\on{GU}}}

\nc{\SL}{{\on{SL}}}
\nc{\SU}{{\on{SU}}}
\nc{\SO}{{\on{SO}}}

\newcommand{\bGr}{{\overline{\Gr}}}
\nc{\bFl}{{\overline{\Fl}}} \nc{\bU}{{\overline{U}}}
\nc{\wGr}{{\widetilde{\Gr}}} \nc{\wJG}{{\widetilde{\calL^+\calG}}}
\nc{\wLG}{{\widetilde{\calL\calG}}}

\nc{\ppars}{(\!(s)\!)}
\newcommand{\ppart}{(\!(t)\!)}
\newcommand{\pparu}{(\!(u)\!)}

\def\xcoch{\mathbb{X}_\bullet}
\def\xch{\mathbb{X}^\bullet}

\title{On the coherence conjecture of Pappas and Rapoport}
\author{Xinwen Zhu}
\address{Xinwen Zhu: Department of Mathematics, Northwestern University, Evanston, IL, 60208, USA.}
\email{xinwenz@math.northwestern.edu}

\begin{document}
\date{December 2010}
\maketitle
\begin{abstract}
We prove the (generalized) coherence conjecture of Pappas and
Rapoport proposed in \cite{PR1}. As a corollary, one of the main
theorems in \cite{PR2}, which describes the geometry of the special
fibers of the local models for ramified unitary groups, holds
unconditionally. Our proof is based on the study of the geometry (in
particular certain line bundles and $\ell$-adic sheaves) of the
global Schubert varieties, which are the equal characteristic
counterparts of the local models.
\end{abstract}
\tableofcontents

\section{Introduction}

The goal of this paper is to prove the coherence conjecture of
Pappas and Rapoport as proposed in \cite{PR1}. The precise
formulation of the conjecture is a little bit technical and will be
given in \S \ref{the conj}. In this introduction, we would like to
describe a vague form of this conjecture, to convey the ideas behind
it and to outline the proofs.

The coherence conjecture was proposed by Pappas and Rapoport in
order to understand the special fibers of local models. Local models
were systematically introduced by Rapoport and Zink in \cite{RZ}
(special cases were constructed earlier by Deligne-Pappas \cite{DP}
and independently by de Jong \cite{dJ}) as a tool to analyze the
\'{e}tale local structure of certain integral models of (PEL-type)
Shimura varieties with parahoric level structures over $p$-adic
fields. Unlike the Shimura varieties themselves, which are usually
moduli spaces of abelian varieties, local models are defined in
terms of linear algebra and therefore are much easier to study. For
example, using local models, G\"{o}rtz (see \cite{Go1,Go2}) proved
the flatness of certain PEL-type Shimura varieties associated to
unramified unitary groups and symplectic groups (some special cases
were obtained in earlier works \cite{CN,dJ,DP}). On the other hand,
a discovery of G. Pappas (cf. \cite{Pa}) showed that the originally
defined integral models in \cite{RZ} are usually not flat when the
groups are ramified. Therefore, nowadays the (local) models defined
in \cite{RZ} are usually called the naive models. In a series of
papers (\cite{PR-1,PR0,PR2}), Pappas and Rapoport investigated the
corrected definition of flat local models. The easiest definition of
these local models is by taking the flat closures of the generic
fibers in the naive local models. Usually, an integral model defined
in this way is not useful since the moduli interpretation is lost
and therefore it is very difficult to study the special fiber, etc
(in fact a considerable part \cite{PR-1,PR0,PR2} is devoted in an
attempt to cutting out the correct closed subschemes inside the
naive models by strengthening the original moduli problem of
\cite{RZ}). Indeed, most investigations of local models so far used
these strengthened moduli problems in a way or another (for a survey
of most progress in this area, we refer to the recent paper
\cite{PRS}).

However, as observed by Pappas and Rapoport in \cite{PR1}, the brute
force definition of the local models by taking the flat closure is
not totally out of control as one might think. Namely, it is known
after G\"{o}rtz' work that the special fibers of the naive models
always embed in the affine flag varieties and that their reduced
subschemes are a union of Schubert varieties. Therefore, it is a
question to describe which Schubert varieties will appear in the
special fibers (of the flat models) and whether the special fibers
are reduced. These questions are reduced to the coherence conjecture
(see \cite{PR1,PR2}, at least in the case the group splits over a
tamely ramified extension), which characterizes the dimension of the
spaces of global sections of certain ample line bundles on certain
union of Schubert varieties. Therefore, we will have a fairly good
understanding of the local models even if we do not know the moduli
problem they represent, provided we can prove the coherence
conjecture.

Let us be a little bit more precise. To this goal, we first need to
recall the theory of affine flag varieties (we refer to \S \ref{loop
grp and flag var} for unexplained notations and more details). Let
$k$ be a field and $\calG$ be a flat affine group scheme of finite
type over $k[[t]]$. Let $G$ be fiber of $\calG$ over the generic
point $F=k\ppart=k[[t]][t^{-1}]$. Then one can define the affine
flag variety $\Fl_\calG=LG/L^+\calG$, which is an ind-scheme, of
ind-finite type (cf. \cite{BL1,F,PR1} and \S \ref{loop grp and flag
var}). When $G$ is an almost simple, simply-connected algebraic
group over $k\ppart$, and $\calG$ is a parahoric group scheme of
$G$, $\Fl_\calG$ is ind-projective and coincides with the affine
flag varieties arising from the theory of affine Kac-Moody groups as
developed in \cite{Ku,Ma} (at least when $G$ splits over a tamely
ramified extension of $k\ppart$). The jet group $L^+\calG$ acts on
$\Fl_\calG$ by left translations and the orbits are finite
dimensional; their closures are called (affine) Schubert varieties.
When $\calG$ is an Iwahori group scheme of $G$, Schubert varieties
are parameterized by elements in the affine Weyl group $W_\aff$ of
$G$ (more generally, if $G$ is not simply-connected, they are
parameterized by elements in the Iwahori-Weyl group
$\widetilde{W}$). For $w\in \widetilde{W}$, we denote the
corresponding Schubert variety by $\Fl_w$.

Let us come back to local models. Let $(G,K,\{\mu\})$ be a triple,
where $G$ is a reductive group over a $p$-adic field $F$, with
finite residue field $k_F$, $K$ is a parahoric subgroup of $G$ and
$\{\mu\}$ is a geometric conjugacy class of one-parameter subgroups
of $G$. Let $E/F$ be the reflex field (i.e. the field of definition
for $\{\mu\}$), with ring of integers $\calO_E$ and residue field
$k_E$. Then for most such triples (at least when $\mu$ is minuscule,
cf. \cite{PRS} for a complete list), one can define the so-called
naive model $\calM^{\on{naive}}_{K,\{\mu\}}$), which is an
$\calO_E$-scheme, whose generic fiber is the flag variety $X(\mu)$
of parabolic subgroups of $G_E$ of type $\mu$. Inside
$\calM^{\on{naive}}_{K,\{\mu\}}$, one defines
$\calM^{\on{loc}}_{K,\{\mu\}}$ as the flat closure of the generic
fiber (for an example of the definitions of such schemes, see \S
\ref{lines on local models}). In all known cases, one can find a
reductive group $G'$ defined over $k\ppart$ and a parahoric group
scheme $\calG$ over $k[[t]]$, such that the special fiber
\[\overline{\calM}^{\on{naive}}_{K,\{\mu\}}:=\calM^{\on{naive}}_{K,\{\mu\}}\otimes k_E\]
embeds into the affine flag variety $\Fl_{\calG}=LG'/L^+\calG$ as a
closed subscheme, which is in addition invariant under the action of
$L^+\calG$. In particular, the reduced subscheme of
$\overline{\calM}^{\on{naive}}_{K,\{\mu\}}$ is a union of Schubert
varieties inside $\Fl_\calG$. Which Schubert variety will appear in
$\overline{\calM}^{\on{naive}}_{K,\{\mu\}}$ usually can be read from
the moduli definition of $\calM^{\on{naive}}_{K,\{\mu\}}$. However,
the special fiber of $\calM^{\on{loc}}_{K,\{\mu\}}$ is more
mysterious, and a lot of work has been done in order to understand
it (we refer to \cite{PRS} (in particular its Section 4) and
references therein for a detailed survey of the current progress).

Here we review two strategies to study
$\calM^{\on{loc}}_{K,\{\mu\}}$. For simplicity, we assume that the
derived group of $G$ is simply-connected and $K$ is an Iwahori
subgroup of $G$ at this moment. In this case, $\calG$ will be an
Iwahori group scheme of $G'$. One can attach to $\{\mu\}$ a subset
$\Adm(\mu)$ in the Iwahori-Weyl group $\widetilde{W}$, usually
called the $\mu$-admissible set (cf. \cite{R} and \S \ref{grp data}
for the definitions). In all known cases, it is not hard to see that
the Schubert varieties $\Fl_w$ for $w\in\Adm(\mu)$ indeed appear in
$\overline{\calM}^{\on{loc}}_{K,\{\mu\}}$, i.e.
\[\calA(\mu):=\bigcup_{w\in\Adm(\mu)}\Fl_w\subset\overline{\calM}^{\on{loc}}_{K,\{\mu\}}.\]
Now, the first strategy to determine the (underlying reduced closed
subscheme of) the special fiber
$\overline{\calM}^{\on{loc}}_{K,\{\mu\}}$ goes as follows. Write
down a moduli functor $\calM'_{K,\{\mu\}}$ which is a closed
subscheme of $\calM^{\on{naive}}_{K,\{\mu\}}$, such that
\[\calM'_{K,\{\mu\}}\otimes E=\calM^{\on{naive}}_{K,\{\mu\}}\otimes E,\quad\quad \overline{\calM}'_{K,\{\mu\}}(\bar{k})=\calA(\mu)(\bar{k}),\]
where $\bar{k}$ is an algebraic closure of $k_E$. Clearly, this will
imply that the reduced subscheme
\begin{equation}\label{reduced loci}
(\overline{\calM}^{\on{loc}}_{K,\{\mu\}})_{\on{red}}=\calA(\mu).
\end{equation}
In fact, much of the previous works about
$\calM^{\on{loc}}_{K,\{\mu\}}$ followed this strategy. However, let
us mention that (so far) the definition of $\calM'_{K,\{\mu\}}$
itself is not group theoretical (i.e. it depends on choosing some
representations of the group $G$). In particular, when $G$ is
ramified, its definition can be complicated. In addition, except a
few cases, it is not known whether
$\calM'_{K,\{\mu\}}=\calM^{\on{loc}}_{K,\{\mu\}}$ in general.

There is another strategy to determine
$\overline{\calM}^{\on{loc}}_{K,\{\mu\}}$, as proposed in
\cite{PR1}. Namely, let us choose an ample line bundle $\calL$ over
$\calM^{\on{navie}}_{K,\{\mu\}}$. Then since by definition
$\calM^{\on{loc}}_{K,\{\mu\}}$ is flat over $\calO_E$ with generic
fiber $X(\mu)$, for $n\gg 0$,
\[\dim_{k_E}\Ga(\calA(\mu),\calL^n)\leq\dim_{k_E}\Ga(\overline{\calM}^{\on{loc}}_{K,\{\mu\}},\calL^n)=\dim_E\Ga(X(\mu),\calL^n).\]
The general expectation (which has been verified in all known cases)
is that
\[\overline{\calM}^{\on{loc}}_{K,\{\mu\}}=\calA(\mu)\]
led Pappas and Rapoport to conjecture the following equivalent
statement
\[\dim_{k_E}\Ga(\calA(\mu),\calL^n)=\dim_E\Ga(X(\mu),\calL^n).\]
Apparently, this conjecture would not be very useful unless one can
say something about the line bundle $\calL$. In fact, the statement
of the conjecture in \cite{PR1} is different and more precise.
Namely, in the \emph{loc. cit.}, they constructed some line bundle
$\calL_1$ on the affine flag variety $\Fl_\calG$ and some line
bundle $\calL_2$ on $X(\mu)$, both of which are explicit and are
given purely in terms of group theory (see \S \ref{the conj} for the
precise construction). Then they conjectured

\medskip

\noindent\bf The Coherence Conjecture. \rm\emph{For $n\gg 0$,}
\[\dim_{\bar{k}}\Ga(\calA(\mu),\calL_1^n)=\dim_E\Ga(X(\mu),\calL_2^n).\]

\medskip

In addition, in \emph{loc. cit.}, for certain groups, they
constructed natural ample line bundles $\calL$ on the corresponding
local models, whose restrictions give $\calL_1$ and $\calL_2$.

\medskip

What makes the coherence conjecture useful? First of all, the
conjecture is group theoretic, i.e. the statement is uniform for all
groups. The non-group theoretic part then is absorbed into the
construction of natural line bundles on local models and the
identification of their restrictions with the group theoretically
constructed line bundles. This is a much simpler problem. An example
is illustrated in \S \ref{lines on local models}. More importantly,
the right hand side in the coherence conjecture is defined over
$\calO_E$ and therefore, it is equivalent to prove that
\[\dim_{\bar{k}}\Ga(\calA(\mu),\calL_1^n)=\dim_{\bar{k}}\Ga(X(\mu),\calL_2^n).\]
Observe that in the above formulation, everything is over the field
$k$ rather than over a mixed characteristic ring. That is, we are
dealing with algebraic geometry rather than arithmetic!

How can we prove this conjecture? Suppose that we can find a scheme
$\bGr_{\calG,\mu}$ (the reason we choose this notation will be clear
soon), which is flat over $\bar{k}[t]$, together with a line bundle
$\calL$ such that its fiber over $0\in\bbA^1$ is
$(\calA(\mu),\calL_1)$ and its fiber over $y\neq 0$ is
$(X(\mu),\calL_2)$, then the coherence conjecture will follow. In
fact, such $\bGr_{\calG,\mu}$ does exist and can be constructed
purely group theoretically. They are the (generalized) equal
characteristic counterparts of local models, which we will call the
global Schubert varieties. Let us briefly indicate the construction
of $\bGr_{\calG,\mu}$ here (the construction of the line bundle
$\calL$, which we ignore here, is also purely group theoretical, see
\S \ref{lines on globaff}). For simplicity, let us assume that $G'$
is split over $k$ (the non-split case will also be considered in the
paper). Let $B$ be a Borel subgroup of $G'$. Then in \cite{G},
Gaitsgory (following ideas of Kottwitz and Beilinson) constructed a
family of ind-schemes $\Gr_{\calG}$ over $\bbA^1$, which is a
deformation from the affine Grassmannian $\Gr_{G'}$ of $G'$ to the
affine flag variety $\Fl_{G'}$ of $G'$. By its construction,
\[\Gr_{\calG}|_{\bbG_m}\cong (\Gr_{G'}\times G'/B)\times\bbG_m\footnote{In the main body of this paper, we will work with a different family so that this extra $G'/B$ factor does not appear.},\quad \Gr_{\calG}|_0\cong\Fl_{G'},\]
where $\Gr_{\calG}|_0$ denotes the scheme theoretic fiber of
$\Gr_{\calG}$ over $0\in \bbA^1$. When $\mu$ is minuscule, the Schubert
variety $\bGr_\mu$ corresponds to $\mu$ in $\Gr_{G'}$ is in fact
isomorphic to $X(\mu)$. In addition, we can ``spread it out" over
$\bbG_m$ as $(\bGr_\mu\times\ast)\times\bbG_m$ to get a closed
subscheme of $\Gr_\calG|_{\bbG_m}$, where $\ast$ is the base point
in $G'/B$. Now define $\bGr_{\calG,\mu}$ as the closure of
$(\bGr_\mu\times\ast)\times\bbG_m$ inside $\Gr_{\calG}$. By
definition, its fiber over $y\neq 0$ is isomorphic to $X(\mu)$. On
the other hand, it is not hard to see that
$\calA(\mu)\subset\bGr_{\calG,\mu}|_0$ (cf. Lemma \ref{easy}).
Therefore, the coherence conjecture will follow if we can show that
$\bGr_{\calG,\mu}|_0=\calA(\mu)$ (and if we can construct the
corresponding line bundle).

At the first sight, it seems the idea is circular. However, it is
not the case. The reason, as we mentioned before, is that
$\bGr_{\calG,\mu}$ now is a scheme over $k$ and we have many more
tools to attack the problem. Observe that to prove that
$\bGr_{\calG,\mu}|_0=\calA(\mu)$, we need to show that
\begin{enumerate}
\item $(\bGr_{\calG,\mu}|_0)_{\on{red}}=\calA(\mu)$ (Theorem \ref{top
fiber});
\item $\bGr_{\calG,\mu}|_0$ is reduced (Theorem \ref{sch fiber}).
\end{enumerate}

Part (1) can be achieved by the calculation of the nearby cycle
$\calZ_\mu=\Psi_{\bGr_{\calG,\mu}}(\bbQ_\ell)$ of the family
$\bGr_{\calG,\mu}$ (see Lemma \ref{support}). Usually, such a
calculation is a hard problem. The miracle here is that if
$\calZ_\mu$ is regarded as an object in the category of Iwahori
equivariant perverse sheaves on $\Fl_{G'}$, it has very nice
properties. Namely, by the main result of \cite{G} (in the case when
$G'$ is split), $\calZ_\mu$ is a central sheaf, i.e. for any other
Iwahori equivariant perverse sheaf $\calF$ on $\Fl_{G'}$, the
convolution product $\calZ_\mu\star\calF$ (see
\eqref{twp1}-\eqref{conv prod} for the definition) is perverse and
\[\calZ_\mu\star\calF\cong\calF\star\calZ_\mu.\]
Then by a result of Arhkipov-Bezrukavnikov \cite[Theorem 4]{AB}, the
above properties put a strong restriction of the support of
$\calZ_\mu$, which will imply Part (1). We shall mention that
although we assume here that $G'$ is split, the same strategy can be
applied to the non-split groups. This is done in \S \ref{PfII},
where we generalize the results of \cite{G} and \cite{AB} to
ramified groups as well. Our arguments are simpler than the
originally arguments in \cite{G,AB}, and will provide the following
technical advantage. As we mentioned above, $\bGr_{\calG,\mu}$
should be regarded as the equal characteristic counterparts of local
models. Therefore, it is natural (and indeed important) to determine
the nearby cycles $\Psi_{\calM^{\on{loc}}_{K,\{\mu\}}}(\bbQ_\ell)$
for the local models. For example, if one could prove that these
sheaves are also central (the \emph{Kottwitz conjecture}\footnote{In
fact, the Kottwitz conjecture is weaker than this statement, and its
significance lies in its use in the Langlands-Kottwitz method for
calculating the Zeta functions of Shimura varieties.}), then one
could conclude \eqref{reduced loci} directly. It turns out the
arguments in \S \ref{PfII} have a direct generalization to the mixed
characteristic situation and in joint work with Pappas \cite{PZ}, we
use it to show the Kottwitz conjecture (some previous cases are
proved by Haines and Ng\^{o} \cite{HN}).

Now we turn to Part (2), which is more difficult. The idea is that
we can assume $\on{char} k>0$ and use the powerful technique of
Frobenius splitting (cf. \cite{MR,BK}). To prove that
$\bGr_{\calG,\mu}|_0$ is reduced, it is enough to prove that it is
Frobenius split. To achieve this goal, we embeds $\bGr_{\calG,\mu}$
into a larger scheme $\bGr_{\calG,\mu,\la}^{BD}$ over $\bbA^1$,
which is a closed subscheme of a version of the Beilinson-Drinfeld
Grassmannian. The scheme $\bGr_{\calG,\mu,\la}^{BD}$ is normal and
its fiber over $0$ is reduced. Then to prove that
\[\bGr_{\calG,\mu}|_0=\bGr_{\calG,\mu}\cap\bGr_{\calG,\mu,\la}^{BD}|_0\]
is Frobenius split, it is enough to construct a Frobenius splitting
of $\bGr_{\calG,\mu,\la}^{BD}$, compatible with $\bGr_{\calG,\mu}$
and $\bGr_{\calG,\mu,\la}^{BD}|_0$. Since
$\bGr_{\calG,\mu,\la}^{BD}$ is normal, it is enough to prove this
for some nice open subscheme $U\subset \bGr_{\calG,\mu,\la}^{BD}$,
such that $\bGr_{\calG,\mu,\la}^{BD}-U$ has codimension two. In
particular, the open subscheme $U$ will not intersect with
$\bGr_{\calG,\mu}|_0$, which is our primary interest. Section
\ref{PfI} is devoted to realizing this idea.

\medskip

Now let us describe the organization of the paper and some other
results proved in it.

In \S \ref{loc}, we review the coherence conjecture of Pappas and
Rapoport. In \S \ref{grp data}, we review the basic theory of
reductive groups over local fields and introduce various notations
used in the rest of the paper. In \S \ref{loop grp and flag var}, we
rapidly recall the main results of \cite{PR1} (and \cite{F})
concerning loop groups and the geometry of their flag varieties. In
\S \ref{the conj}, we state the main theorems (Theorem \ref{MainI}
and \ref{MainII}) of our paper, which give a modified version of
original coherence conjecture of Pappas and Rapoport (see Remark
\ref{correction} for the reason of the modification).

In \S \ref{globS}, we introduce the main geometric object we are
going to study in the paper, namely, the global Schubert varieties.
They are varieties projective over the affine line $\bbA^1$, which
are the counterparts of local models in the equal characteristic
situation. In \S \ref{globaffGrass}, we define the global affine
Grassmannian over a curve for general (non-constant) group schemes.
After the work of \cite{PR1,PR3,H}, this construction is now
standard. In \S \ref{grpsch}, we construct a special Bruhat-Tits
group scheme over $C=\bbA^1$, i.e. a group scheme which is only
ramified at the origin. Let us remark that similar constructions are
also considered in \cite{HNY,Ri}. In \S \ref{glob Sch}, we apply the
construction of the global affine Grassmannian to the group scheme
we consider in the paper. We introduce the global Schubert variety
$\bGr_{\calG,\mu}$, which is associated to a geometric conjugacy
class of 1-parameter subgroup $\{\mu\}$ of $G$, over a ramified
cover $\tilde{C}$ of $C$. We then state another main theorem
(Theorem \ref{Fibers}) which asserts that the special fiber of
$\bGr_{\calG,\mu}$ is $\calA^Y(\mu)$, and first show that the
variety $\calA^Y(\mu)$ is contained in this special fiber (Lemma
\ref{easy}). In \S \ref{lines on globaff}, we explain that our
assertion about the special fiber of $\bGr_{\calG,\mu}$ is
equivalent to the coherence conjecture. The key ingredient is a
certain line bundle on the global affine Grassmannian, namely, the
pullback of the determinant line bundle along the closed embedding
\[\Gr_\calG\to\Gr_{\Lie(\calG)}.\]
We calculate its central charges at each fiber (which turn out to be
twice of the dual Coxeter number) and find the remarkable fact that
the central charge of line bundles on the global affine
Grassmannians are constant along the curve (Proposition \ref{cal of
cc}).

In \S \ref{some prop}, we make some preparations towards the proof
of our main theorem. We study two basic geometrical structures of
$\bGr_{\calG,\mu}$: (i) in \S \ref{affine chart}, we will construct
certain affine charts of $\bGr_{\calG,\mu}$, which turn out to be
isomorphic to affine spaces over $\tilde{C}$; and (ii) in \S \ref{Gm
action}, we will construct a $\bbG_m$-action on $\bGr_{\calG,\mu}$,
so that the map $\bGr_{\calG,\mu}\to \tilde{C}$ is
$\bbG_m$-equivariant, where $\bbG_m$ acts on $\tilde{C}=\bbA^1$ by
natural dilatation. To establish (i), we will need to first
construct the global root subgroups of $\calL\calG$ as in \S
\ref{globrootgrp}.

The next two sections are then devoted to the proof of the theorem
concerning the special fiber of $\bGr_{\calG,\mu}$, as has been
already outlined above. The first part of the proof, presented in \S
\ref{PfI}, concerns the scheme theoretic structure of the special
fiber. Namely, we prove that it is reduced. This is achieved by the
technique of Frobenius splitting. As a warm up, we prove in \S
\ref{factor} that Theorem \ref{MainI} is a special case of Theorem
\ref{MainII}, which should be well-known to experts. Then we
introduce the Beilinson-Drinfeld Grassmannian and the convolution
Grassmannian and reduce Theorem \ref{sch fiber} to Theorem
\ref{Fsplit}. In \S \ref{sp}, we prove a special case of Theorem
\ref{Fsplit} by studying the affine flag variety associated to a
special parahoric group scheme. Recall that a result of
Beilinson-Drinfeld (cf. \cite[4.6]{BD}) asserts that the Schubert
varieties in the affine Grassmannian are Gorenstein. We examine in
\S \ref{sp} to what extend this result holds for ramified groups
(i.e. reductive groups split over a ramified extension). It turns
out this result extends to all affine flag varieties associated to
special parahorics except in the case the special parahoric is a
parahoric of the ramified odd unitary group $\SU_{2n+1}$, whose
special fiber has reductive quotient $\on{SO}_{2n+1}$ (Theorem
\ref{Goren}). In this exceptional case, no Schubert variety of
positive dimension in the corresponding affine flag variety is
Gorenstein (Remark \ref{non-Goren}).  In \S \ref{PfII}, we give the
second part of the proof, which asserts that topologically, the
special fiber of $\bGr_{\calG,\mu}$ coincides with $\calA^Y(\mu)$.
This is achieved by the description of the support of the nearby
cycle (for the intersection cohomology sheaf) of this family. In the
case when the group is split, this follows the earlier works of
\cite{G} and \cite{AB}. In \S \ref{central sheaves} and \S \ref{Waki fil}, we generalize
their results to ramified groups, with certain simplifications of
the original arguments.

The paper has two appendices. The first one, \S \ref{lines on local
models}, calculates the line bundles on the local models for the
ramified unitary groups. The study of these local models are the
main motivation for Pappas and Rapoport to make the coherence
conjecture. Since their original conjecture is not as stated in our
main theorem, we explain in this appendix why our main theorem is
correct for the applications to local models. The second appendix
(\S \ref{recol}) collects and strengthens some results, which already exist in literature, in a form
needed in the main body of the paper.

\medskip

\noindent\bf Notations. \rm Let $k$ be a field, and fix $\bar{k}$ to
be an algebraic closure of $k$. We will denote by
$k^s\subset\bar{k}$ the separable closure of $k$ in $\bar{k}$.

If $X$ be a $Y$-scheme and $V\to Y$ is a morphism, the base change
$X\times_YV$ is denoted by $X_V$ or $X|_V$. If $V=\Spec R$, it is
sometimes also denoted by $X_R$. If $V=x=\Spec k$ is a point, then
it is sometimes also denoted by $(X)_x$.

For a vector bundle $\calV$ on a scheme $V$, we denote by
$\det(\calV)$ the top exterior power of $\calV$, which is a line
bundle.

If $A$ is an affine algebraic group (not necessarily a torus) over a
field $k$, we denote by $\xch(A)$ (resp. $\xcoch(A)$) its character
group (resp. cocharacter group) over $k^s$. The Galois group
$\Ga=\Gal(k^s/k)$ acts on $\xch(A)$ (resp. $\xcoch(A)$) and the
invariants (resp. coinvariants) are denoted by $\xch(A)^\Ga$ (resp.
$\xch(A)_\Ga, \xcoch(A)^\Ga, \xcoch(A)_\Ga$).

If $\calG$ is a flat group scheme over $V$, the trivial
$\calG$-torsor (i.e. $\calG$ itself regarded as a $\calG$-torsor by
right multiplication) is denoted by $\calE^0$. For a $\calG$-torsor
$\calE$, we use $\ad\ \calE$ to denote the associated adjoint
bundle. If $\calP$ is a $\calG$ torsor and $X$ is a scheme over $V$
with an action of $\calG$, we denote the twisted product by
$\calP\times^\calG X$, which is the quotient of $\calP\times_VX$ by
the diagonal action $\calG$.

If $G$ is a reductive group over a field, we denote by $G_\der$ its
derived group, $G_\s$ the simply-connected cover of $G_\der$ and
$G_\ad$ its the adjoint group.

\medskip

\noindent\bf Acknowledgement. \rm The author thanks D. Gaitsgory, G.
Pappas, M. Rapoport, J.-K. Yu for useful discussions, and G. Pappas
and M. Rapoport for reading an early draft of this paper. The work
of the author is supported by the NSF grant DMS-1001280.

\section{Review of the local picture, formulation of the
conjecture}\label{loc}
\subsection{Group theoretical data}\label{grp data}
Let $k$ be an algebraically closed field. Let $\calO=k[[t]]$ and
$F=k\ppart$. Let $\Ga=\Gal(F^s/F)$ be the inertial (Galois) group,
where $F^s$ is the separable closure of $F$. Let us emphasize that
we choose a uniformizer $t$. Let $G$ be a connected reductive group
over $F$. In this paper, unless otherwise stated,  $G$ is assumed to split over a \emph{tamely} ramified extension
$\tilde{F}/F$. It is called a ramified group if it is non-split over $F$.

Let $S$ be a maximal $F$-split torus of $G$. Let $T=\calZ_G(S)$ be
the centralizer of $S$ in $G$, which is a maximal torus of $G$ since
$G$ is quasi-split over $F$ (as $F$ is a field of cohomological dimension one, this follows from \cite[Theorem 1.9]{St}). Let us choose a rational Borel subgroup
$B\supset T$. Let $H$ be a split Chevalley group over $\bbZ$ such
that $H\otimes F^s\cong G\otimes F^s$. We need to choose this
isomorphism carefully. Let us fix a pinning $(H,B_H,T_H,X)$ of $H$
over $\bbZ$. Let us recall that this means that $B_H$ is a Borel
subgroup of $H$, $T_H$ is a split maximal torus contained in $B_H$,
and $X=\Sigma_{\tilde{a}\in\Delta} X_{\tilde{a}}\in\Lie B$, where
$\Delta$ is the corresponding set of simple roots,
$\tilde{U}_{\tilde{a}}$ is the root subgroup corresponding to
$\tilde{a}$ and $X_{\tilde{a}}$ is a generator in the rank one free
$\bbZ$-module $\Lie \tilde{U}_{\tilde{a}}$. Let us choose an
isomorphism $(G,B,T)\otimes_F\tilde{F}\cong
(H,B_H,T_H)\otimes_{\bbZ}\tilde{F}$, where $\tilde{F}/F$ is a cyclic
extension such that $G\otimes\tilde{F}$ splits. This induces an
isomorphism of the root data
$(\xch(T_H),\Delta,\xcoch(T_H),\Delta^\vee)\cong(\xch(T),\Delta,\xcoch(T),\Delta^\vee)$.
Let $\Xi$ be the group of pinned automorphisms of $(H,B_H,T_H,X)$. The natural map from $\Xi$ to the group of the automorphisms of
the root datum $(\xch(T_H),\Delta,\xcoch(T_H),\Delta^\vee)$ is an isomorphism (\cite[Proposition 7.1.6]{C}).

Now the action of $\Gamma=\Gal(\tilde{F}/F)$ on
$G\otimes_F\tilde{F}$ induces a homomorphism $\psi:\Gamma\to\Xi$.
Then we can always choose an isomorphism
\begin{equation}\label{pinned isom}
(G,B,T)\otimes_{F}\tilde{F}\cong(H,B_H,T_H)\otimes_\bbZ\tilde{F}
\end{equation}
such that the action of $\ga\in \Gamma$ on the left hand side
corresponds to $\psi(\ga)\otimes\ga$. In the rest of the paper, we
fix such an isomorphism. This determines a point $v_0$ in $A(G,S)$,
the apartment associated to $(G,S)$ (\cite{BT})\footnote{More
precisely, $v_0$ is a point in the apartment associated to the
adjoint group $(G_\ad,S_\ad)$. But since in the paper, we only use
the combinatorial structures of $A(G,S)$, we will not distinguish it
from the one associated to the adjoint group.}. This is a special
point of $A(G,S)$, which in turn gives a parahoric group scheme
$\calG_{v_0}$ over $\calO$, namely
\begin{equation}\label{constr of special}
\calG_{v_0}:=((\Res_{\calO_{\tilde{F}}/\calO}(H\otimes\calO_{\tilde{F}}))^{\Ga})^0.
\end{equation}
Let us explain the notations. Here $\Res$ stands for the Weil
restriction, so that
$\Res_{\calO_{\tilde{F}}/\calO}(H\otimes\calO_{\tilde{F}})$ is a
smooth group scheme over $\calO$ (cf. \cite[2.2]{Ed}), with an
action of $\Ga$. The notation $(-)^\Ga$ stands for taking the
$\Ga$-fixed point subscheme. Under our tameness assumption,
$\tilde{\calG}_{v_0}:=(\Res_{\calO_{\tilde{F}}/\calO}(H\otimes\calO_{\tilde{F}}))^{\Ga}$
is smooth by \cite[3.4]{Ed}. Finally, $(-)^0$ stands for taking the
neutral connected component. Thereofre, $\calG_{v_0}$ and
$\tilde{\calG}_{v_0}$ have the same generic fiber and the special
fiber of $\calG_{v_0}$ is the neutral connected component of the
special fiber of $\tilde{\calG}_{v_0}$.

Recall that $A(G,S)$ is an affine space under $\xcoch(S)_\bbR$. For
every facet $\sigma\subset A(G,S)$, let $\calG_\sigma$ be the
parahoric group scheme over $\calO$ (in particular, the special
fiber of $\calG_\sigma$ is connected). Let us choose a special
vertex $v\in A(G,S)$ (e.g $v_0$), and identify $A(G,S)$ with
$\xcoch(S)_\bbR$ via this choice. Let $\mathbf{a}$ be the unique alcove in
$A(G,S)$, whose closure contains the point $v$, and is contained in
the finite Weyl chamber determined by our chosen Borel. This
determines a set of simple affine roots $\al_i, i\in\bold S$, where
$\bold S$ is the set of vertices of the affine Dynkin diagram
associated to $G$.

Let $\widetilde{W}$ be the Iwahori-Weyl group of $G$ (cf.
\cite{HR}), which acts on $A(G,S)$. This is defined to be
$\calN_G(S)(F)/\ker\kappa$, where $\calN_G(S)$ is the normalizer of
$S$ in $G$, and
\begin{equation}\label{kot}
\kappa:T(F)\to\xcoch(T)_{\Ga}
\end{equation}
is the Kottwitz homomorphism (cf. \cite[\S 7]{Ko}). One has the
following exact sequence
\begin{equation}\label{affine Weyl}
1\to \xcoch(T)_{\Ga}\to\widetilde{W}\to {W_0}\to 1,
\end{equation}
where ${W_0}$ is the relative Weyl group of $G$ over $F$. In what
follows, we use $t_\la$ to denote the translation element in
$\widetilde{W}$ given by $\la\in\xcoch(T)_\Ga$ from the above map
\eqref{affine Weyl}\footnote{\label{fn}Note that under the sign
convention of the Kottwitz homomorphism in \cite{Ko}, $t_\la$ acts
on $A(G,S)$ by $v\mapsto v-\la$.}. But occasionally, we also use
$\la$ itself to denote this translation element if no confusion is
likely to arise. The pinned isomorphism \eqref{pinned isom}
determines a set of positive roots $\Phi^+=\Phi(G,S)^+$ for $G$.
There is a natural map $\xcoch(T)_\Ga\to\xcoch(S)_\bbR$. We define
\begin{equation}\label{plus}
\xcoch(T)_\Ga^+=\{\la\mid (\la,a)\geq 0 \mbox{ for }
a\in\Phi^+\}.
\end{equation}
Our choice of the special vertex $v$ of $A(G,S)$ gives a splitting
of the exact sequence and, therefore we can write $w=t_\la w_f$ for
$\la\in \xcoch(T)_\Ga$ and $w_f\in {W_0}$.

Let $W_\aff$ be the affine Weyl group of $G$, i.e. the Iwahori-Weyl
group of $G_\s$, which is a Coxeter group. One has
\[1\to \xcoch(T_{\s})_{\Ga}\to\ W_\aff\to {W_0}\to 1,\]
where $T_\s$ is the inverse image of $T$ in $G_\s$. One can write
$\widetilde{W}=W_\aff\rtimes \Omega$, where $\Omega$ is the subgroup
of $\widetilde{W}$ that fixes the chosen alcove $\mathbf{a}$. This gives
$\widetilde{W}$ a quasi-Coxeter group structure. Hence it makes
sense to talk about the length of an element $w\in\widetilde{W}$ and
there is a Bruhat order on $\widetilde{W}$. Namely, if we write
$w_1=w'_1\tau_1, w_2=w'_2\tau_2$ with $w'_i\in W_\aff,
\tau_i\in\Omega$, then $\ell(w_i)=\ell(w'_i)$ and $w_1\leq w_2$ if
and only if $\tau_1=\tau_2$ and $w'_1\leq w'_2$. A lot of the
combinatorics of the Iwahori-Weyl group arises from the study of the
restriction of the length function and the Bruhat order to
$\xcoch(T)_\Ga\subset\widetilde{W}$. Some of them will be reviewed
in \S \ref{cb of IW}.

Now let us recall the definition of the \emph{admissible set} in the
Iwahori-Weyl group. Let $\bar{W}$ be the absolute Weyl group of $G$,
i.e. the Weyl group for $(H,T_H)$. Suppose that
$\mu:(\bbG_m)_{\tilde{F}}\to G\otimes \tilde{F}$ gives a geometric
conjugacy class of 1-parameter subgroups. It determines a
$\bar{W}$-orbit in $\xcoch(T)$. One can associate $\{\mu\}$ a
${W_0}$-orbits $\Lambda$ in $\xcoch(T)_{\Ga}$ as follows. Choose a
Borel subgroup of $G$ containing $T$, and is defined over $F$. This
gives a unique element in this $\bar{W}$-orbit, still denoted by
$\mu$, which is dominant w.r.t. this Borel subgroup. Let $\bar{\mu}$
be its image in $\xcoch(T)_{\Ga}$, and let $\Lambda={W_0}\bar{\mu}$.
It turns out $\Lambda$ does not depend on the choice of the rational
Borel subgroup of $G$, since any two such $F$-rational Borels that
contain $T$ will be conjugate to each other by an element in
${W_0}$. For $\mu\in\xcoch(T)$, define the admissible set
\begin{equation}\label{Adm}
\Adm(\mu)=\{w\in\widetilde{W}\mid w\leq t_\la, \mbox{ for some }
\la\in\Lambda\}.
\end{equation}

Under the map $\xcoch(T)_{\Ga}\to\widetilde{W}\to
\widetilde{W}/W_\aff\cong\Omega$, the set $\Lambda$ maps to a single
element (cf. \cite[Lemma 3.1]{R}), denoted by $\tau_\mu$. Define
\[\Adm(\mu)^\circ=\tau_\mu^{-1}\Adm(\mu).\]
For $Y\subset \bold S$ any subset, let $W^Y$ denote the subgroup of
$W_\aff$ generated by $\{r_{i}, i\in\bold S-Y\}$, where $r_i$ is the
simple reflection corresponding to $i$. Then set
\[\Adm^Y(\mu)= W^Y\Adm(\mu)W^Y\subset \widetilde{W},\]
and
\[\Adm^Y(\mu)^\circ=\tau_\mu^{-1}\Adm^Y(\mu).\]
Note that $\Adm^Y(\mu)^\circ\subset W_\aff$, and this subset only depends on the image of $\mu$ under $\xcoch(T)\to\xcoch(T_\ad)$, where
$T_\ad$ is the image of $T$ in $G_\ad$. 

\subsection{Loop groups and their flag varieties}\label{loop grp and flag var}
Let $\sigma\subset A(G,S)$ be a facet. Let
\[\Fl_\sigma=LG/L^+\calG_{\sigma}\]
be the (partial) flag variety of $LG$. Let us recall that $LG$ is
the loop group of $G$, which represents the functor which associates
to every $k$-algebra $R$ the group $G(R\ppart)$, $L^+\calG_{\sigma}$
is the jet group of $\calG_{\sigma}$, which represents the functor
which associates to every $k$-algebra $R$ the group
$\calG_{\sigma}(R[[t]])$, and $\Fl_\sigma=LG/L^+\calG_{\sigma}$ is
the fpqc quotient. Let us also recall that $LG$ is represented by an
ind-affine scheme, $L^+\calG_{\sigma}$ is represented by an affine
scheme, which is a closed subscheme of $LG$, and $\Fl_\sigma$ is
represented by an ind-scheme, ind-projective over $k$. Denote by
$I=L^+\calG_{\mathbf{a}}$ the Iwahori subgroup of $LG$, and denote $\Fl_{\mathbf{a}}$ by
$\Fl$, which we call the affine flag variety of $G$. If $G$ splits
over $F$, so that $G=H\otimes F$ and \eqref{pinned isom} is the
natural isomorphism, then the special vertex $v_0$ is hyperspecial,
and corresponds to the parahoric group scheme $H\otimes k[[t]]$.
Then we denote $\Fl_{v_0}$ by $\Gr_H$ and call it the affine
Grassmannian of $H$. Let $Y\subset \bold S$ be a subset, and
$\sigma_Y\subset A(G,S)$ be the facet such that $\al_i(\sigma_Y)=0$
for $i\in\bold S-Y$. Observe that $\sigma_{\bold S}=\mathbf{a}$ is the chosen
alcove. We also denote $\Fl_{\sigma_Y}$ by $\Fl^Y$ for simplicity.

Let us recall that the $I$-orbits of $\Fl$ are parameterized by
$\widetilde{W}$. In general, the $L^+\calG_{\sigma_Y}$-orbits of
$\Fl^{Y'}$ are parameterized by $W^Y\setminus\widetilde{W}/W^{Y'}$,
where $W^Y$ is the Weyl group of $\calG_{\sigma_Y}\otimes k$. For
$w\in\widetilde{W}$, let ${^Y\Fl^{Y'}_w}\subset \Fl^{Y'}$ denote the
corresponding Schubert variety, i.e. the closure of the
$L^+\calG_{\sigma_Y}$-orbit through $w$. If $Y=Y'$, then we simply
denote it by $\Fl^Y_w$. If $G$ is split, and $\calG=H\otimes k[[t]]$
is a hyperspecial model, recall that $L^+\calG$-orbits of $\Gr_H$
are parameterized by
$\bar{W}\setminus\widetilde{W}/\bar{W}\cong\xcoch(T)^{+}$, the set
of dominant coweights of $G$. For $\mu\in\xcoch(T)^{+}$, let
$\bGr_\mu$ be the corresponding Schubert variety in $\Gr_H$.

\quash{More precisely, let
\begin{equation}\label{tla}
t_\la:\Spec K\to \bbG_m\stackrel{\la}{\to}T\subset H,
\end{equation}
then $\bGr_\la$ is the closure of $L^+H$-orbit through $(t_\la \mod
L^+H)$, which is the same as the closure of the $B$-orbit, since we
assume that $\la$ is dominant.}

Let us recall the following result of \cite{F,PR1}.
\begin{thm}\label{NSV}
Let $p=\on{char}k$. Assume that $p\nmid |\pi_1(G_\der)|$, where
$G_\der$ is the derived group of $G$. Then the Schubert variety
$\Fl^Y_w$ is normal, has rational singularities, and is
Frobenius-split if $p>0$.
\end{thm}

For $\mu\in\xcoch(T)$, let
\begin{equation}\label{nAmu}
\calA^Y(\mu)^\circ=\bigcup_{w\in\Adm^Y(\mu)^\circ}{^{Y^\circ}\Fl^{Y}_{\s,w}},
\end{equation}
where $\sigma_{Y^\circ}=\tau_\mu^{-1}(\sigma_Y)$, and where
${^{Y^\circ}\Fl^{Y}_{\s,w}}$ is the union of Schubert varieties
(more precisely, the closure of
$L^+\calG_{\sigma_{Y^\circ}}$-orbits) in the partial affine flag
variety $\Fl^{Y}_\s=LG_\s/L^+\calG_{\sigma_{Y}}$. Then
$\calA^Y(\mu)^\circ$ is a reducible subvariety of $\Fl^{Y}_\s$, with
irreducible components
\[{^{Y^\circ}\Fl^{Y}_{\s,\tau_\mu^{-1}t_\la}}, \quad \la\in\Lambda\subset\xcoch(T)_{\Ga}\subset\widetilde{W}.\]
Observe that $\calA^Y(\mu)^\circ$ only depends on the image of $\mu$  under $\xcoch(T)\to\xcoch(T_\ad)$.

When $p\nmid|\pi_1(G_\der)|$, it is also convenient to consider
\begin{equation}\label{Amu}
\calA^Y(\mu)=\bigcup_{w\in\Adm^Y(\mu)}\Fl^Y_w.
\end{equation}
Choosing a lift $g\in G(F)$ of $\tau_\mu\in\widetilde{W}$ and
identifying $\Fl_{\s}^Y$ with the reduced part of the neutral
connected component of $\Fl^Y$ (see \cite[\S 6]{PR1}), we can define
a map $\Fl_\s^Y\to\Fl^Y,\ x\mapsto gx$. Clearly, this map induces an
isomorphism $\calA^Y(\mu)^\circ\cong\calA^{Y}(\mu)$.

In particular, if $G=H\otimes F$ is split and $\sigma_Y=v_0$ is the
hyperspecial vertex corresponding to $H\otimes\calO$, then
$\calA^Y(\mu)^\circ$ is denoted by $\Gr_{\leq\mu}$, so that if
$p\nmid|\pi_1(G_\der)|$, then we have the isomorphism
$\Gr_{\leq\mu}\cong\bGr_\mu$.

\medskip

We also need to review the Picard group of $\Fl$. For simplicity, we
assume that $G$ is simple, simply-connected, absolutely simple. In
this case $\Fl$ is connected. For each $i\in\bold S$, let $P_i$ be
the corresponding parahoric subgroup scheme such that $L^+P_i\supset I$ so that
$L^+P_i/I\cong\bbP^1$. This $\bbP^1$ maps naturally to $\Fl$ via
$L^+P_i\to LG$, and the image will be denoted as $\bbP^1_i$. Then it is
known (\cite[\S 10]{PR1}) that there is a unique line bundle
$\calL(\epsilon_i)$ on $\Fl$, whose restriction to the $\bbP^1_i$ is
$\calO_{\bbP^1}(1)$, and whose restrictions to other $\bbP^1_j$s
with $j\neq i$ is trivial. Then there is an isomorphism
\[\Pic(\Fl)\cong\bigoplus_{i\in\bold S} \bbZ \calL(\epsilon_i).\]
Let us write $\otimes_i\calL(\epsilon_i)^{n_i}$ as $\calL(\sum_i
n_i\epsilon_i)$. As explained in \emph{loc. cit.}, the $\epsilon_i$
can be thought of as the fundamental weights of the Kac-Moody group
associated to $LG$, and therefore, $\Pic(\Fl)$ is identified with
the weight lattice of the corresponding Kac-Moody group.

There is also a morphism
\begin{equation}\label{central charge}
c:\Pic(\Fl)\to\bbZ
\end{equation}
called the central charge. If we identify $\calL\in\Pic(\Fl)$ with a
weight of the corresponding Kac-Moody group, then $c(\calL)$ is just
the restriction of this weight to the central $\bbG_m$ in the
Kac-Moody group. Explicitly,
\begin{equation}\label{cc}c(\calL(\epsilon_i))=a_i^\vee,\end{equation} where $a_i^\vee
(i\in\bold S)$ are defined as in \cite[6.1]{Kac}. The kernel of
$c$ can be described as follows. Let $s$ denote the closed point of
$\Spec \calO$, and let $(\calG_{\mathbf{a}})_s$ denote the special fiber of
$\calG_{\mathbf{a}}$. Recall that for any $k$-algebra $R$, $\Fl(R)$ is the set
of $\calG_{\mathbf{a}}$-torsors on $\Spec R[[t]]$ together with a
trivialization over $\Spec R\ppart$. Therefore, by restriction of
the $\calG_{\mathbf{a}}$-torsors by $t\mapsto 0$ to $\Spec R\subset\Spec
R[[t]]$, we obtain a natural morphism $\Fl\to\bbB (\calG_{\mathbf{a}})_s$ (here
$\bbB (\calG_{\mathbf{a}})_s$ is the classifying stack of $(\calG_{\mathbf{a}})_s$), which
induces $\xch((\calG_{\mathbf{a}})_s)\cong\Pic(\bbB(\calG_{\mathbf{a}})_s)\to\Pic(\Fl)$.
We have the short exact sequence
\begin{equation}\label{ker of c}0\to\xch((\calG_{\mathbf{a}})_s)\to\Pic(\Fl)\stackrel{c}{\to}\bbZ\to 0.\end{equation}

Now let $Y\subset\bold S$ be a non-empty subset. Observe that if
$\calL(\sum n_i\epsilon_i)$ is a line bundle on $\Fl$, with $n_i=0$
for $i\in\bold S-Y$, then this line bundle is the pullback of a
unique line bundle along $\Fl\to\Fl^Y$, denoted by
$\calL^Y(\sum_{i\in Y} n_i\epsilon_i)$. In this way, we have
\begin{equation}\label{rev:pic}
\Pic(\Fl^Y)\cong\bigoplus_{i\in Y}\bbZ\calL(\epsilon_i).
\end{equation}
The central charge of a line bundle $\calL$ on $\Fl^Y$ is defined to
be the central charge of its pullback to $\Fl$, i.e. the image of
$\calL$ under $\Pic(\Fl^Y)\to\Pic(\Fl)\stackrel{c}{\to}\bbZ$.
Observe that $\calL^Y(\sum_{i\in Y}n_i\epsilon_i)$ is ample on
$\Fl^Y$ if and only if $n_i>0$ for all $i\in Y$.

In the case $G=H\otimes F$ is split, the central charge map induces
an isomorphism $c:\Pic(\Gr_H)\cong\bbZ$. We will denote by
$\calL_{b}$ the ample generator of the Picard group of $\Gr_H$.
Observe that, for $Y=\{i\}$ not special, the ample generator of
$\Pic(\Fl^Y)$ has central charge $a_i^\vee$, which is in general
greater than one. That is, the composition
$\Pic(\Fl^Y)\to\Pic(\Fl)\stackrel{c}{\to}\bbZ$ is injective but not
surjective in general.

\subsection{The coherence conjecture}\label{the conj}

Now we formulate the coherence conjecture of Pappas and Rapoport.
However, the original conjecture, as stated in \emph{loc. cit.}
needs to be modified (see Remark \ref{correction}).

Assume that $G$ is simple, absolutely simple, simply-connected and
splits over a tamely ramified extension $\tilde{F}/F$. Let $\{\mu\}$
be a geometric conjugacy class of 1-parameter subgroups
$(\bbG_m)_{\tilde{F}}\to G_{\ad}\otimes{\tilde{F}}$. First assume
that $\mu$ is minuscule. Let $P(\mu)$ be the corresponding maximal
parabolic subgroup of $H$, and let $X(\mu)=H/P(\mu)$ be the
corresponding partial flag variety of $H$. Let $\calL(\mu)$ be the
ample generator of the Picard group of $X(\mu)$. Then define
\[h_{\mu}(a)=\dim H^0(X(\mu),\calL(\mu)^a).\] If
$\mu=\mu_1+\cdots+\mu_n$ is a sum of minuscule coweights, let
$h_\mu=h_{\mu_1}\cdots h_{\mu_n}$. The following is the main theorem
of this paper, which is a modified version of the original coherence
conjecture of Pappas and Rapoport in \cite{PR1}.

\begin{Thm}\label{MainI}\emph{Let $\mu=\mu_1+\cdots+\mu_n$ be a sum of minuscule coweights,
then for any $Y\subset \bold S$, and ample line bundle $\calL$ on
$\Fl^Y$, we have
\[\dim H^0(\calA^Y(\mu)^\circ,\calL^a)=h_\mu(c(\calL)a),\]
where $c(\calL)$ is the central charge of $\calL$.}
\end{Thm}

This theorem is a consequence of the following more general theorem.

\begin{Thm}\label{MainII}\emph{Let $\mu\in \xcoch(T_\ad)$. Then for any $Y\subset\bold S$,
and ample line bundle $\calL$ on $\Fl^Y$, we have
\[\dim H^0(\calA^Y(\mu)^\circ,\calL)=\dim
H^0(\Gr_{\leq\mu},\calL_{b}^{c(\calL)}).\]}
\end{Thm}

Since Theorem \ref{MainI} is not the same as what Pappas and
Rapoport originally conjectured, and their conjecture is aimed at
studying the local models of Shimura varieties, we will explain why
this is the correct theorem for applications to local models in \S
\ref{lines on local models}. Let us remark that if $G$ is split of
type $A$ or $C$, Theorem \ref{MainI} is proved in \cite{PR1}, using
the previous results on the local models of Shimura varieties (cf.
\cite{Go1,Go2,PR0}). However, it seems that Theorem \ref{MainII} is
new even for symplectic groups.

One consequence of our main theorem (see \S \ref{lines on local models}) is that

\begin{cor}The statement of Theorem 0.1 in \cite{PR2} holds unconditionally.
\end{cor}

Our main theorem can be also applied to local models of other types
(for example for the (even) orthogonal groups) to deduce some
geometrical properties of the special fibers. This will be done in
\cite{PZ}.

\begin{rmk}\label{correction}The original coherence conjecture in
\cite{PR1} needs to be modified. This is due to a miscalculation in
\cite[10.a.1]{PR1}. Namely, when $G$ is simply-connected, the affine
flag variety of $G$ (denoted by $\calF_G$ temporarily) embeds into
the affine flag variety of $H$ (denoted by $\calF_H$ temporarily).
Therefore there is a restriction map
$\Pic(\calF_H)\to\Pic(\calF_G)$, which was described explicitly in
\emph{loc. cit.}. This description is wrong in the case when the group is a non-split even unitary group (and in some other cases). Adjusting the work in \emph{loc. cit.} to account for this produces the modified coherence conjecture that we show in this paper. Let us
remark that the same miscalculation led an incorrect example in
\cite{H} Remark 19 (4), and an incorrect statement in the last
sentence of the first paragraph in p. 502 of \emph{loc. cit.} (see
Proposition \ref{cal of cc}).
\end{rmk}

\section{The global Schubert varieties}\label{globS}
Theorem \ref{MainII} will be a consequence of the geometry of the
global Schubert varieties, which will be introduced in what follows.
Global Schubert varieties are the function field counterparts of the
local models.

\subsection{The global affine Grassmannian}\label{globaffGrass}

Let $C$ be a smooth curve over $k$, and $\calG$ be a smooth \emph{affine} group
scheme over $C$.  Let $\Gr_\calG$ be the global affine Grassmannian
over $C$. Let us recall the functor it represents. For every
$k$-algebra $R$,
\begin{equation}
\Gr_\calG(R)=\left\{(y,\calE,\beta)\ \left|\ \begin{split}&y:\Spec R\to C,\ \calE \mbox{ is a } \calG\mbox{-torsor on } C_R,\\
&\beta:\calE|_{C_R-\Ga_y}\cong \calE^0|_{C_R-\Ga_y} \mbox{ is a
trivialization}\end{split}\right.\right\},
\end{equation}
where $\Ga_y$ denotes the graph of $y$. This is a formally smooth
ind-scheme over $C$ (\cite[Proposition 5.5]{PZ}).

We also have the jet group
$\calL^+\calG$ and the loop group $\calL\calG$ of $\calG$. To define it, we need some notations.
Let $y: \Spec R\to C$. We
will denote by $\Gamma_y\subset C_R$ the closed subscheme
given by the graph of $y$ and consider  the formal
completion of $C_R$ along $\Gamma_y$, which is an
affine formal scheme and following
  \cite[2.12]{BD} we can also consider the affine scheme $\hat\Gamma_y$
  given by the relative spectrum of the ring of regular functions on that completion.
  There is a natural closed immersion $\Gamma_y\to \hat\Gamma_y$ and we will
denote by $\hat\Gamma_y^\circ:=\hat\Gamma_y-\Gamma_y$ the complement of
the image. In our paper, we will soon specialize to the case $C=\bbA^1=\Spec
k[v]$ so that $y:\Spec R\to C$ is given by $v\mapsto y\in R$ and
therefore $\Ga_y=\Spec R[v]/(v-y)$ and $\hat{\Ga}_y\simeq\Spec
R[[w]]$ and the map $p:\hat{\Ga}_y'\to C_R$ is given by $v\mapsto
w+y$. We will often  write
$\hat\Gamma_y= \Spec R[[v-y]] $. Then
$\hat\Gamma_y^\circ=\Spec R[[v-y]][(v-y)^{-1}]$.

Now, we define $\calL^+\calG$ and $\calL\calG$. For any $k$-algebra $R$,
\begin{equation}\label{J}
\calL^+\calG(R)=\left\{(y,\beta)\ \left|\ y:\Spec R\to C, \beta\in
\calG(\hat{\Ga}_y)\right.\right\},
\end{equation}
and 
\begin{equation}\label{L}
\calL\calG(R)=\left\{(y,\beta)\ \left|\ y:\Spec R\to C, \beta\in
\calG(\hat{\Ga}^\circ_y)\right.\right\}.
\end{equation}
The former is a scheme formally smooth (but not of finite type) over $C$, and the latter is a formally smooth ind-scheme over $C$.

Let us describe the fibers of $\calL\calG,\calL^+\calG,\Gr_\calG$
over $C$. Let $x\in C$ be a closed point. Let $\calO_x$ denote the
completion of the local ring of $C$ at $x$ and $F_x$ be the
fractional field of $\calO_x$. Then
\[(\calL\calG)_x\cong L(\calG_{F_x}),\quad (\calL^+\calG)_x\cong L^+(\calG_{\calO_x}),\quad (\Gr_\calG)_x\cong \Gr_{\calG_{\calO_x}}:=L(\calG_{F_x})/L^+(\calG_{\calO_x}).\]

\quash{
Strictly speaking, we do not need the following remark in the
sequel. However, it helps us clarify the proof of Proposition
\ref{sla}.
\begin{rmk}\label{pitfall}
Let $f:\calG\to\calG'$ be two smooth group schemes over $C$ such
that $\calG|_{C-\{x\}}\cong\calG'|_{C-\{x\}}$. Then clearly the
natural morphism $\calL f:\calL\calG\to\calL\calG'$ will induce
isomorphisms $\calL\calG|_{C-\{x\}}\cong\calL\calG|_{C-\{x\}}$ and
$(\calL\calG)_x\cong(\calL\calG')_x$. However, the morphism $\calL
f$ itself is not necessarily an isomorphism.
\end{rmk}}

The groups $\calL\calG$ and $\calL^+\calG$ naturally act on
$\Gr_\calG$. To see this, let us use the descent lemma of
Beauville-Laszlo (see \cite{BL}, or rather a general form of this
lemma given in \cite[Theorem 2.12.1]{BD}) to show
\begin{lem}\label{second dfn of Gr}The natural map
\[\Gr_\calG(R)\to\left\{(y,\calE,\beta)\ \left|\ \begin{split}&y:\Spec R\to C, \calE \mbox{ is a } \calG\mbox{-torsor
on }\\& \hat{\Ga}_y, \ \beta:\calE|_{\hat{\Ga}^\circ_y}\cong
\calE^0|_{\hat{\Ga}^\circ_y} \mbox{ is a
trivialization}\end{split}\right.\right\}\] is a bijection for each
$R$.
\end{lem}
Then $\calL\calG$ and $\calL^+\calG$ act on $\Gr_\calG$ by changing
the trivialization $\beta$. The trivial $\calG$-torsor gives
$\Gr_\calG\to C$ a section $e$. Then we have the projection
\begin{equation}\label{pr}
\pr: \calL\calG\to \calL\calG\cdot e=\Gr_\calG.
\end{equation}

We need the following lemma in the sequel.
\begin{lem}\label{base change}The formations of $\Gr_\calG,\calL\calG,\calL^+\calG$
commute with any \'{e}tale base change, i.e. if $f:C'\to C$ is
\'{e}tale, then $\Gr_\calG\times_CC'\cong\Gr_{\calG\times_CC'}$,
etc. In addition, the action of $\calL\calG$ on $\Gr_\calG$ also
commutes with any \'{e}tale base change.
\end{lem}
\begin{proof}We have the following observation. Let $y':\Spec R\to C'$ be an $R$-point of $C'$ and $f(y):\Spec R\to
C$ be the corresponding $R$-point of $C$.  Since $f$ is \'{e}tale, the morphism obtained from $f$ by completing along $y'$ and $y$ gives an isomorphism of affine
formal schemes which induces an isomorphism $\hat\Ga_y\simeq\hat\Ga_{y'}$ between the affine spectra of their coordinate rings. In addition, this isomorphism restricts to an isomorphism $\hat{\Ga}^\circ_{y'}\simeq\hat{\Ga}^\circ_{y}$. The lemma
now follows.
\end{proof}

\subsection{The group scheme}\label{grpsch}
We will be mostly interested in the case that $\calG$ is a
Bruhat-Tits group scheme over $C$. Let us specify the meaning of this term. Let
$\eta$ denote the generic point of $C$. Then a smooth group scheme
$\calG$ over $C$ is called a Bruhat-Tits group scheme if
$\calG_\eta$ is (connected) reductive, and for any closed point $y$
of $C$, $\calG_{\calO_y}$ is a parahoric group scheme of
$\calG_{F_y}$.

Now let us specify the Bruhat-Tits group scheme that will be
relevant to us. Let $G_1$ be an almost simple, absolutely simple and
simply-connected, and split over a tamely ramified extension
$\tilde{F}/F$, as in the coherent conjecture. Then we can assume
that $\tilde{F}/F$ is cyclic of order $e=1,2,3$. Let $\ga$ be a
generator of $\Ga=\Gal(\tilde{F}/F)$. For technical reasons, which
is apparent from the statement of Theorem \ref{NSV}, we need the
following well-known result.

\begin{lem}There is a connected reductive group $G$ over $F$, which splits over $\tilde{F}/F$,
such that $G_\der\cong G_1$ and $\xcoch(T)\to\xcoch(T_\ad)$ is
surjective. Here $T$ is a maximal torus of $G$ as in \S \ref{grp
data}.
\end{lem}
For example, if $G_1=\SL_n$ or $\on{Sp}_{2n}$, then $G$ can be chosen as $\GL_n$ and $\on{GSp}_{2n}$ respectively.

We let $(H,B_H,T_H,X)$ be a split pinned group over $\bbZ$, together
with an isomorphism $(G,B,T)\otimes_F\tilde{F}\cong
(H,B_H,T_H)\otimes\tilde{F}$ as in \S \ref{grp data}. Let us choose
the special vertex $v_0$ to identify $A(G,S)$ with $\xcoch(S)_\bbR$,
and $\mathbf{a}$ be the chosen alcove in $A(G,S)$ as in \S \ref{grp data}.
Let $Y\subset\bold S$ as before.

Let $[e]:\bbA^1\to \bbA^1$ be the ramified cover given by $y\to
y^e$. To distinguish these two $\bbA^1$s, let us denote it as
$[e]:\tilde{C}\to C$. The origin of $C$ is denoted by $0$ and the
origin of $\tilde{C}$ is denoted by $\tilde{0}$. Write
$C^\circ=C-\{0\}$ and
$\tilde{C}^{\circ}=\tilde{C}-\{\tilde{0}\}$. Observe that $\Ga$
acts on $H\times\tilde{C}$ naturally. Namely, it acts on the first
factor by pinned automorphisms, and the second by transport of
structures. Let
\[\calG|_{C^{\circ}}=(\Res_{\tilde{C}^{\circ}/C^\circ}(H\times\tilde{C}))^{\Ga}.\]
Then $\calG_{F_0}\cong G$ after choosing some $F_0\cong F$. Now, gluing $\calG|_{C^\circ}$ and $\calG_{\sigma_Y}$ along the the fpqc cover $C=C^\circ\cup\Spec\calO_0$ (see \cite[Lemma 5]{H} for the detailed discussion of the descent theory in this case), we get a group scheme $\calG$
over $C$, satisfying
\begin{enumerate}
\item $\calG_\eta$ is connected reductive  with connected center,
splits over a tamely ramified extension, such
that $(\calG_\eta)_\der$ is simple, absolutely simple, and
simply-connected;

\item For some choice of isomorphism $F_0\cong F$, $\calG_{F_0}\cong G$;

\item For any $y\neq 0$, $\calG_{\calO_y}$ is hyperspecial, (non-canonically) isomorphic to $H\otimes\calO_y$;

\item $\calG_{\calO_0}=\calG_{\sigma_Y}$ under the isomorphism $\calG_{F_0}\cong G$.
\end{enumerate}

A more detailed account of the construction of this group scheme is given in \S \ref{globrootgrp}.
Let us mention that
similar group schemes have been constructed in \cite{HNY,Ri}. For
this group scheme $\calG$, we know that the fiber of $\Gr_\calG$
over $y\neq 0$ is isomorphic to the affine Grassmannian $\Gr_H$ of
$H$, and the fiber over $0$ is isomorphic to the affine flag variety
$\Fl^Y$ of $G$. Likewise, the fiber of $\calL^+\calG$ over $y\neq 0$
is isomorphic to $L^+H$ and the fiber over $0$ is isomorphic to
$L^+\calG_{\sigma_Y}$.

Let $\calT$ be the subgroup scheme of $\calG$, such that
\begin{enumerate}
\item $\calT_\eta$ is a maximal torus of $\calG_\eta$;

\item For any $y\neq 0$, $\calT_{\calO_y}$ is a split torus;

\item $\calT_{F_0}$ is the torus $T$ and $\calT_{\calO_0}$ is the connected N\'{e}ron model
of $\calT_{F_0}$.
\end{enumerate}

We can construct $\calT$ as the neutral connected component of
\begin{equation}\label{globNeron}
\tilde{\calT}=(\Res_{\tilde{C}/C}(T_H\times\tilde{C}))^{\Ga}.
\end{equation}
Note that $\calT$ embeds into $\calG$ naturally. Indeed, under our tameness assumption, $\calT$ is the connected N\'{e}ron model of $(\Res_{\tilde{C}^\circ/C^\circ}(T_H\times\tilde{C}))^{\Ga}$, and $\calT(\calO_0)\subset\calG(\calO_0)$ then the claim follows by the construction of parahoric group schemes as in  \cite[5.2]{BT2}.

\subsection{The global Schubert variety}\label{glob Sch}
It turns out that it is more convenient to base change everything
over $C$ to $\tilde{C}$. Let $u$ (resp. $v$) denote a global
coordinate of $\tilde{C}$ (resp. $C$) such that the map
$[e]:\tilde{C}\to C$ is given by $v \mapsto u^e$. Recall that $0\in
C(k)$ (resp. $\tilde{0}\in\tilde{C}(k)$) is given by $v=0$ (resp.
$u=0$). The crucial step toward the construction of the global
Schubert varieties is the following proposition.

\begin{prop}\label{sla}
For each $\mu\in\xcoch(\calT_\eta)\cong\xcoch(T_H)$, there is a
section
\begin{equation*}
s_\mu:\tilde{C}\to\calL\calT\times_C\tilde{C}
\end{equation*}
such that for any $\tilde{y}\in \tilde{C}(k)$ the element
\[s_\mu(\tilde{y})\in(\calL\calT)_{y}(k)=\calT_{F_{y}}(F_{y}),\quad\quad y=[e](\tilde{y})\]
maps under the Kottwitz homomorphism
$\kappa:\calT_{F_y}(F_y)\to\xcoch(\calT_\eta)_{\Gal(F_y^s/F_y)}$ to
the image of $\mu$ under the natural projection
$\xcoch(\calT_\eta)\to\xcoch(\calT_\eta)_{\Gal(F_y^s/F_y)}$.
\end{prop}
The proposition is obvious for split groups. But for the ramified
groups, the proof is a little bit complicated and only the statement of the proposition will be used
in the main body of the paper. Therefore, those who are only
interested in split groups can skip the proof.
\begin{proof}
Let us first review how to construct an element in $t_\mu\in
T(k\ppart)$ whose image under the Kottwitz homomorphism \eqref{kot}
is $\mu$ under the map $\xcoch(T)\to\xcoch(T)_\Ga$. Let
$k\ppars/k\ppart$ be a finite separable extension of degree $n$ so
that $T_{k\ppars}$ splits, where $s^n=t$. Then $\la(s)\in
T(k\ppars)$. By the construction of the Kottwitz homomorphism (cf. \cite[\S 7]{Ko}), we
can take $t_\la$ to be the image of $\la(s)$ under the norm map
$T(k\ppars)\to T(k\ppart)$.

Now we construct $s_\mu$. Let $\tilde{\calT}$ is as in
\eqref{globNeron}. We will first
construct a section $s_\mu:\tilde{C}\to\calL\tilde{\calT}$ and then
prove it indeed factors as
$s_\mu:\tilde{C}\to\calL\calT\to\calL\tilde{\calT}$. 

Let $\Ga_{[e]}$ denote the graph of $[e]:\tilde{C}\to C$. By
definition,
\[\Hom_C(\tilde{C},\calL\tilde{\calT})=\Hom_C(\hat{\Ga}^\circ_{[e]},\tilde{\calT})=\Hom(\hat{\Ga}^\circ_{[e]}\times_C\tilde{C},T_H)^\Ga,\]
where $\Ga$ acts on $\hat{\Ga}^\circ_{[e]}\times_C\tilde{C}$
via the action on the second factor.

Recall that we have the global coordinates $u,v$ and the map
$[e]:\tilde{C}\to C$ is given by $v\mapsto u^e$. Then
$\calO_{\hat{\Ga}^\circ_{[e]}}\cong k[u](\!(v-u^e)\!)$.
Therefore, the ring of functions on
$\hat{\Ga}^\circ_{[e]}\times_C\tilde{C}$ can be written as
\[A=k[u_1](\!(v-u_1^e)\!)\otimes_{k[v]}k[u_2],\]
where the map $k[v]\to k[u_2]$ is given by $v\mapsto u_2^e$. Let
$\ga$ be a generator of $\Ga=\Aut(\tilde{C}/C)$ acting on $u_2$ as
$u_2\mapsto \xi u_2$, where $\xi$ is a primitive $e$'th root of
unit. For $i=1,\ldots,e$, the element $(\xi^i\otimes u_2-u_1\otimes
1)$ is invertible in $A$, and therefore gives a morphism
\[x_i:\hat{\Ga}^\circ_{[e]}\times_C\tilde{C}\to\bbG_m.\]
Clearly $x_i\circ\ga =x_{i+1}$ (as usual, $x_{i+e}=x_i$).

Now choose a basis $\omega_1,\ldots,\omega_\ell$ of $\xch(T_H)$. Let
us define
\[s_\mu:\hat{\Ga}^\circ_{[e]}\times_C\tilde{C}\to T_H\]
as
\[\omega_j(s_\mu)=x_1^{(\mu,\ga\omega_j)}x_2^{(\mu,\ga^2\omega_j)}\cdots x_e^{(\mu,\ga^e\omega_j)}.\]
Clearly, $s_\mu$ is independent of the choice of
$\omega_1,\ldots,\omega_\ell$ (however, it depends on the global
coordinate $u$ on $\tilde{C}$). Furthermore, $s_\mu$ is
$\Ga$-equivariant. Therefore, we constructed a section
$s_\mu:\tilde{C}\to\calL\tilde{\calT}$.

Now we prove that this section indeed factors as
$s_\mu:\tilde{C}\to\calL\calT\to\calL\tilde{\calT}$. 
In other words,
the morphism $\hat{\Ga}^\circ_{[e]}\to\tilde{\calT}$ factors as
$\hat{\Ga}^\circ_{[e]}\to\calT\to\tilde{\calT}$. By definition,
$\calT$ is the neutral connected component of $\tilde{\calT}$.
Therefore, it is enough to prove that the image of
$\hat{\Ga}^\circ_{[e]}|_0\to\tilde{\calT}|_0$ lands in the
neutral connected component of $\tilde{\calT}|_0$. Observe that
$\hat{\Ga}^\circ_{[e]}|_0\cong \Spec k(\!(u_1)\!)$. Let
$\tilde{C}_0$ be the fiber of $\tilde{C}\to C$ over $0$ so that
$\tilde{C}_0\cong k[u]/u^e$ with a $\Ga$-action. It has a unique
closed point $\tilde{0}$. Recall that
$\tilde{T}|_0=(\Res_{\tilde{C}_0/k} (T_H\times\tilde{C}_0))^\Ga$ and
therefore, there is a canonical map $\epsilon:\tilde{\calT}|_0\to
T_H^\Ga$ given by adjunction, making the following diagram commute
\[\begin{CD}
\Hom_C(\hat{\Ga}^\circ_{[e]},\tilde{\calT})@>>>\Hom(\hat{\Ga}^\circ_{[e]}|_0,\tilde{\calT}|_0)@>\epsilon>>\Hom(\hat{\Ga}^\circ_{[e]}|_0,T_H^\Ga)\\
@|@|@|\\
\Hom(\hat{\Ga}^\circ_{[e]}\times_C\tilde{C},T_H)^\Ga@>>>\Hom(\hat{\Ga}^\circ_{[e]}|_0\times\tilde{C}_0,T_H)^\Ga@>>>\Hom(\hat{\Ga}^\circ_{[e]}|_0\times\{\tilde{0}\},T_H^\Ga)
\end{CD}\]
In our case $\epsilon(s_\mu):\hat{\Ga}^\circ_{[e]}|_0 \to T_H$
is given by
\[\omega_j(\epsilon(s_\mu))=(-u_1)^{(\sum_{\ga\in\Ga}\ga\mu,\omega_j)}.\]
In other words, $\epsilon(s_\mu)$ is the composition
\[\hat{\Ga}^\circ_{[e]}|_0\stackrel{-u_1}{\longrightarrow}\bbG_m\stackrel{\sum_{\ga\in\Ga}\ga\mu}{\longrightarrow}T_H.\]
Since for any $\Ga$-invariant coweight $\mu$, the image
$\mu:\bbG_m\to T_H^\Ga$ lands in the neutral connected component of
$T_H^\Ga$ (the torus part), $s_\mu:\tilde{C}\to\calL\tilde{\calT}$
factors through $\tilde{C}\to\calL\calT\to\calL\tilde{\calT}$.

Finally, let us check that
$s_\mu:\tilde{C}\to\calL\calT\times_C\tilde{C}$ satisfies the
desired properties as claimed in the proposition.

Let $\tilde{y}\in \tilde{C}(k)$ be a closed point given by $u\mapsto
\tilde{y}\in k$. Then $s_\mu(\tilde{y})$ corresponds to
$s_\mu(\tilde{y}):\Spec k(\!(v-\tilde{y}^e)\!)\otimes_{k[v]}k[u_2]\to T_H$
given by
\[\omega_j(s_\mu(\tilde{y}))=\prod_{i=1}^e(\xi^i1\otimes u_2-y)^{(\mu,\ga^i\omega_j)}.\]
If $\tilde{y}=0$, the assertion of the proposition follows directly
from the review of the construction of $t_\mu$ at the beginning. If
$\tilde{y}\neq 0$, let $w=1\otimes u_2-y$. Then
\[\prod_{i=1}^e(\xi^i1\otimes u_2-y)^{(\mu,\ga^i\omega_j)}=w^{(\mu,\omega_j)}f(w)\]
where \[f(w)=\prod_{i=1}^{e-1}(\xi^i1\otimes
u_2-y)^{(\mu,\ga^i\omega_j)}\in k[[w]]^\times.\] Therefore, as an
element in $T_H(k(\!(w)\!))$, which is canonically isomorphic to
$\calL\calT_{y}$, $s_\mu(\tilde{y})$ maps to $\mu$ under the
Kottwitz homomorphism.
\end{proof}

\begin{rmk}Note that the natural map $\calL\calT\to\calL\tilde\calT$ induces isomorphisms $\calL\calT|_{C^\circ}\simeq\calL\tilde\calT|_{C^\circ}$ and $\calL\calT|_0\to\calL\tilde\calT|_0$. But itself is not an isomorphism.
\end{rmk}

\begin{rmk}For a general $\mu$, there is no such section $C\to\calL\calT$ satisfying the property of the proposition.
This is the reason that we want to base change everything over $C$
to $\tilde{C}$. However, if $\mu\in\xcoch(T)$ is defined over $F$,
then $s_\mu$ indeed descents to a section $C\to\calL\calT$. This
means that in this case the variety $\bGr_{\calG,\mu}$ defined
below, which a priori is a variety over $\tilde{C}$, descends to a
variety over $C$. One can summarize this by saying that
$\bGr_{\calG,\mu}$ is defined over the "reflex field" of the geometric conjugacy class $\{\mu\}$ (which is the same of the field of definition of $\mu$ as $G$ is quasi-split over $F$ (\cite[Lemma 1.1.3]{Ko1})). The
same phenomenon appears in the theory of Shimura varieties.
\end{rmk}

The composition of $s_\mu$ and the natural morphism (see \eqref{pr})
$\pr:\calL\calT\to\Gr_\calT$  (resp. $\calL\calT\to\calL\calG$)
gives a section $\tilde{C}\to\Gr_\calT\times_C\tilde{C}$ (resp.
$\tilde{C}\to\calL\calG\times_C\tilde{C}$), which is still denoted
by $s_\mu$.

The construction of $\tilde{C}\to\calL\calT\times_C\tilde{C}$ will
depend on the choice of the global coordinate $u$ of $\tilde{C}$,
but the section $s_\mu:\tilde{C}\to \Gr_\calT\times_C\tilde{C}$ does
not. Indeed, there is the following moduli interpretation of such
section. Recall that $\Gr_{\calT}$ is ind-proper over $C$
(\cite{H}), and therefore, $s_\mu$ is uniquely determined by a
section
$\tilde{C}^{\circ}\to\Gr_\calT\times_C\tilde{C}^{\circ}\cong\Gr_{T_H\times\tilde{C}^{\circ}}$
(by Lemma \ref{base change}). Then this section, under the moduli
interpretation of $\Gr_{T_H\times\tilde{C}^{\circ}}$, is given as
follows: let $\Delta$ be the diagonal of $\tilde{C}^{\circ}\times\tilde{C}^\circ$,
and $\calO_{(\tilde{C}^{\circ})^2}(\mu\Delta)$ be the $T_H$-torsor
on $(\tilde{C}^{\circ})^2$, such that for any weight $\nu$ of
$T_H$, the associated line bundle is
$\calO_{\tilde{C}^{\circ}}((\mu,\nu)\Delta)$. This $T_H$-torsor
has a canonical trivialization away from $\Delta$.

\begin{lem}\label{moduli interpretation of sla}
The map $s_\mu:\tilde{C}^{\circ}\to\Gr_\calT$ corresponds to
$(\calE,\beta)$, where $\calE$ is the $T_H$-torsor
$\calO_{(\tilde{C}^{\circ})^2}(\mu\Delta)$, and $\beta$ is its
canonical trivialization over $(\tilde{C}^{\circ})^2-\Delta$.
\end{lem}
\begin{proof}The Kottwitz homomorphism
$\kappa:LT_H(k)\to\xcoch(T_H)$ induces an isomorphism $\Gr_H(k)\cong
LT_H(k)/L^+T_H(k)$. On the other hand, recall that if we fix a point
$x$ on the curve $\tilde{C}$, we can interpret $\Gr_{T_H}$ as the
set of $(\calE,\beta)$, where $\calE$ is an $T_H$-torsor and $\beta$
is a trivialization of $\calE$ away from $x$. Under this
interpretation, any $t_\mu\in\xcoch(T_H)$ is interpreted as the
$T_H$-torsor $\calO_{\tilde{C}}(\mu x)$\footnote{The reason that
$t_\mu$ represents $\calO_{\tilde{C}}(\mu x)$ rather than
$\calO_{\tilde{C}}(-\mu x)$ is due to the original sign convention
of the Kottwitz homomorphism in \cite{Ko}.}, with its canonical
trivialization away from $x$. Then the lemma is clear.
\end{proof}

By composing with the natural morphism $\Gr_\calT\to\Gr_\calG$, we obtain a
section of $\Gr_\calG\times_C\tilde{C}$, still denoted by $s_\mu$.

\medskip

\noindent\bf Notation. \rm In what follows, we denote
$\Gr_\calG\times_C\tilde{C}$ (resp. $\calL^+\calG\times_C\tilde{C}$,
resp. $\calL\calG\times_C\tilde{C}$) by $\wGr_\calG$ (resp. $\wJG$,
resp. $\wLG$).

\begin{dfn}\label{globschubert}For each $\mu\in\xcoch(\calT_\eta)\cong\xcoch(T_H)$, the global Schubert variety
$\bGr_{\calG,\mu}$ is the minimal $\wJG$-stable irreducible closed
subvariety of $\wGr_\calG$ that contains $s_\mu$.
\end{dfn}
Let us emphasize that $\bGr_{\calG,\mu}$ is not a subvariety of
$\Gr_\calG$. Rather, it lies in $\Gr_\calG\times_C\tilde{C}$. Recall
that for any $\mu\in\xcoch(T)$, one defines a subset
$\Adm^Y(\mu)\subset\widetilde{W}$ as in \eqref{Adm}. The main geometric property
of $\bGr_{\calG,\mu}$ which we will prove in this paper is as follows.

\begin{Thm}\label{Fibers}\emph{Assume that the group $G$ and the group scheme $\calG$ are as in \S \ref{grpsch}. Let $y$ be a closed point of $\tilde{C}$. Then
\[(\bGr_{\calG,\mu})_y\cong\left\{\begin{array}{ll}\bigcup_{w\in\Adm^Y(\mu)}\Fl^Y_w & y=\tilde{0}\\ \bGr_\mu & y\neq \tilde{0}.\end{array}\right.\]
In particular, all the fibers are reduced.}
\end{Thm}

We first prove the easy part of the theorem.

\begin{lem}
$(\bGr_{\calG,\mu})_y\cong\bGr_\mu$ for $y\neq \tilde{0}$.
\end{lem}
\begin{proof}
Write $\tilde{C}^{\circ}=\tilde{C}-\tilde{0}$. We want to show
that $\bGr_{\calG,\mu}|_{\tilde{C}^{\circ}}$ is isomorphic to
$\bGr_\mu\times\tilde{C}^{\circ}$. First we have a canonical
isomorphism
\begin{equation}\label{triv1}
\calG\times_C\tilde{C}^{\circ}\cong H\times\tilde{C}^{\circ}
\end{equation}
and therefore by Lemma \ref{base change}
\begin{equation}\label{triv2}
\Gr_\calG\times_C\tilde{C}^{\circ}\cong\Gr_{H\times\tilde{C}^{\circ}},\quad
\calL\calG\times_C\tilde{C}^{\circ}\cong \calL(H\times\tilde{C}).
\end{equation}
Secondly, $\tilde{C}^{\circ}\cong\bbG_m$ which admits a global
coordinate $u$ so that $\calL(H\times\tilde{C}^{\circ})\cong
LH\times\tilde{C}^{\circ}$ and
$\Gr_{H\times\tilde{C}^{\circ}}\cong\Gr_H\times\tilde{C}^{\circ}$.
Finally, by Lemma \ref{moduli interpretation of sla}, the section
$s_\mu:\tilde{C}^{\circ}\to\Gr_\calG\times_C\tilde{C}^{\circ}\cong
\Gr_H\times\tilde{C}^{\circ}$ satisfies
$s_\mu(\tilde{C}^{\circ})\subset\bGr_\mu\times\tilde{C}^{\circ}$.
\end{proof}
Using this lemma, we see that it is enough to make the following
convention.

\medskip

\noindent\bf Convention. \rm When we discuss $\bGr_{\calG,\mu}$, we
will assume that $\mu\in\xcoch(T_H)$ is dominant with respect to the
chosen Borel $B_H$ as in \S \ref{grp data} and \S \ref{grpsch}.

\medskip

At this moment, we can also see that
\begin{lem}\label{easy}
The scheme
$(\bGr_{\calG,\mu})_{\tilde{0}}\subset(\Gr_{\calG})_0\cong\Fl^Y$
contains $\Fl^Y_w$ for $w\in\Adm^Y(\mu)$.
\end{lem}
\begin{proof}
Clearly, it is enough to show that
$\Fl^Y_\la\subset(\bGr_{\calG,\mu})_{\tilde{0}}$ for any $\la$ in
$\Lambda$, where $\Lambda$ is the ${W_0}$-orbit in $\xcoch(T)_\Ga$
containing $\mu$ as constructed in \S \ref{grp data}. Observe that
$\bGr_{\calG,\mu}$ is the flat closure of
$\bGr_{\calG,\mu}|_{\tilde{C}^{\circ}}$ in $\wGr_\calG$, since
the later is clearly $\wJG$-stable. Then the claim follows from
that for any $\la\in\xcoch(T_H)$ in the $\bar{W}$-orbit of $\mu$,
$s_\la(\tilde{0})\in\Fl^Y_\la$ and
$s_\la(\tilde{C}^{\circ})\subset
\bGr_{\calG,\mu}|_{\tilde{C}^{\circ}}\cong\bGr_\mu\times\tilde{C}^{\circ}$.
\end{proof}

To prove the theorem, it is remains to show that

\begin{thm}\label{top fiber}
The underlying reduced subscheme of $(\bGr_{\calG,\mu})_{\tilde{0}}$
is $\bigcup_{w\in\Adm^Y(\mu)}\Fl^Y_w$.
\end{thm}

\begin{thm}\label{sch fiber}
$(\bGr_{\calG,\mu})_{\tilde{0}}$ is reduced.
\end{thm}

By the same argument as in \cite[9.2.1]{PZ}, $\bGr_{\calG,\mu}$ is normal. We conjecture that it is Cohen-Macaulay as well.

\section{Line bundles on $\Gr_\calG$ and $\Bun_\calG$}\label{lines on globaff}
This subsection explains why Theorem \ref{Fibers} and Theorem
\ref{MainII} are equivalent to each other. The key ingredients are
the line bundles on the global affine Grassmannian $\Gr_\calG$.
Observe that $\Gr_\calG$ can be disconnected. This will create some
complications in trying to determine the line bundle on $\Gr_\calG$ directly.
Instead, we will pass to its group scheme $\calG_\der$ (defined below), whose generic
fiber then is simply-connected so that we can apply the results of Heinloth
\cite{H} directly.

\subsection{Line bundles on $\Gr_\calG$ and $\Bun_\calG$}\label{sub;lines on globaff}

In this subsection, we temporary assume that $C$ is a smooth curve
over $k$ and $\calG$ is a Bruhat-Tits group scheme over $C$ such
that $\calG_\eta$ is almost simple, absolutely simple, and
simply-connected.

\begin{prop}\label{cal of cc}
Let $\calL$ be a line bundle on $\Gr_{\calG}$. Then the function
$c_\calL$ that associates to every $y\in C(k)$ the central charge of
the restriction of $\calL$ to $(\Gr_{\calG})_y$ is constant.
\end{prop}
This proposition implies that the statement in the last sentence of
the first paragraph in p. 502 of \cite{H} is not correct.
\begin{proof}
Let $\Pic(\Gr_{\calG}/C)$ denote the relative sheaf of Picard
groups over $C$. As explained in \cite{H}, this is an \'{e}tale
sheaf over $C$. Let
$D=\on{Ram}(\calG)$ be the set of points of $C$ such that for every
$y\in\on{Ram}(\calG)$, the fiber $\calG_y$ is not semisimple. This is a finite set. Then there is a short exact sequence
\begin{equation}\label{global Pic}
1\to \prod_{y\in D}\xch(\calG_y)\to\Pic(\Gr_{\calG}/C)\to \frakc\to 1,
\end{equation}
where  $\frakc$ is a constructible sheaf, with all fibers isomorphic to
$\bbZ$ and is constant on $C-D$.

According to the description of the sheaf $\frakc$ in Remark 19 (3)
of \emph{loc. cit.}, if $\calL$ is a line bundle on $\Gr_\calG$ such
that $c_\calL(y)=0$ for some $y\in C(k)$, then $c_\calL=0$.
Therefore, to prove the proposition, it is enough to construct one
line bundle $\calL_{2c}$ on $\Gr_{\calG}$, such that
$c_{\calL_{2c}}$ is constant on $C$.

Let $\calV_0=\Lie\calG$ be the Lie algebra of $\calG$. This is a
locally free $\calO_C$-module on $C$ of rank $\dim_\eta\calG_\eta$,
on which $\calG$ acts by the adjoint representation. This induces a
morphism $\calG\to\GL(\calV_0)$, and therefore a morphism
$i:\Gr_{\calG}\to\Gr_{\GL(\calV_0)}$. Let $\calL_{\det}$ denote the
determinant line bundle on $\Gr_{\GL(\calV_0)}$. Let us recall its
construction: We want to associate to every $\Spec
R\to\Gr_{\GL(\calV_0)}$ a line bundle on $\Spec R$ in a compatible
way. Recall that a morphism $\Spec R\to\Gr_{\GL(\calV_0)}$ represents a
morphism $y\in C(R)$, a vector bundle $\calV$ on $C_R$ and an
isomorphism $\calV|_{C_R-\Ga_y}\cong\calV_0|_{C_R-\Ga_y}$. There
exists some $N$ large enough such that
\[\calV_0(-N\Ga_y)\subset\calV\subset\calV_0(N\Ga_y)\]
and $\calV_0(N\Ga_y)/\calV$ is $R$-flat. Then the line bundle on $\Spec R$ is
\[\det(\calV_0(N\Ga_y)/\calV)\otimes\det(\calV_0(N\Ga_y)/\calV_0)^{-1},\]
which is independent of the choice of $N$ up to a canonical isomorphism.

The pullback $i^*\calL_{\det}$ is a line bundle on $\Gr_{\calG}$,
which will be our $\calL_{2c}$. To see this is the desired line
bundle, we need to calculate its central charge when restricted to
each $y\in C(k)$. Let $D=\on{Ram}(\calG)$. First consider $y\in
C-D$. Then the map $i_y:(\Gr_{\calG})_y\to(\Gr_{\GL(\calV_0)})_y$ is
just \[\Gr_{H}\to\Gr_{\GL(\Lie H)},\] where $H$ is the split
Chevalley group over $\bbZ$ such that $G\otimes k(\eta)^s\cong
H\otimes k(\eta)^s$. It is well known that in this case
$i^*_y\calL_{\det}$ over $y$ has central charge $2h^\vee$, where
$h^\vee$ is the dual Coxeter number of $H$ (in fact, this statement
is a consequence of the following argument).

It remains to calculate the central charge of $\calL_{2c}$ over
$y\in D$. Without loss of generality, we can assume that $D$
consists of one point, denoted by $0$. So let $y=0$ and $G=\calG_{F_0}$. Then the closed embedding
$i_0:(\Gr_{\calG})_0\to(\Gr_{\GL(\calV_0)})_0$ is just
\[LG/L^+\calG_{\calO_0}\to\Gr_{\GL(\Lie \calG_{\calO_0})}.\] Let us
first assume that $\calG_{\calO_0}$ is an Iwahori group scheme of
$\calG_{F_0}$. Write $I= L^+\calG_{\calO_0}$ and $\Fl=LG/I$ as usual. We
claim that in this case
\begin{lem}\label{adj} We have an isomorphism
$i^*_0\calL_{\det}\cong \calL(2\sum_{i\in\bold S}\epsilon_i)$.
\end{lem}
Assuming this fact, we find the central charge of
$i^*_0\calL_{\det}$ is $2\sum_{i\in \bold S}a_i^\vee$. By checking
all the affine Dynkin diagrams, we find that \[\sum_{i\in\bold
S}a_i^\vee=h^\vee.\] In fact, we find that for affine Dynkin diagrams
$X^{(r)}_N$, where $X=A,B,C,D,E,F,G$ and $r=1,2,3$, the sum $\sum
a_i^\vee$ is independent of $r$ (see \cite[Remark 6.1]{Kac}), and it
is well-known (or by definition) that for $r=1$, $\sum
a_i^\vee=h^\vee$. Therefore, the proposition follows in this case.

Now we prove Lemma \ref{adj}. This is equivalent to proving that the
restriction of $i^*_0\calL_{\det}$ to each $\bbP^1_j$ (whose
definition is given in \S \ref{loop grp and flag var}) is isomorphic
to $\calO_{\bbP^1}(2)$. Recall that a $k$-point $gI\in \Fl$ corresponds to a pro-algebraic subgroup of $L\calG_{F_0}$ given by $I':=gIg^{-1}$, which is the jet group of an Iwahori group scheme of $\calG_{F_0}$. By abuse of notation, we still denote this Iwahori group scheme by $I'$.
Then $\Fl\to\Gr_{\GL(\Lie \calG_{\calO_0})}$ maps an
Iwahori group scheme $I'$ of $G$ to its Lie algebra $\Lie I'$, which
is a free $\calO_0$-module, together with the canonical isomorphism
$\Lie I'\otimes F_0\cong \Lie G\cong\Lie I\otimes F_0$.

For $j\in\bold S$, let $P_j$ be the minimal parahoric (but not
Iwahori) group scheme corresponding to $j$. Then the subscheme
$\bbP^1_j\subset\Fl$ classifies the Iwahori group schemes of  $G$
that map to $P_j$. Let $P_j^u\to P_j$ be the "unipotent radical" of
$P_j$. More precisely, $P_j^u$ is smooth over $\calO_0$ with
$P_j^u\otimes F_0=G$ and the special fiber of $P_j^u$ maps onto the
unipotent radical of the special fiber of $P_j$. If $I'$ is an
Iwahori group scheme of $G$ that maps to $P_j$, then
\[\Lie P_j^u\subset \Lie I'\subset \Lie P_j.\]
Let $\bar{P}^{\on{red}}_j$ be the reductive quotient of the special
fiber of $P_j$. Then $\bar{P}^{\on{red}}_j$ is isomorphic to
$\GL_2,\SL_2$ or $\SO_3$ over $k$. Let
$\Gr(2,\Lie\bar{P}^{\on{red}}_j)\cong \bbP^2$ denote the
Grassmannian (over $k$) of $2$-planes in the three dimensional
vector space $\Lie\bar{P}^{\on{red}}_j$. We have the following
commutative diagram
\[\begin{CD}
\bbP_j^1@>>>\Gr(2,\Lie\bar{P}^{\on{red}}_j)\\
@VVV@VVV\\
\Fl@>>>\Gr_{\GL(\Lie \calG_{\calO_0})}
\end{CD}\]
where $\bbP_j^1\to \bbG(2,\Lie\bar{P}^{\on{red}}_j)$ is given by
\[I'\mapsto (\Lie I'/\Lie P_j^u\subset \Lie P_j/\Lie P_j^u)\cong\Lie\bar{P}^{\on{red}}_j)\]
and $\Gr(2,\Lie\bar{P}^{\on{red}}_j)\to\Gr_{\GL(\Lie I)}$ is given by realizing that $\Gr(2,\Lie\bar{P}^{\on{red}}_j)$ represents the free $\calO_0$-modules that are in between $\Lie P^u_j$ and $\Lie P_j$. Observe that the degree of the map $\bbP^1_j\to \Gr(2,\Lie\bar{P}^{\on{red}}_j)\simeq\bbP^2$  is two as it is just the map that sends a Borel subgroup of $\SL_2$ to the two-dimensional vector subspace of $\fraks\frakl_2$ given by the Lie algebra of the Borel subgroup.

By construction, the restriction of $\calL_{\det}$ to
$\Gr(2,\Lie\bar{P}^{\on{red}}_j)$ is the (positive) determinant
line bundle on $\bbG(2,\Lie\bar{P}^{\on{red}}_j)$, or
$\calO_{\bbP^2}(1)$. Therefore, the restriction of $\calL_{\det}$ to
$\bbP_j^1$ is isomorphic to $\calO_{\bbP^1}(2)$. This finishes the
proof of Lemma \ref{adj} and therefore the proposition in the case
$\calG_{\calO_0}$ is Iwahori.

\medskip

Now let $\calG_{\calO_0}$ be a general parahoric group scheme. Let
$\calG'$ be the group scheme over $C$ together with $\calG'\to\calG$
which is an isomorphism over $C-\{0\}$ and $\calG'_{\calO_0}$ is
Iwahori. Let $\calV_0=\Lie\calG$ and $\calV'_0=\Lie\calG'$. We have
the natural map \[p:(\Gr_{\calG'})_0\to(\Gr_{\calG})_0\] induced
from $\calG'\to\calG$ and the maps
\[i:\Gr_{\calG}\to\Gr_{\GL(\calV_0)},\quad i':\Gr_{\calG'}\to\Gr_{\GL(\calV'_0)}.\]
Let $\calL_{\det}$ (reps. $\calL'_{\det}$) be the determinant line bundle on $\Gr_{\GL(\calV_0)}$ (reps. $\Gr_{\GL(\calV'_0)}$). We need to show that $p^*i^*_0\calL_{\det}$ and $i'^*_0\calL'_{\det}$
have the same central charge (observe that these two line bundles
are not isomorphic). From this, we conclude that the central charge
of $i^*\calL_{\det}$ is also constant along $C$.

Let us extend $\calG$ and $\calG'$ to group schemes over the
complete curve $\bar{C}$ such that
$\calG|_{\bar{C}-\{0\}}=\calG'|_{\bar{C}-\{0\}}$. Let $\Bun_{\calG}$
(resp. $\Bun_{\calG'}$) be the moduli stack of $\calG$-torsors
($\calG'$-torsors) on $\bar{C}$.  Let $\calG_0,\calG'_0$ be the
restriction of the two group schemes over $0\in C$, and let $P$ be
the image of $\calG'_0\to\calG_0$. This is indeed a Borel subgroup
of $\calG'_0$. Recall that by restricting a $\calG'$-torsor to $0\in
C$, we obtain a map $(\Gr_{\calG'})_0\stackrel{r}{\to}\bbB
\calG'_0$, and we have the similar map for $\calG$. Then we have the
following diagram with both squares Cartesian
\begin{equation}\label{car}
\xymatrix{
(\Gr_{\calG'})_0\ar[r]\ar[d]&\Bun_{\calG'}\ar^r[r]\ar[d]&\bbB \calG'_0\ar[r]&\bbB P\ar[d]\\
(\Gr_{\calG})_0\ar[r]&\Bun_{\calG}\ar^r[rr]&&\bbB \calG_0.}
\end{equation}
Indeed, it is clear that the left square is Cartesian because
$\calG|_{\bar{C}-\{0\}}=\calG'|_{\bar{C}-\{0\}}$. The fact that
the second square is Cartesian is established in Proposition
\ref{cart}.

Let
$y:\Spec R\to(\Gr_{\calG'})_0$ be a morphism given by $(\calE,\beta)$,
where $\calE$ is a $\calG'$-torsor on $C_R$. Then we have the
natural short exact sequence
\[0\to\ad\calE\to \ad(\calE\times^{\calG'}\calG)\to\calE\times^{\calG'}(\Lie\calG/\Lie\calG')\to 0.\]
On the other hand, $p:(\Gr_{\calG'})_0\to(\Gr_{\calG})_0$ is a
relatively smooth morphism since $\bbB P\to\bbB \calG'_0$ is smooth.
Let $\calT_p$ denote the relative tangent sheaf. We claim that
$\calE\times^{\calG'}(\Lie\calG/\Lie\calG')\cong y^*\calT_p$, where
$y^*\calT_p$ is the sheaf on $\Spec R$, regarded as a sheaf on $C_R$
via the closed embedding $\{0\}\times \Spec R=:\{0\}_R\to C_R$. But
this follows from \eqref{car} and
\[\calE\times^{\calG'}(\Lie\calG/\Lie\calG')\cong \calE|_{\{0\}_R}\times^{\calG'_0}(\Lie \calG_0/\Lie P)\cong(\calE|_{\{0\}_R}\times^{\calG'_0}P)\times^{P}(\Lie\calG_0/\Lie P).\]
Therefore, we have
\begin{equation}\label{An exact sequence}
0\to\ad\ \calE\to\ad(\calE\times^{\calG''}\calG')\to y^*\calT_p\to
0,
\end{equation}

Let us finish the proof that $p^*i^*_0\calL_{\det}$ and $i'^*_0\calL'_{\det}$ have the same central charge and therefore the proof of the proposition. From the above lemma,
\begin{equation}\label{compare line bundles}
p^*i^*_0\calL_{\det}\cong i'^*_0\calL'_{\det}\otimes \det(\calT_p).
\end{equation}
So it is enough to prove that $\det\calT_p$ as a line bundle on
$(\Gr_{\calG'})_0$ has central charge zero. But from \eqref{car},
$\det\calT_p$ is a pullback of some line bundle from $\bbB P$, and
hence from $\bbB\calG'_0$, which has zero central charge by
\eqref{ker of c}.
\end{proof}

Now, we assume that $C$ is a complete curve and let $\Bun_\calG$ be the moduli stack of $\calG$-torsors on $C$. Let $\Pic(\Bun_{\calG})$ be the Picard group of rigidified line bundles (trivialized over the trivial $\calG$-torsor) on $\Bun_{\calG}$. Let $D=\on{Ram}(\calG)$. Observe that $\prod_{y\in C(k)}\xch(\calG_y)=\prod_{y\in D}\xch(\calG_y)$. Fix $0\in C(k)$. Let $\Fl^Y=L\calG_{F_0}/L^+\calG_{\calO_0}$, which is a partial affine flag variety of $\calG_{F_0}$. According to \cite[\S 7]{H}, we have the following commutative diagram
\[\begin{CD}
0@>>> \prod_{y\in C(k)}\xch(\calG_y)@>>> \Pic(\Bun_{\calG})@>>> \bbZ@>>> 0\\
@.@VVV@VVV@VVV@.\\
0@>>>\xch(\calG_0)@>>>\Pic(\Fl^Y)@>c>>\bbZ@>>>0
\end{CD}\]
The left vertical arrow is the projection to the factor
corresponding to $0$ and the right vertical arrow is injective (but
not necessarily surjective). Probably, one can show that
$\Pic(\Bun_\calG)\to\bbZ$ is in fact given by
$\Pic(\Bun_\calG)\to\Ga(C,\Pic(\Gr_\calG/C))\to\Ga(C,\frakc)\cong\bbZ$
and the right vertical arrow is the natural restriction map
$\Ga(C,\frakc)\to \frakc|_0$.  Here we will not need to show this. We can however, use the above diagram to show that,
for any $\calL\in\Pic(\Fl^Y)$,
a certain tensor power of it will descend to a line bundle on
$\Bun_\calG$. Therefore we conclude

\begin{cor}\label{line and central charge}
Let $C$ be a smooth but (not necessarily complete) curve and let $\calG$
be a Bruhat-Tits group scheme over $C$ such that $\calG_\eta$ is
almost simple, absolutely simple and simply-connected. Let $H$ be
the split Chevalley group over $\bbZ$ such that $\calG\otimes
k(\eta)^s\cong H\otimes k(\eta)^s$. Let $0\in C(k)$ and let $\calL$
be a line bundle on $\Fl^Y=L\calG_{F_0}/L^+\calG_{\calO_0}$. Then
there is a line bundle on $\Gr_{\calG}$, whose restriction to
$(\Gr_{\calG})_0\cong\Fl^Y$ is isomorphic to $\calL^n$ for some
$n\geq 1$, and whose restriction to $(\Gr_{\calG})_y\cong\Gr_{H}
(y\not\in\on{Ram}(\calG))$ is isomorphic to $\calL_{b}^{nc(\calL)}$, where $\calL_b$ is the ample generator of $\Pic(\Gr_H)\cong\bbZ$.
\end{cor}

\begin{proof}Let $\bar{C}$ be a complete curve containing $C$. We extend $\calG$ to a
Bruhat-Tits group scheme over $\bar{C}$. Then some tensor power
$\calL^n$ of $\calL$ descends to a line bundle $\calL'$ on
$\Bun_{\calG}$. Let $h_{\on{glob}}:\Gr_{\calG}\to\Bun_{\calG}$ be
the natural projection. Then $h_{\on{glob}}^*\calL'$ is a line
bundle on $\Gr_{\calG}$ whose restriction to $(\Gr_{\calG})_0$ is
isomorphic to $\calL^n$, and whose restriction to
$(\Gr_{\calG})_y\cong\Gr_{H}$ $(y\not\in\on{Ram}(\calG))$ has
central charge $c(\calL^n)$, and therefore is isomorphic to
$\calL_b^{nc(\calL)}$.
\end{proof}

\subsection{Theorem \ref{Fibers} is equivalent to Theorem \ref{MainII}}
Let us begin with a general construction. Let $\calG$ be a
Bruhat-Tits group scheme over a curve $C$. Then away from a finite
subset $D\subset C$, $\calG|_{C-D}$ is reductive. Let
$\calG_{\der}|_{C-D}$ be the derived group of $\calG|_{C-D}$ so that
for $y\in C(k)$, $(\calG_{\der})_{F_y}$ is the derived group of
$\calG_{F_y}$ (SGA III, Expos\'{e} XXII 6.2]). It is known that there is a canonical bijection
between the facets in the building of $(\calG_{\der})_{F_y}$ and
those in the building of $\calG_{F_y}$, and under this bijection,
the corresponding parahoric group scheme for $(\calG_{\der})_{F_y}$
maps to the corresponding parahoric group scheme for $\calG_{F_y}$. For example, see \cite{HR} Proposition 3 and its proof for the last statement.
Therefore, we can extend $\calG_{\der}|_{C-D}$ to a Bruhat-Tits
group scheme over $C$ together with a morphism $\calG_\der\to\calG$,
such that for all $y\in D$,
$(\calG_{\der})_{\calO_y}\to\calG_{\calO_y}$ is the morphism of
parahoric group schemes given by the facet determined by
$\calG_{\calO_y}$.

\begin{dfn}\label{der}
The group scheme $\calG_\der$ together with the morphism
$\calG_\der\to\calG$ is called the derived group scheme of $\calG$.
\end{dfn}

Now let us specialize the group $\calG$ to be the Bruhat-Tits group
scheme over $C=\bbA^1$ as defined in \S \ref{grpsch}. Let us denote
$\calG_1=\calG_\der$ for simplicity. Let $C^{\circ}=C-\{0\}$.
Observe that $\calG_1|_{C^{\circ}}$ is reductive and
$(\calG_1)_{F_0}\cong G_1=G_\der$ and for $y\neq 0$, $(\calG_1)_{\calO_y}$
is hyperspecial for $H_{\der}\otimes \calO_y$. In addition,
$(\calG_1)_\eta$ is simply-connected.

Let us explain why Theorem \ref{Fibers} and Theorem \ref{MainII} are
equivalent. The natural morphism $\calG_1\to\calG$ induces a morphism
$\Gr_{\calG_1}\to\Gr_{\calG}$. One can show that this is a closed
embedding (we will not use this fact; it follows however, a posteriori, from
the argument below). But at least it follows directly from
\cite[\S 6]{PR1} that both $(\Gr_{\calG_1})_0\to(\Gr_{\calG})_0$ and
$\Gr_{\calG_1}|_{C^{\circ}}\to\Gr_{\calG}|_{C^{\circ}}$ are
closed immersions. These induce isomorphisms from $(\Gr_{\calG_1})_0$ and $(\Gr_{\calG_1})_{C^\circ}$ to the
reduced subschemes of the neutral connected component of $(\Gr_{\calG})_0$ and of $\Gr_{\calG}|_{C^\circ}$ respectively. 
Let $\mu\in\xcoch(T)$, and let $\bGr_{\calG,\mu}$ be the
corresponding global Schubert variety as in \S \ref{glob Sch}.
Recall the section $s_\mu$ from Proposition \ref{sla}. Regard it as
a section of $\wLG$, which acts on $\wGr_\calG$. Then
\[s_\mu^{-1}\bGr_{\calG,\mu}|_{\tilde{C}^{\circ}}\subset\wGr_{\calG_1}|_{\tilde{C}^{\circ}}.\]
This follows from $t_\mu^{-1}\bGr_\mu\subset \Gr_{H_\der}$ for any
$\mu\in\xcoch(T_H)$, where $t_\mu$ is considered as any lifting of
$t_\mu\in\widetilde{W}$ to $T_H(F)$. Let $\Gr_{\calG_1,\leq\mu}$ be
the flat closure of
$s_\mu^{-1}\bGr_{\calG,\mu}|_{\tilde{C}^{\circ}}$ in
$\wGr_{\calG'}$. We have the natural map
\[\Gr_{\calG_1,\leq\mu}\to s_\mu^{-1}\bGr_{\calG,\mu},\]
which induces a closed embedding
$(\Gr_{\calG_1,\leq\mu})_{\tilde{0}}\to
(s_\mu^{-1}\bGr_{\calG,\mu})_{\tilde{0}}$ since $\Fl^Y_\s\to\Fl^Y$
is a closed embedding. By flatness, this necessarily implies that
$(\Gr_{\calG_1,\leq\mu})_{\tilde{0}}\cong
(s_\mu^{-1}\bGr_{\calG,\mu})_{\tilde{0}}$. To see this, let $x\in (\Gr_{\calG_1,\leq\mu})_{\tilde{0}}$ and $y$ be its image in $(s_\mu^{-1}\bGr_{\calG,\mu})_{\tilde{0}}$. Let $A$ and $B$ be their local rings respectively. Let $u$ be a local coordinate around $\tilde{0}$. Then $B\to A$ is injective, since $B[u^{-1}]\to A[u^{-1}]$ is an isomorphism and $B$ has no $u$-torsion. On the other hand, $B/uB\to A/uA$ is surjective. This implies that $B/uB=A/uA$.

Let $\tau_\mu$ be the image of $\mu$ in
$\Omega\cong\xcoch(T)_\Ga/\xcoch(T_\s)_\Ga$ and let $Y^\circ\subset
\bold S$ so that $\sigma_{Y^\circ}=\tau_\mu^{-1}(\sigma_Y)$ as
before. Let $g\in G_1(F)$ be a lifting of $t_{-\mu}\tau_\mu\in
W_\aff$. Then since $\Fl^Y_w\in (\bGr_{\calG,\mu})_{\tilde{0}}$ for
$w\in\Adm^Y(\mu)$ (see Lemma \ref{easy}),
$g({^{Y^\circ}\Fl_{\s,w}^Y})\subset
(\Gr_{\calG_1,\leq\mu})_{\tilde{0}}$ for $w\in \Adm^Y(\mu)^\circ$. In
other words,
$\calA^{Y}(\mu)^\circ\subset(\Gr_{\calG_1,\leq\mu})_{\tilde{0}}$. 

Let $\calL$ be an ample line bundle on $\Fl_\s^Y$. Suppose that its certain tensor power
$\calL^n$ extends to a line bundle on $\Gr_{\calG_1}$ by Corollary \ref{line and central
charge}. Then we have
\[\dim\Ga((\Gr_{\calG_1,\leq\mu})_y,\calL_b^{nc(\calL)})=\dim\Ga((\Gr_{\calG_1,\leq\mu})_{\tilde{0}},\calL^n)\geq \dim\Ga(\calA^{Y}(\mu)^\circ,\calL^n)\]
by the flatness and the fact that $H^1({^Y\Fl_w^{Y'}},\calL)=0$ for
any Schubert variety ${^Y\Fl_w^{Y'}}$ and any ample line bundle
$\calL$. In addition, the equality holds if and only if
$\calA^{Y}(\mu)^\circ=(\Gr_{\calG_1,\leq\mu})_{\tilde{0}}$. Clearly,
for $y\neq\tilde{0}$, $(\Gr_{\calG_1,\leq\mu})_y=g\Gr_{\leq\mu}$ and
therefore
\[\Ga((\Gr_{\calG_1,\leq\mu})_y,\calL_b^{nc(\calL)})\cong\Ga(\Gr_{\leq\mu},\calL_b^{nc(\calL)})\]
by Corollary \ref{line and central charge}. Therefore, Theorem
\ref{MainII} implies Theorem \ref{Fibers}. Conversely, Theorem
\ref{Fibers} implies that the statement of Theorem \ref{MainII}
holds for $\calL^n,\calL^{2n},\ldots$. Therefore we have the equality of Euler characteristic $\chi(\Gr_{\leq\mu},\calL_b^{mc(\calL)})=\chi(\calA^{Y}(\mu)^\circ,\calL^m)$ for any $m$ as both are polynomial in $m$. 
But it is well-known that both $\calL_b^{mc(\calL)}$ and $\calL^m$ have  no higher cohomology. (Charateristic $p>0$ case follows from Frobenius splitting, and characteristic zero case follows from the semicontinuity, see \cite{Ma} for details.)
Therefore, the statement of Theorem \ref{MainII} also holds for $\calL$.

\medskip

To finish this section, let us mention the following observation. 
\begin{cor}\label{change char}
If Theorem \ref{MainII} (equivalently, Theorem \ref{Fibers}) holds for one prime $p\nmid e$, then it holds for
all $p\nmid e$ as well as in the case $\on{char}k=0$.
\end{cor}
\begin{proof}
Recall that the affine flag varieties and Schubert varieties are
defined over $W(k)$, the ring of Witt vectors of $k$, and the formation
commutes with base change (\cite{F,PR1}). In addition, line bundles are also defined over $W(k)$. (After identifying the affine flag varieties with those arising from Kac-Moody theory (\cite{PR1}), this follows from \cite[XVIII]{Ma}. In fact, they are even defined over $\bbZ'$, where $\bbZ$
is obtained from $\bbZ$ by adding $e$th roots of unity and inverting
$e$.) By the vanishing of corresponding $H^1$ (reason mentioned above),  both sides  are free
$W(k)$-modules and the formation of cohomology commutes with base change. 
The corollary follows.
\end{proof}

\section{Some properties of $\bGr_{\calG,\mu}$}\label{some prop}
In this section, we study two basic geometrical structures of
$\bGr_{\calG,\mu}$: (i) in \S \ref{affine chart}, we will construct
certain affine charts of $\bGr_{\calG,\mu}$, which turn out to be
isomorphic to affine spaces over $\tilde{C}$; and (ii) in \S \ref{Gm
action}, we will construct a $\bbG_m$-action on $\bGr_{\calG,\mu}$,
so that the map $\bGr_{\calG,\mu}\to \tilde{C}$ is
$\bbG_m$-equivariant, where $\bbG_m$ acts on $\tilde{C}=\bbA^1$ by
natural dilatation. To establish (i), we will need to first
construct the global root subgroups of $\calL\calG$ as in \S
\ref{globrootgrp}. We shall remark that the proofs of these results for $G$ split are very straightforward. It is only when $G$ is not split that some
complicated discussion is needed. Those who are only interested in
split groups can skip this section.

\subsection{global root groups}\label{globrootgrp}
We will introduce certain ``root subgroups" of $\calL\calG$ (more precisely, of $\calL^+\calG$, see Remark \ref{local root group} (i) ), whose
fibers over $0\in C$ is the usual root subgroups of the loop group
$LG$ as constructed in \cite[9.a,9.b]{PR1}.

Let us first review the shape of root groups of $G$. Let
$(H,B_H,T_H,X)$ be a pinned Chevalley group over $\bbZ$ as in \S
\ref{grp data}. In particular, $H_\der$ is simply-connected. Let $\Xi$ be the group of pinned automorphisms of $H_\der$,
which is simple, almost simple, simply-connected by our assumption.
So $\Xi$ is a cyclic group of order $1,2$ or $3$. Let
$\tilde{\Phi}=\Phi(H,T_H)$ be the set of roots of $H$ with respect
to $T_H$. For each $\tilde{a}\in\Phi(H,T_H)$, let
$\tilde{U}_{\tilde{a}}$ denote the corresponding root group. Then
for each $\ga\in\Xi$, one has an isomorphism
$\ga:\tilde{U}_{\tilde{a}}\cong\tilde{U}_{\ga\tilde{a}}$. The
stabilizer of $\tilde{a}$ in $\Xi$ is either trivial or the whole
group. Let us choose a Chevalley-Steinberg system of $H$, i.e. for
each $\tilde{a}\in\Phi(H,T_H)$, an isomorphism
$x_{\tilde{a}}:\bbG_a\cong \tilde{U}_{\tilde{a}}$ over $\bbZ$. In
addition, we require that:
\begin{enumerate}
\item if $\tilde{a}\in\Delta$ is a simple root, then
$X_{\tilde{a}}=dx_{\tilde{a}}(1)$, where
$X=\sum_{\tilde{a}\in\Delta}X_{\tilde{a}}$;
\item if the stabilizer
of $\tilde{a}$ in $\Xi$ is trivial, then $\ga\circ
x_{\tilde{a}}=x_{\ga\tilde{a}}$ for any $\ga\in\Xi$.
\end{enumerate}
Note that if $\gamma$ stabilizes
$\tilde{a}$, it is not necessarily always the case that $\ga\circ x_{\tilde{a}}=x_{\tilde{a}}$, as can be
seen for $\SL_3$. In this case, one obtains a quadratic character
\begin{equation}\label{quad}\chi_{\tilde{a}}:\Xi\to\Aut_\bbZ(\bbG_a)=\{\pm 1\}\end{equation}
such that $\ga\circ
x_{\tilde{a}}=x_{\tilde{a}}\circ\chi_{\tilde{a}}(\ga)$. Of course,
this can happen only if the order of $\Xi$ is $2$.

Recall that $\Ga=\Aut(\tilde{C}/C)$ is a group of order $e=1,2,3$,
which acts on $H$ via pinned automorphisms and the corresponding map
$\Ga\to\Xi$ is injective.

Let $j:\Phi(H,T_{H})\to \Phi(G,S)$ be the restriction of the root
systems. For $a\in \Phi(G,S)$, let
\[\eta(a)=\{\tilde{a}\in\Phi(H,T_{H})|j(\tilde{a})=ma, m\geq 0\}.\]
This is a subset of $\Phi(H,T_{H})$ satisfying the condition of \cite[5.1.16]{C}. Let $\tilde{U}_{\eta(a)}$ be the closed subgroup scheme of $H$ as defined in \emph{loc. cit.}. As a scheme,
$\tilde{U}_{\eta(a)}\cong\prod_{\tilde{a}\in\eta(a)}\tilde{U}_{\tilde{a}}$,
where the product is taken over any given order (which we fix from
now on) on $\eta(a)$. Informally, this is the subgroup of $H$ generated by
$\tilde{U}_{\tilde{a}}, \tilde{a}\in\eta(a)$. This subgroup is invariant under $\Xi$.  Then
$(\Res_{\tilde{F}/F}\tilde{U}_{\tilde{a}})^\Ga$ is the root group of
$G$ corresponding to $a$.

For an integer $n$, let us denote by $\bbG_{a,n,\tilde{C}}$ the
group scheme over $\tilde{C}$, which is the $n$th congruent group
scheme of $\bbG_{a,\tilde{C}}$. In other words,
$\bbG_{a,n+1,\tilde{C}}$ is the dilatation of $\bbG_{a,n,\tilde{C}}$
along the trivial subgroup in the fiber over $\tilde{0}$ (see
\cite[\S 3.2]{BLR} or \S \ref{deform}). More concretely,
$\bbG_{a,n,\tilde{C}}=\Spec k[u,t_n]\simeq\bbG_{a,\tilde{C}}$ and
the map $\bbG_{a,n+1,\tilde{C}}\to\bbG_{a,n,\tilde{C}}$ is given by
$t_n\mapsto ut_{n+1}$. We also have the congruent group schemes
$\tilde{U}_{\tilde{a},n,\tilde{C}}$ of
$\tilde{U}_{\tilde{a},\tilde{C}}$. The Chevalley-steinberg
isomorphism $x_{\tilde{a}}:\bbG_a\to\tilde{U}_{\tilde{a}}$ induces
the isomorphism
\[x_{\tilde{a},n}:\bbG_{a,\tilde{C}}\simeq\bbG_{a,n,\tilde{C}}\to
\tilde{U}_{\tilde{a},n,\tilde{C}}\] making the following diagram
commutative
\begin{equation}\label{root; comm}
\begin{CD}\bbG_{a,\tilde{C}}@>x_{\tilde{a},n+1}>>\tilde{U}_{\tilde{a},n+1,\tilde{C}}\\
@Vt_n\mapsto ut_{n+1}VV@VVV \\
\bbG_{a,\tilde{C}}@>x_{\tilde{a},n}>>\tilde{U}_{\tilde{a},n,\tilde{C}}
\end{CD}
\end{equation}

\medskip

Our goal is to construct some global root groups for $\calL\calG$.
For the purpose, we describe a construction of 
$\calG$.

Let us normalize the valuation so that $u$ has value $1/e$. Then we
embed $A(G,S)$ into $A(H,T_H)$. Let $x\in \sigma_Y$ be a point. It
determines a parahoric group scheme $\tilde{\calG}_x$ of $H\otimes\tilde{F}$, and $\calG_{\sigma_Y}$ is the neutral connected component of $(\Res_{\tilde{F}/F}\tilde{\calG}_x)^\Ga$. (One can see the claim as follow: Lifting $x$ to a point in the extended building of $G$, then $(\tilde{\calG}_x(\tilde{\calO}))^\Gamma\subset G(F)$ is the stabilizer of this point. On the other hand, by \cite[2.2, 3.4]{Ed} $(\Res_{\tilde{F}/F}\tilde{\calG}_x)^\Ga$ is smooth. Therefore, its neutral connected component is the parahoric group scheme of $G$ given $x$.)

We extend $\tilde{\calG}_x$ to a
group scheme $\tilde{\calG}$ over $\tilde{C}$ as in \S \ref{grpsch},  so that $\tilde{\calG}|_{\tilde{C}^\circ}=H\times\tilde{C}^\circ$ and $\tilde{\calG}|_{\calO_0}=\tilde{\calG}_x$ (under the identification $\tilde{F}\cong\tilde{F}_0$). 
From the
construction, $\tilde{\calG}$ contains
\[\prod_{\tilde{a}\in\Phi(H,T_H)^-}\tilde{U}_{\tilde{a},\lceil e\tilde{a}(v_0-x)\rceil,\tilde{C}}\times T_{H,\tilde{C}}\times \prod_{\tilde{a}\in\Phi(H,T_H)^+}\tilde{U}_{\tilde{a},\lceil e\tilde{a}(v_0-x)\rceil,\tilde{C}}\] as
a fiberwise dense open subscheme (\cite[2.2.10, 3.9.4]{BT2}), where $\lceil y
\rceil$ denotes the smallest integer that are $\geq y$. Observe that
since $x$ is fixed under the action of $\Gamma$, for $a\in\Phi$, the
closed subgroup scheme
$\prod_{\tilde{a}\in\eta(a)}\tilde{U}_{\tilde{a},\lceil e \tilde{a}(x-v_0)\rceil,\tilde{C}}$
of $\tilde{\calG}$ is invariant under the action of $\Ga$. Let
\[U_{a,\sigma_Y,C}=(\Res_{\tilde{C}/C}\prod_{\tilde{a}\in\eta(a)}\tilde{U}_{\tilde{a},\lceil e\tilde{a}(v_0-x)\rceil,\tilde{C}})^{\Ga},\]
which does not depend on $x$. By \cite[2.2, 3.4]{Ed}, $U_{a,\sigma_Y,C}$ is smooth. In addition, a check for $\SL_2$ and $\SU_3$ cases shows that $U_{a,\sigma_Y,C}$ is connected.
Then
$(U_{a,\sigma_Y,C})_{F_0}$ is the root group of $\calG_{F_0}\cong
G$ corresponding to $a$, and for $y\neq 0$,
$(U_{a,\sigma_Y,C})_y\cong\tilde{U}_{\eta(a)}$ non-canonically. In addition,
\[\prod_{a\in\Phi^{nd,-}}U_{a,\sigma_Y,C}\times(\Res_{\tilde{C}/C}T_{H,\tilde{C}})^{\Ga,0}\times\prod_{a\in\Phi^{nd,+}}U_{a,\sigma_Y,C}\]
is a fiberwise dense open subscheme of $\calG$, where
$\Phi^{nd}\subset\Phi=\Phi(G,S)$ denote the set of non-divisible
roots, i.e. $a\in\Phi^{nd}$ if $a/2\not\in\Phi$. Given an affine root $\al$ of $G$ with vector part $a$, the corresponding root subgroup of $\calL\calG$ will be constructed as a closed subgroup scheme of $\calL
U_{a,\sigma_Y,C}$.

Recall that we constructed the special vertex $v_0$ in \S \ref{grp
data}. In \S \ref{grp data}, we use this vertex to identify $A(G,S)$
with $\xcoch(S)_\bbR$. Then we can write affine roots as $a+m$,
where $a\in\Phi(G,S)$ and $m\in\frac{1}{e}\bbZ$. Let $a+m$ be an
affine root such that $em\geq \lceil ea(v_0-x)\rceil$. Let us
construct a closed immersion
\begin{equation}\label{root subgroup}
x_{a+m}:\bbG_{a,C}\to\calL U_{a,\sigma_Y,C}.
\end{equation}
Let us describe of $x_{a+m}$ at the level of $R$-points, where $R$ is a $k$-algebra. Recall we
write $C=\Spec k[v], \tilde{C}=\Spec k[u]$, such that
$[e]:\tilde{C}\to C$ is given by $v\mapsto u^e$. Let $y:\Spec R\to C$ be an $R$-point of $C$. We
identify $\Hom_{C}(\Spec R,\bbG_{a,C})$ with $R$ in an obvious
manner. We thus need to construct a map (functorial in
$R$)
\[x_{a+m}:R\to \Hom_C(\Spec R,\calL U_{a,\sigma_Y,C}).\]

The graph of $y:\Spec R\to C$ is $\Ga_{y}=\Spec R[v]/(v-y)$ and
$\hat{\Ga}^\circ_{y}=\Spec R(\!(v-y)\!)$. Now, by definition
\begin{equation}\label{pts}\Hom_{C}(\Spec R,\calL U_{a,\sigma_Y,C})=\Hom(\Spec
R(\!(v-y)\!)\times_{C}\tilde{C},\tilde{U}_{a,\sigma_Y,\tilde{C}})^\Ga,\end{equation} where $\Ga$
acts on $\Spec R(\!(v-y)\!)\times_{C}\tilde{C}$ via the action on
$\tilde{C}$, and acts on $\tilde{U}_{a,\sigma_Y,\tilde{C}}:=\prod_{\tilde{a}\in\eta(a)}\tilde{U}_{\tilde{a},\lceil
e\tilde{a}(v_0-x)\rceil,\tilde{C}}$ as above.

Let us introduce the following notation. Each element
$s\in R(\!(v-y)\!)\otimes_{k[v]}k[u]$
determines a morphism $\Spec R(\!(v-y)\!)\times_{C}\tilde{C}\to
\bbG_{a,\tilde{C}}$, and let
\[x_{\tilde{a},n}(s):\Spec R(\!(v-y)\!)\times_{C}\tilde{C}\to \tilde{U}_{\tilde{a},n,\tilde{C}}\]
denote the composition of this morphism with
$x_{\tilde{a},n}:\bbG_{a,\tilde{C}}\to\tilde{U}_{\tilde{a},n,\tilde{C}}$.

Now we construct $x_{a+m}$. There are two cases.

\begin{enumerate}
\item[(i)] $2a\not\in\Phi(G,S)$. In this case, $\Ga$ acts
transitively on $\eta(a)$. There are two subcases.

\item[(ia)] $\eta(a)=\tilde{a}$, so that $\Ga$ fixes $\tilde{a}$ and $\tilde{U}_{\eta(a)}=\tilde{U}_{\tilde{a}}$.
Define
\[x_{a+m}(r)=x_{\tilde{a},\lceil ea(v_0-x)\rceil}(r\otimes u^{em-\lceil ea(v_0-x)\rceil}).\]
Since $a+m$ is an affine root, $\Ga$ acts on $u^{em-\lceil
ea(v_0-x)\rceil}$ exactly via the quadratic character
$\chi_{\tilde{a}}$ as defined in \eqref{quad}, $x_{a+m}(r)$ is an
element in \eqref{pts}.

\item[(ib)] $\Ga$ acts simply transitively on $\eta(a)$. Choose $\tilde{a}\in\eta(a)$ and
$\ga\in\Ga$ a generator. Using the isomorphism
$\prod_{i=1}^{e}\tilde{U}_{\ga^i(\tilde{a})}\cong\tilde{U}_{\eta(a)}$,
one defines
\[x_{a+m}(r)=\prod_{i=1}^{e}x_{\ga^i(\tilde{a}),\lceil ea(v_0-x)\rceil}(r\otimes\ga^i(u)^{em-\lceil ea(v_0-x)\rceil}).\]
Since for $\tilde{a},\tilde{a}'\in\eta(a)$, the groups
$\tilde{U}_{\tilde{a}}$ and $\tilde{U}_{\tilde{a}'}$ commute, and
therefore $x_{a+m}(r)$ is an element in \eqref{pts}.

\item[(ii)] $2a\in \Phi(G,S)$, so that
$\eta(a)=\{\tilde{a},\tilde{a}',\tilde{a}+\tilde{a}'\}$. In this
case, $\on{char}k\neq 2$, $e=2$, and the group is the odd unitary
group. In addition, the quadratic character
$\chi_{\tilde{a}+\tilde{a}'}$ is non-trivial. Recall that for any
$s,s'$,
\begin{equation}\label{comm}x_{\tilde{a}}(s)x_{\tilde{a}'}(s')=x_{\tilde{a}'}(s')x_{\tilde{a}}(s)x_{\tilde{a}+\tilde{a}'}(\pm
ss'),\end{equation} where $\pm$ depends on
$x_{\tilde{a}},x_{\tilde{a}'},x_{\tilde{a}+\tilde{a}'}$, but not on
$s,s'$. Define
\begin{eqnarray*}x_{a+m}(r)&=&x_{\tilde{a},\lceil e\tilde{a}(v_0-x)\rceil}(r\otimes
u^{em-\lceil e\tilde{a}(v_0-x)\rceil})\times
\\
&&x_{\tilde{a}',\lceil e\tilde{a}(v_0-x)\rceil}((-1)^{em-\lceil
e\tilde{a}(v_0-x)\rceil}r\otimes
u^{em-\lceil e\tilde{a}(v_0-x)\rceil})\times\\
&&x_{\tilde{a}+\tilde{a}',\lceil
2e\tilde{a}(v_0-x)\rceil}(\mp(-1)^{em-\lceil
e\tilde{a}(v_0-x)\rceil}\frac{1}{2}r^2\otimes u^{2em-\lceil
2e\tilde{a}(v_0-x)\rceil)})
\end{eqnarray*}
where $\mp$ is the sign \bf opposite \rm the sign $\pm$ in
\eqref{comm}. Using \eqref{comm}, it is clear that $x_{a+m}(r)$ is
again an element in \eqref{pts}.
\end{enumerate}

We have completed the construction of \eqref{root subgroup}. Note that they are independent of the choice of $x\in\sigma_Y$ by \eqref{root; comm}.
In addition, over $0\in C$ (i.e. by setting $y=0$), the map \eqref{root
subgroup} reduces to an isomorphism of $\bbG_a$ and the root
subgroup of $LG$ corresponding to $a+m$, as constructed in
\cite[9.a,9.b]{PR1}. This motivates us to define
\begin{dfn}
Let $a+m$ be an affine root of $G$ such that $em\geq \lceil
ea(v_0-x)\rceil$. The subgroup scheme
$\calU_{a+m}=x_{a+m}(\bbG_{a,C})$ is called root subgroup of
$\calL\calG$ corresponding to $a+m$.
\end{dfn}

\begin{rmk}\label{local root group}
(i) Note that in the above definition, the requirement $em\geq\lceil ea(v_0-x)\rceil$ is
necessary, as we need $r\otimes u^{em-\lceil ea(v_0-x)\rceil}$ to be an element in $R(\!(v-y)\!)\otimes_{k[v]}k[u]$. Note that in fact $\calU_{a+m}\subset\calL^+U_{a,\sigma_Y,C}$. If
$f:\calG'\to\calG$ is a map of Bruhat-Tits group schemes, then
$\calL f x_{a+m}=x_{a+m}$ if $x_{a+m}$ is defined for $\calG'$ (and
therefore for $\calG$).

(ii) By taking the fibers
$U_{a+m}=(\calU_{a+m})_0\subset(\calL\calG)_0\cong LG$, we obtain
the root subgroups of $LG$. Note that, however, as $R(\!(v)\!)\otimes_{k[v]}k[u]=R(\!(u)\!)$, we could drop the requirement
$em\geq\lceil e\tilde{a}(v_0-x)\rceil$ and $U_{a+m}\subset LG$ is defined for
\emph{all} affine roots of $G$ . If we do not identify $A(G,S)$ with
$\xcoch(S)_\bbR$ via $v_0$, we write them as $U_\al$, where $\al$ is
an affine root.
\end{rmk}

The following lemma about the root subgroups for (global)
loop groups is the counterpart of a well-known fact about
the root subgroups of Kac-Moody groups. To describe it, let us use the following notation. for a group (ind)-scheme $\calU$ over $C$ and $y\in C(R)$ and $R$-point, $\calU(R)$ will denote the group of $R$-points of $\calU$ over $y$.

\begin{lem}\label{commutator}
Let $R$ be a $k$-algebra and $y\in C(R)$. Let $a+m,b+n$ ($a\not\in\bbR b$) be two affine roots of $G$ such that $\calU_{a+m},\calU_{b+n}$ are defined. 
Then the commutator
$[\calU_{a+m}(R),\calU_{b+n}(R)]$ is contained in the group generated by
$\calU_{(pa+qb)+(pm+qn)}(R)$, where $p,q\in\bbZ_{>0}$ such that
$(pa+qb)+(pm+qn)$ is also an affine root of $G$ (the groups
$\calU_{(pa+qb)+(pm+qn)}$ are clearly defined for $\calG$).
\end{lem}
\begin{proof}
Let us define a subset $\Psi_{a,b}\subset\tilde{\Phi}=\Phi(H,T_H)$ as
\[\Psi_{a,b}=\{\tilde{a}\in\tilde{\Phi}\mid j(\tilde{a})=pa+qb \mbox{
for }
p,q\in\bbZ_{>0}\}=\bigcup_{pa+qb\in\Phi^{nd},p,q>0}\eta(pa+qb),\]
where $j:\tilde{\Phi}\to\Phi$. For $\tilde{a}\in\Psi_{a,b}$ such
that $j(\tilde{a})=pa+qb$, let $k(\tilde{a})=pm+qn$. Using the same
notation as above, let us define
\[\tilde{\calU}_{\tilde{a}+k(\tilde{a})}\subset\calL\Res_{\tilde{C}/C}\tilde{U}_{a,\sigma_Y,\tilde{C}},\]
where $\tilde{U}_{a,\sigma_Y,\tilde{C}}=\tilde{U}_{\tilde{a},\lceil e\tilde{a}(v_0-x)\rceil,\tilde{C}})$, to
be the group over $C$, whose $R$-points over $y:\Spec R\to C$ are
given by
\[\{x_{\tilde{a},\lceil e\tilde{a}(v_0-x)\rceil}(r\otimes u^{ek(\tilde{a})-\lceil e\tilde{a}(v_0-x)\rceil}), r\in R\}\subset\Hom_C(\Spec R, \calL(\Res_{\tilde{C}/C}\tilde{U}_{a,\sigma_Y,\tilde{C}})).\]
Let $p,q\in\bbZ_{>0}$, and let $\tilde{\calU}_{\eta(pa+qb),pm+qn}$
be the group generated by $\tilde{\calU}_{\tilde{a}+k(\tilde{a})}$
for those $\tilde{a}\in\eta(pa+qb)\subset\Psi_{a,b}$. Then over
$y\in C(R)$,
\[\calU_{(pa+qb)+(pm+qn)}(R)=\tilde{\calU}_{\eta(pa+qb),pm+qn}(R)\cap\calL U_{pa+qb,\sigma_Y,C}(R).\]
Let
$\tilde{\calU}_{\Psi_{a,b},m,n}$ be the group generated by
$\tilde{\calU}_{\tilde{a}+k(\tilde{a})}, \tilde{a}\in\Psi_{a,b}$.
Recall that for the fixed Chevalley-Steiberg system
$\{x_{\tilde{a}}, \tilde{a}\in\tilde{\Phi}\}$, and for two roots
$\tilde{a},\tilde{b}\in\tilde{\Phi}$, there exists $c(p,q)\in\bbZ$
for $p,q\in\bbZ_{>0}$ such that for any ring $R$ and $r,s\in R$, the
commutator $[x_{\tilde{a}}(r),x_{\tilde{b}}(s) ]$ can be written as
$[x_{\tilde{a}}(r),x_{\tilde{b}}(s) ]=\prod_{p\tilde{a}+q\tilde{b}\in\tilde{\Phi},p,q>0}x_{p\tilde{a}+q\tilde{b}}(c(p,q)r^ps^q)$ (\cite[Proposition 5.1.14]{C}).
Therefore, the commutator of
$[\tilde{\calU}_{\tilde{a}+k(\tilde{a})},\tilde{\calU}_{\tilde{b}+k(\tilde{b})} ]$
is contained in the group generated by
$\tilde{\calU}_{p\tilde{a}+q\tilde{b}+k(p\tilde{a}+q\tilde{b})}$,
where $p,q\in\bbZ_{>0}$ and $p\tilde{a}+q\tilde{b}\in\tilde{\Phi}$.
Now we can apply \cite[Proposition 6.1.6]{BT}, with the pair
$\tilde{\calU}_{\tilde{a}+k(\tilde{a})}(R)\subset \tilde{U}_{\tilde{a}}$ playing the role $Y_a\subset U_a$ in \emph{loc. cit.}. Then  we have
\[\tilde{\calU}_{\Psi_{a,b},m,n}(R)\stackrel{\simeq}{\leftarrow}\prod_{\tilde{a}\in\Psi_{a,b}}\tilde{\calU}_{\tilde{a}+k(\tilde{a})}(R)\stackrel{\simeq}{\rightarrow}\prod_{pa+qb\in\Phi^{nd},p,q>0}\tilde{\calU}_{\eta(pa+qb),pm+qn}(R).\]
Here the first isomorphism is obtained by setting $\Psi^{\rm red}$ in \emph{loc. cit.} as $\Psi_{a,b}$, and the second isomorphism is obtained by setting $\Psi^{\rm red}=\eta(pa+qb)$ for all $pa+qb\in \Phi^{nd}, p, q>0$.

Next, let $\calL U_{(a,b)}$ be the group generated by $\calL
U_{pa+qb,\sigma_Y,C}$, $pa+qb\in\Phi, p,q\in\bbZ_{>0}$. Again by
\emph{loc.cit.}, for $y\in C(R)$, there is a bijection
\[\prod_{pa+qb\in\Phi^{nd}, p,q>0}\calL U_{pa+qb,\sigma_Y,C}(R)\stackrel{\simeq}{\rightarrow} \calL U_{(a,b)}(R).\]
Combining the above two isomorphisms, we thus obtain that
\[\begin{split}&(\calU_{a+m}(R),\calU_{b+n}(R))\subset\tilde{\calU}_{\Psi_{a,b},m,n}\cap\calL U_{(a,b)}(R)\\
=&\prod_{pa+qb\in\Phi^{nd},p,q>0}(\tilde{\calU}_{\eta(pa+qb),pm+qn}(R)\cap\calL
U_{pa+qb,\sigma_Y,C}(R))\\
=&\prod_{pa+qb\in\Phi^{nd},p,q>0}\calU_{(pa+qb)+(pm+qn)}(R).\end{split}\]
\end{proof}

\subsection{Some affine charts of $\bGr_{\calG,\mu}$}\label{affine
chart} We introduce certain affine charts of $\bGr_{\calG,\mu}$,
which turn out to be isomorphic to affine spaces. Let
$\Lambda=W_0\mu\subset\xcoch(T)_\Ga$ as before, and let
$\la\in\Lambda$. Denote $\Phi_\la\subset\Phi(G,S)$ to be the subset
\begin{equation}\label{Phila}\Phi_\la=\{a+m\mid (a,\la)>0, a+m \mbox{ affine root }, 0\leq m-\frac{\lceil
ea(v_0-x)\rceil}{e}<(a,\la)\}.\end{equation} 
By Lemma \ref{length}, this is a set with $(2\rho,\mu)$ elements (recall that $\mu\in\xcoch(T)^+$).
For each $a+m$, $\calU_{a+m}$ is defined and is a subgroup of $\calL^+\calG$. 

Let us endow a total order on the set $\Phi_\la$ as follows: First fix an order on $\{a\mid
(a,\la)>0\}\cap\Phi^{nd}$. Then we can extend it to an order on $\{a\mid
(a,\la)>0\}$ by requiring if $a,2a\in\{a\mid (a,\la)>0\}$, then
$a<2a<b$ for any $b\in\{a\mid (a,\la)>0\}\cap\Phi^{nd}$ such that
$a<b$. Finally, we can give an order on $\Phi_\la$ by requiring
$a+m<b+n$ if either $a<b$ or $a=b,m<n$.

Now, consider $\prod_{\Phi_\la}\calU_{a+m}\to\calL^+\calG$ given by multiplication. Here and everywhere else the fiber products are over $C$. This is a closed immersion. In fact, let $\Psi\subset \Phi(G,S)$ be the image of the map $\Phi_\la\to \Phi(G,S)$ by taking the vector part of an affine root. Then $\Psi\cap(-\Psi)=\emptyset$. Therefore $\prod_{a\in\Psi}U_{a,\sigma_Y,C}\to \calG$ is a closed embedding. On the other hand, for $a\in\Psi$, it is not hard to see that the morphism $\prod_{m}\calU_{a+m}\prod_{m'}\calU_{2a+m'}\to \calL^+ U_{a,\sigma_Y,C}$ is a closed immersion. The claim follows.

Let us denote by
$\calU_{\Phi_\la}\subset\calL^+\calG$ the image of the above map. This is a closed subscheme of $\calL^+\calG$. We claim that $\calU_{\Phi_\la}$ is indeed a closed subgroup scheme of $\calL^+\calG$.

\begin{lem}Let $R$ be a $k$-algebra  and $y\in C(R)$. Then the $R$-points of $\calU_{\Phi_\la}$ over $y$ form the subgroup of $\calL^+\calG(R)$ generated by $\calU_{a+m}(R), a+m\in\Phi_\la$. 
\end{lem}
\begin{proof}Let us denote the subgroup generated by $\calU_{a+m}(R), a+m\in\Phi_\la$ by $\langle \calU_{a+m}(R) \rangle$.
By Lemma \ref{commutator}, the collection of groups
$\{\calU_{a+m}(R), a+m\in\Phi_\la\}$ satisfies the condition as
required by \cite[Lemma 6.1.7]{BT}. Then by \emph{loc. cit.}, we have
\[\calU_{\Phi_\la}(R)=\prod_{a+m\in\Phi_\la}\calU_{a+m}(R)\cong \langle \calU_{a+m}(R) \rangle.\]
The lemma follows.
\end{proof}

Recall the section $s_\la:\tilde{C}\to\wLG$ as constructed in the
paragraph after Proposition \ref{sla}\footnote{More precisely, we
need to choose a lifting of $\tilde{\la}\in\xcoch(T)$ of $\la$, but
to simply the notation, we denote this lifting by $\la$.}. Consider
\[c_\la: \calU_{\Phi_\la}\times_C\tilde{C}\to\bGr_{\calG,\mu}, \quad g\mapsto gs_\la.\]
\begin{prop}\label{cla}
The morphism $c_\la$ is an open immersion.
\end{prop}
\begin{proof}We first show that the 
stabilizer of $s_\la$ in $\calU_{\Phi_\la}$ is trivial. 
Recall that
$\calL\calG$ acts on $\Gr_\calG$, and the stabilizer of the section
$e:C\to\Gr_\calG$ (defined by the trivial $\calG$-torsor) is
$\calL^+\calG$. Therefore, the stabilizer in $\wLG$ of the section
$s_\la$ is $s_\la(\wJG) s_\la^{-1}$. Therefore, it is enough to
prove that $\wJG\cap
s_\la^{-1}(\calU_{\Phi_\la}\times_C\tilde{C})s_\la$ is trivial, or
equivalently \[(\calL^+ U_{a,\sigma_Y,C}\times_C\tilde{C})\cap
s_\la^{-1}(\calU_{a+m}\times_C\tilde{C})s_\la\] is trivial for all
$a+m\in\Phi_\la$.

Let us analyze the $R$-points of
$s_\la^{-1}(\calU_{a+m}\times_C\tilde{C})s_\la$ over $y:\Spec
R\to\tilde{C}$. Recall that $s_\la(y)$ is given by the
$\Ga$-equivariant map
\[s_\la(y):\Spec R(\!(v-y^e)\!)\otimes_{k[v]}k[u]\to T_H\]
such that for any weight $\omega$ of $T_H$, the composition $\omega
s_\la(y)$ (which is determined by an invertible element in
$R(\!(v-y^e)\!)\otimes_{k[v]}k[u]$) is $\prod_{i=1}^e(1\otimes
\ga^i(u)-y\otimes 1)^{(\la,\ga^i\omega)}$. Note that for any
$\tilde{a}\in\Phi(H,T_H)$ such that $j(\tilde{a})=a\in\{a\mid
(a,\la)>0\}$,
\[
\prod_{i=1}^e(1\otimes\ga^i(u)-y\otimes
1)^{(-\la,\ga^i(\tilde{a}))}u^{em-\lceil
e\tilde{a}(v_0-x)\rceil}\not\in R[[v-y^e]]\otimes_{k[v]}k[u],
\]
as $em-\lceil
e\tilde{a}(v_0-x)\rceil<e(\la,a)$, which implies $\calL^+
U_{a,\sigma_Y,C}\times_C\tilde{C}\cap
s_\la^{-1}(\calU_{a+m}\times_C\tilde{C})s_\la$ is the identity for
all $a+m\in\Phi_\la$. 

Therefore, the stabilizer of $s_\la$ in $\calU_{\Phi_\la}$ is trivial. Then, $c_\la$ is a monomorphism of irreducible varieties over $k$ of the same dimension. We show that $U=c_\la(\calU_{\Phi_\la}\times_C\tilde{C})$ is open. 
For simplicity, we write $\bGr_{\calG,\mu}\to\tilde{C}$ as $f:X\to \tilde{C}$, and write the group scheme $\calU_{\Phi_\la}\times_C\tilde{C}$ by $\calU$.
As $\dim\calU=\dim X=(2\rho,\mu)+1$, and  $U$ is constructible, it contains an non-empty open subset of $X$. Let $W\subset U$ be the maximal open subset of $X$ contained in $U$. Then  $W$ is $\calU$-stable. In particular, $W=\calU(W\cap s_\la(\tilde{C})$. We claim that $s_\la(\tilde{C})\subset W$, which implies that $W=U$. Otherwise, $\tilde{C}- f(W\cap s_\la(\tilde{C}))$ consist of finitely many points $x_1,\ldots,x_n$. Then $W\cap f^{-1}(x_i)=\emptyset$. Note that $U_{x_i}:=f^{-1}(x_i)\cap U$ is just the orbit of $s_\la(x_i))$ under $\calU_{x_i}$ in $f^{-1}(x_i)$. As itself contains a non-empty open subset of $f^{-1}(x_i)$, $U_{x_i}$ is open in $f^{-1}(x_i)$. Let $Z_{x_i}=f^{-1}(x_i)-U_{x_i}$. This is a closed subset of $f^{-1}(x_i)$.

Let $Y=f^{-1}(f(W\cap s_\la(\tilde{C})))$. Then $W$ is open dense in $Y$. Let $D=Y-W$, which is a proper closed subset of $Y$, and let $\bar{D}$ be its closure in $X$. Then $\bar{D}$ flat over $\tilde{C}$. Therefore $\bar{D}_{x_i}:=f^{-1}(x_i)\cap\bar{D}$ is a proper closed subset of $f^{-1}(x_i)$, of dimension strictly smaller than $(2\rho,\mu)$. Therefore $U_{x_i}\not\subset \bar{D}_{x_i}$. Now, $X-\bar{D}-\bigcup_i Z_{x_i}$ is open, contained in $U$, and is strictly larger than $W$. This is a contradiction.

We therefore have proved that $c_\la$ is a monomorphism, which maps onto an open subset of $\bGr_{\calG,\mu}$. Finally we show that $c_\la:\calU\to \bGr_{\calG,\mu}$ is an open immersion. It is enough to show that $c_\la: \calU\to \wGr_{\calG}$ is a locally closed embedding. Let $V=\wGr_\calG-(\bGr_{\calG,\mu}-U)$, considered as an open sub-indscheme of $\wGr_\calG$. Then it is enough to show that $c_\la:\calU\to V$ is a closed immersion. But this can be checked fiberwise over $\tilde{C}$: Over a point $x\in\tilde{C}$, orbit maps are alway locally closed immersions. As the image of $c_\la$ is closed in $V$, over each point $x\in\tilde{C}$, $(c_\la)_x:\calU_x\to V_x$ is a closed immersion. The proposition follows.
\end{proof}

In what follows, we denote the image of $c_\la$ ($\la\in\Lambda$) by
$U_\la$, so that $U_\la$ is affine open in $\bGr_{\calG,\mu}$ which
is smooth over $\tilde{C}$ (indeed an affine space over
$\tilde{C}$). Note that
\[(U_\la)_{\tilde{0}}=(\calU_{\Phi_\la})_{\tilde{0}}s_\la(0)=L^+\calG_{\mathbf{a}}t_\la\]
is exactly the $L^+\calG_{\mathbf{a}}$-orbit through $t_\la$ in $\Fl^Y$.

\subsection{A $\bbG_m$-action on $\bGr_{\calG,\mu}$}\label{Gm
action} Let $\calG$ be a group scheme over $C$ as in \S \ref{grpsch}.
Let $f:\wGr_{\calG}\to \tilde{C}$ be the structural map. We
construct a natural $\bbG_m$-action on $\wGr_{\calG}$, which lifts
the natural action of $\bbG_m$ on $\tilde{C}$ via dilatations. In
addition, each $\bGr_{\calG,\mu}$ is stable under this
$\bbG_m$-action.

The construction of the $\bbG_m$-action on $\wGr_{\calG}$ is
straightforward. Recall that the global coordinate on $\tilde{C}$ is
$u$ and on $C$ is $v$, and that the map $[e]:\tilde{C}\to C$ is
given by $v\mapsto u^e$. Recall that an $R$-point of $\wGr_{\calG}$
is given by $u\mapsto y$ and a $\calG$-torsor $\calE$ on $C_R$,
which is trivialized over $C_R-\Ga_{[e](y)}$. Let $r\in R^\times$ be
an $R$-point of $\bbG_m$. We need to construct a new $\calG$-torsor
on $C_R$, together with a trivialization over $C_R-\Ga_{[e](ry)}$.
Indeed, let $r^e:C_R\to C_R$ given by $v\mapsto r^ev$. It maps
$\Ga_{[e](y)}$ to $\Ga_{[e](ry)}$. Then the pullback of $\calE$
along $r^{-e}$ is an $(r^{-e})^*\calG$-torsor on $C_R$, together
with a trivialization on $C_R-\Ga_{[e](ry)}$. Therefore, to complete
the construction, it is enough to show that $(r^{-e})^*\calG$ is
canonically isomorphic to $\calG$ as group schemes over $C_R$. Let
us remark that the same construction will give an action of $\bbG_m$ on
$\wLG$ (resp. $\wJG$), compatible with the dilatations on $\tilde{C}$.
Furthermore, the action of $\wLG$ (resp. $\wJG$) on $\wGr_{\calG}$
is $\bbG_m$-equivariant.

Let us define the action of $\bbG_m$ on $C=\Spec k[v]$ via
$(r,v)\mapsto r^ev$. Observe that $\mu_e\subset\bbG_m$ acts
trivially on $C$ via this action.

\begin{lem}Given the action of $\bbG_m$ on $C$ as above, there is a natural action of
$\bbG_m$ on $\calG$, such that $\calG\to C$ is $\bbG_m$-equivariant.
\end{lem}
\begin{rmk}However, the natural dilatation on $C$ could not lift to
$\calG$.
\end{rmk}
\begin{proof}As has been explained in \S \ref{globrootgrp}, there is a group scheme $\tilde{G}$ over $\tilde{C}$, satisfying $\tilde{G}|_{\tilde{C}^\circ}=H\times\tilde{C}^\circ$ and $\tilde{\calG}_{\tilde{\calO}_0}$ is a parahoric group scheme of $H\otimes\tilde{F}_0$, given by a point $x\in A(H,T_H)^\Gamma$, such that
$\calG$ is the neutral connected component of
$(\Res_{\tilde{C}/C}\tilde{\calG})^{\Ga}$. To prove the proposition,
it is enough to prove that there is a natural $\bbG_m$ action on
$\tilde{\calG}$, compatible with the rotation on $\tilde{C}$. In
addition, this $\bbG_m$-action should be compatible with the action of $\Ga$ on
$\tilde{\calG}$.

Let $m,p:\bbG_m\times\tilde{C}\to\tilde{C}$ be the action map and
the projection map respectively. We need show that there is an
isomorphism of group schemes $p^*\tilde{\calG}\cong
m^*\tilde{\calG}$ over $\bbG_m\times\tilde{C}$, satisfying the
usual compatibility conditions. Since $\bbG_m$ naturally acts on
$\tilde{\calG}|_{\tilde{C}^{\circ}}=H\times\tilde{C}^{\circ}$
by acting via rotation on the second factor, there is a natural isomorphism
\[c:p^*\tilde{\calG}|_{\bbG_m\times\tilde{C}^{\circ}}\cong m^*\tilde{\calG}|_{\bbG_m\times\tilde{C}^{\circ}},\]
which is compatible with the $\Ga$-actions. We need to show that this
uniquely extends to an isomorphism over $\bbG_m\times\tilde{C}$.
Then it will be automatically compatible with the $\Ga$-actions. Indeed, the uniqueness is clear since
$p^*\tilde{\calG}$ (resp. $m^*\tilde{\calG}$) is flat over
$\bbG_m\times\tilde{C}$, so that
$\calO_{p^*\tilde{\calG}}\subset\calO_{p^*\tilde{\calG}}[u^{-1}]$
(resp.
$\calO_{m^*\tilde{\calG}}\subset\calO_{m^*\tilde{\calG}}[u^{-1}]$).
We need to prove that the map
$c:\calO_{m^*\tilde{\calG}}[u^{-1}]\to\calO_{p^*\tilde{\calG}}[u^{-1}]$
indeed sends $\calO_{m^*\tilde{\calG}}$ to $\calO_{p^*\tilde{\calG}}$.
But this can be checked over each closed point of $\bbG_m$.
Therefore, it remains to prove that for every $r\in\bbG_m(k)$, the
isomorphism of
$r^*\tilde{\calG}|_{\tilde{C}^{\circ}}\cong\tilde{\calG}|_{\tilde{C}^{\circ}}$
extends to an isomorphism $r^*\tilde{\calG}\cong\tilde{\calG}$. We
can replace $\tilde{C}$ by $\calO_{\tilde{0}}$. By \cite[Proposition
1.7.6]{BT2}, it is enough to prove that the isomorphism
$r:\tilde{\calG}(F_{\tilde{0}})\to\tilde{\calG}(F_{\tilde{0}})$
induces an isomorphism
$\tilde{\calG}(\calO_{\tilde{0}})\to\tilde{\calG}(\calO_{\tilde{0}})$.
But it is clear that each root subgroup of $L(H\otimes\tilde{F}_0)$ with respect to $(H\otimes \tilde{F}_0, T_H\otimes \tilde{F}_0)$ as constructed in \S \ref{globrootgrp} (see Remark \ref{local root group} (ii)) is invariant under this $\bbG_m$-action. Therefore, for any $x\in A(H,T_H)$, the corresponding parahoric group of $H\otimes \tilde{F}_0$ is invariant under this $\bbG_m$-action.
\end{proof}

It remains to show that each $\bGr_{\calG,\mu}$ is invariant under
this $\bbG_m$-action. It is enough to show that the section
$s_\mu:\tilde{C}^{\circ}\to\wGr_{\calT}\subset\wGr_{\calG}$ is invariant under this
$\bbG_m$-action. Recall that $s_\mu:\tilde{C}^{\circ}\to
\Gr_\calT\times_C\tilde{C}^{\circ}\cong\Gr_{T_H\times\tilde{C}^{\circ}}$
is given by the $T_H$-bundle
$\calO_{\tilde{C}^{\circ}}(\mu\Delta)$ with its canonical
trivialization away from $\Delta$ (see Lemma \ref{moduli
interpretation of sla}).  From this moduli interpretation, it is
clear that $s_\mu$ is $\bbG_m$-invariant. 

\medskip

By restriction to $(\wGr_\calG)_{\tilde{0}}\cong\Fl^Y$, we obtain an
action of $\bbG_m$ on $\Fl^Y$ (and therefore on $\Fl^Y_{\s}$). As is shown in \cite{PR1}, the affine
flag variety $\Fl^Y_{\s}$ coincides with the affine flag variety in the
Kac-Moody setting. Under this identification, the above
$\bbG_m$-action on $\Fl^Y_{\s}$ corresponds to the action of the extra
one-dimensional torus (usually called the rotation torus) in the
maximal torus of the affine Kac-Moody group. We do not make the
statement precise. Instead, we mention

\begin{lem}\label{Gm stable}
Each Schubert variety in $\Fl^Y$ is invariant under this action of
$\bbG_m$ on $\Fl^Y$.
\end{lem}
\begin{proof}
Since $\bbG_m$ acts on $\calG$, it acts on $L^+\calG_{\calO_0}$.
Clearly, it also acts on $L^+\calT_{\calO_0}$, and therefore acts on
the normalizer $N_{G(F_0)}(\calT(\calO_0))$ of $\calT(\calO_0)$
in $\calG(F_0)$. Since
$N_{G(F_0)}(\calT(\calO_0))/\calT(\calO_0)\cong\widetilde{W}$ is
discrete, the induced $\bbG_m$-action fixes every element. The lemma
follows.
\end{proof}

\section{Proofs I: Frobenius splitting of global Schubert
varieties}\label{PfI} In this section, we prove Theorem \ref{sch
fiber} assuming Theorem \ref{top fiber}. We also deduce Theorem
\ref{MainI} from Theorem \ref{MainII}.

\subsection{Factorization of affine Demazure modules}\label{factor}
In this subsection, we show that Theorem \ref{MainII} implies
Theorem \ref{MainI}. This proof is essentially contained in \cite{Z1}.
Here, we repeat the arguments since they serve as a
prototype for the following subsections.

Let $H$ be a split Chevalley group over $k$ such that $H_\der$ is
almost simple, simply-connected, as assumed in \S \ref{grp data}. Let $\Gr_H$ be the affine
Grassmannian of $H$ and $\calL_b$ be the line bundle on
$(\Gr_H)_{\rm red}$ (the reduced ind-subscheme of $\Gr_H$), which restricts
to the ample generator of the Picard group of each of connected
component (which is isomorphic to $\Gr_{H_\der}$). We have the
following two assertions.

\begin{lem}Let $\mu\in\xcoch(T_H)$ be a minuscule coweight, so that $\bGr_\mu\cong X(\mu)=H/P(\mu)$, where $P(\mu)$ is the maximal parabolic subgroup corresponding to $\mu$. Then the restriction of $\calL_b$ to $\bGr_\mu$ is the ample generator of the Picard group of $X(\mu)$.
\end{lem}
\begin{proof}
Let us use the following notation. For $\nu$ a dominant weight of
$P(\mu)$, let $\calL(\nu)$ be the line bundle on $H/P(\mu)$, such
that $\Ga(H/P(\mu),\calL(\nu))^*$ is the Weyl module of $H$ of
highest weight $\nu$.

First assume that $\on{char} k=0$. Let us fix a normalized invariant
form $(\cdot,\cdot)_{norm}$ on $\xcoch(T_H)$ so that the square of
the length of short coroots is two. Note that this invariant form
may not be unique if $H$ is not semisimple.  For a coweight
$\mu\in\xcoch(T_H)$ of $T_H$, let $\mu^*$ be the image of $\mu$
under $\xcoch(T_H)\to\xch(T_H)$ induced by this form. In other
words, $(\mu^*,\la)=(\mu,\la)_{norm}$. Let $t_\mu\in
\Gr_{T_H}(k)\subset \Gr_H(k)$ be the corresponding point as in the
proof of Lemma \ref{moduli interpretation of sla}. Now assume that
$\mu$ is dominant. For every positive root $\tilde{a}$ of $H$, corresponding an $\SL_2\subset H$, let
$\SL_2 t_\mu\subset H t_\mu\cong H/P(\mu)$ be the corresponding rational curve. Then according to \cite[Lemma 2.2.2]{Z1}, the
restriction of $\calL_b$ to this rational curve has degree  $\frac{2(\mu,\tilde a)}{(\tilde a,\tilde a)}$. Therefore, the restriction of $\calL_b$ to $H t_\mu\cong H/P(\mu)$ is isomorphic to
$\calL(\mu^*)$. Note that for $\mu$ minuscule, $(\mu,\mu)_{norm}=2$, and therefore $\mu^*$ is the corresponding minuscule weight. The lemma follows in this case.

To prove the lemma in the case $\on{char}k>0$, observe that
everything is defined over $\bbZ$ (see \cite{F} where it is proven
that the Schubert varieties are defined over $\bbZ$ and commute with
base change). It is well-known that
$\Pic(H/P(\mu)_\bbZ)\cong\Pic(H/P(\mu)_k)\cong\bbZ$. The lemma
follows.
\end{proof}

The following proposition is essentially equivalent \cite[Theorem 1]{FL}, whose proof is of combinatoric nature. Here we reproduce a proof given in \cite[Theorem 1.2.2]{Z1}.
\begin{prop}Let $\calL$ be a line bundle on $(\Gr_H)_{\rm red}$, whose restriction to each connected component of $\Gr_H$ has the same central charge. Then
\[H^0(\bGr_{\mu+\la},\calL)\cong H^0(\bGr_\mu,\calL)\otimes H^0(\bGr_\la,\calL).\]
\end{prop}

\begin{proof}Recall that $H^1(\bGr_\mu,\calL)=0$ since $\bGr_\mu$ is Frobenius split and $\calL$ is ample. Therefore, it is enough to prove
the proposition for $\calL^n,\calL^{2n},\ldots$ and some $n\geq 1$.
Therefore we can replace $\calL$ by $\calL^n$, we can assume that
the central charge of $\calL$ is $2h^\vee$, i.e.
$\calL=\calL_b^{2h^\vee}$. Then $\calL^n$ is the pullback of the
$n$-tensor of the determinant line bundle $\calL_{\det}^n$ of
$\Gr_{\GL(\Lie H)}$ along $i:\Gr_H\to\Gr_{\GL(\Lie H)}$, as has been
discussed in the proof of Theorem \ref{line and central charge}. Let
us choose a complete curve (e.g. $\bar{C}=C\cup\{\infty\}$) and let
$\Bun_H$ be the moduli stack of $H$-bundles on the curve. Then we
know that $\calL$ is the pullback along $\Gr_H\to \Bun_H$ of a line bundle on $\Bun_H$ (which in turn is the pullback along $\Bun_H\to\Bun_{\GL(\Lie H)}$ of the determinant line bundle on $\Bun_{\GL(\Lie H)}$). Denote this
line bundle on $\Bun_H$ as $\omega^{-1}$ (in
fact this is the anti-canonical bundle of $\Bun_H$). 

Consider the convolution affine Grassmannian $\Gr_{H\times
C}^{Conv}$ over $C$, defined as
\[\Gr_{H\times C}^{Conv}(R)=\left\{(y,\calE,\calE',\beta,\beta')\ \left|\ \begin{split}&y\in C(R),\ \calE,\calE' \mbox{ are two } H \mbox{-torsors on } C_R,\\&\beta:\calE|_{C_R-\Ga_y}\cong \calE^0|_{C_R-\Ga_y}\mbox{ is a trivialization},\\&\mbox{and }
\beta':\calE'|_{(C-\{0\})_R}\cong
\calE|_{(C-\{0\})_R}.\end{split}\right.\right\}.\] This is an
ind-scheme formally smooth over $C$, and by the same argument as in \cite[Proposition 5]{G} we have
\[\Gr_{H\times C}^{Conv}|_{C^{\circ}}\cong \Gr_{H\times C^{\circ}}\times \Gr_H, \quad (\Gr_{H\times C}^{Conv})_0\cong \Gr_H\tilde{\times}\Gr_H,\]
where $\Gr_H\tilde{\times}\Gr_H:=LH\times^{L^+H}\Gr_H$ is the local
convolution Grassmannian. In addition, $\Gr_{H\times C}^{Conv}$ is a
fibration over $\Gr_{H\times C}$ by sending
$(y,\calE,\calE',\beta,\beta')$ to $(y,\calE,\beta)$, with the
fibers isomorphic to $\Gr_H$.

Now $\bGr_\mu\times\bGr_\la$ extends naturally to a closed variety
of $\Gr_{H\times C^{\circ}}\times \Gr_H$. The closure of this
variety in $\Gr_{H\times C}^{Conv}$ is denoted as $\bGr_{H\times
C,\mu,\la}^{Conv}$. As is proven in \cite[1.2.2]{Z1}, for $y\neq 0$,
$(\bGr_{H\times C,\mu,\la}^{Conv})_y\cong \bGr_\mu\times\bGr_\la$
and $(\bGr_{H\times C,\mu,\la}^{Conv})_0\cong
\bGr_\mu\tilde{\times}\bGr_\la$, where
$\bGr_\mu\tilde{\times}\bGr_\la$ is the ``twisted product" of
$\bGr_\mu$ and $\bGr_\la$ (see \emph{loc. cit.} or \eqref{twp} below
for the precise definition).

Let $h:\Gr_{H\times C}^{Conv}\to \Bun_H$ be the map sending
$(y,\calE,\calE',\beta,\beta')$ to $\calE'$. Then as explained in
\cite[1.2.2]{Z1}, $h^*(\omega^{-1})^n$, when restricted to $\Gr_{H\times
C}^{Conv}|_{C^{\circ}}$, is isomorphic to
$\calL^n\boxtimes\calL^n$, whereas over $(\Gr_{H\times
C}^{Conv})_0$, it is isomorphic to $m^*\calL^n$, where
$m:\Gr_H\tilde{\times}\Gr_H\to\Gr_H$ is the natural convolution
map (which is obtained from multiplication in the loop group). Therefore, \[H^0(\bGr_\mu,\calL^n)\otimes
H^0(\bGr_\la,\calL^n)\cong
H^0(\bGr_\mu\tilde{\times}\bGr_\la,m^*\calL^n)\cong
H^0(\bGr_{\mu+\la},\calL^n).\] The last isomorphism is due to the
fact $\calO_{\bGr_{\mu+\la}}\cong
m_*\calO_{\bGr_\mu\tilde{\times}\bGr_\la}$.
\end{proof}
Clearly, these two assertions together with Theorem \ref{MainII} will imply Theorem \ref{MainI}.

\subsection{Reduction of Theorem \ref{sch fiber} to Theorem \ref{Fsplit}}
In this subsection, we prove Theorem \ref{sch fiber}, assuming
Theorem \ref{top fiber}. The key ingredient is the Frobenius
splitting of varieties in characteristic $p$.

We begin with introducing more ind-schemes. Let $\calG$ be the group scheme over $C$ as in \S \ref{grpsch}. In particular, $\calG_{\calO_0}=\calG_{\sigma_Y}$. Let $\Gr_{\calG}^{BD}$ be
the Beilinson-Drinfeld affine Grassmannian for $\calG$ over $C$.
That is, for every $k$-algebra $R$,
\begin{equation}\label{BD Grass}
\Gr^{BD}_\calG(R)=\left\{(y,\calE,\beta)\ \left|\ \begin{split}&y\in
C(R), \calE \mbox{ is a } \calG\mbox{-torsor on }C_R,\mbox{ and}\\&
\beta:\calE|_{C^{\circ}_R-\Ga_y}\cong \calE^0|_{C^{\circ}_R-\Ga_y} \mbox{
is a trivialization}
\end{split}\right.\right\}.
\end{equation}
This is formally smooth ind-scheme ind-proper over $C$ (the ind-representability of $\Gr^{BD}_\calG$ is explained in the proof of Theorem 10.5 of \cite{PZ}). Again, by the same argument as in \cite[Proposition 5]{G}, we have
\[\Gr_\calG^{BD}|_{C^{\circ}}\cong\Gr_{\calG}|_{C^{\circ}}\times (\Gr_\calG)_0,\quad (\Gr_{\calG}^{BD})_0\cong(\Gr_\calG)_0\cong\Fl^Y.\]

We also need the convolution affine
Grassmannian $\Gr^{Conv}_\calG$. The functor it represents is as
follows. Let $R$ be a $k$-algebra.
\begin{equation}\label{Conv Grass}
\Gr^{Conv}_\calG(R)=\left\{(y,\calE,\calE',\beta,\beta')\ \left|\
\begin{split}& y\in C(R),\ \calE,\calE' \mbox{ are two }
\calG\mbox{-torsors on } C_R,
\\&  \beta:\calE|_{C_R-\Ga_y}\cong \calE^0|_{C_R-\Ga_y}
\mbox{ is a }\mbox{trivialization}, \\& \mbox{and }
\beta':\calE'|_{C^{\circ}_R}\cong\calE|_{C^{\circ}_R}\end{split}\right.\right\}.
\end{equation}
The ind-representibility of $\Gr^{Conv}_\calG$ can be seen from another construction of
$\Gr_\calG^{Conv}$. Namely, there is a $L^+\calG_{\calO_0}$-torsor
$\Gr_{\calG,\underline{0}}$ over $\Gr_\calG$ whose $R$-points
classify
\begin{equation}\label{B-torsor}
\Gr_{\calG,\underline{0}}(R)=\left\{(y,\calE,\beta,\ga)\ \left|\
\begin{split}& (y,\calE,\beta)\in\Gr_\calG(R), \mbox{ and a trivialization}\\ &
\ga:\calE|_{\widehat{\{0\}\times \Spec
R}}\cong\calE^0|_{\widehat{\{0\}\times \Spec
R}}\end{split}\right.\right\}
\end{equation}
where $\widehat{\{0\}\times \Spec R}$ is spectrum of its coordinate ring of the completion of $C_R$
along $\{0\}\times\Spec R$. Then
\[\Gr_\calG^{Conv}\cong \Gr_{\calG,\underline{0}}\times^{L^+\calG_{\calO_0}}\Fl^Y.\] 
The projection 
\[\pi:\Gr^{Conv}_\calG\to\Gr_\calG\] sends $(y,\calE,\calE',\beta,\beta')$ to $(y,\calE,\beta)$.

There is a natural map
\begin{equation}\label{p}
m:\Gr^{Conv}_\calG\to\Gr^{BD}_\calG
\end{equation}
sending $(y,\calE,\calE',\beta,\beta')$ to
$(y,\calE',\beta\circ\beta')$. This is a morphism over $C$, which is
an isomorphism over $C-\{0\}$. Over $0$, this morphism is the local
convolution morphism
\begin{equation}\label{loc conv}
m:\Fl^Y\tilde{\times}\Fl^Y:=LG\times^{L^+\calG_{\calO_0}}\Fl^Y\to\Fl^Y,
\end{equation}
given by the natural multiplication of the loop group.

In addition, there is a section
\[z:\Gr_\calG\to\Gr^{Conv}_\calG\]
given by sending $(y,\calE,\be)$ to $(y,\calE,\calE,\be,\id)$.
Therefore, via $z$ (resp. $m\circ z$), $\Gr_\calG$ is realized as closed
subschemes of $\Gr^{Conv}_\calG$ (resp. $\Gr^{BD}_\calG$). 

\medskip

Let $w\in W^Y\setminus\widetilde{W}/W^Y$ be an element in the
extended affine Weyl group and let $\Fl_w^Y$ denotes the
corresponding Schubert variety in $\Fl^Y$. Then
$\bGr_\mu\times\Fl_w^Y\subset\Gr_H\times\Fl^Y$ extends to a variety
\[\bGr_{\calG,\mu}|_{\tilde{C}^{\circ}}\times\Fl^Y_w\subset
(\Gr_\calG\times_C\tilde{C}^{\circ})\times\Fl^Y\cong\Gr^{BD}_{\calG}\times_C\tilde{C}^{\circ}\cong\Gr^{Conv}_\calG\times_C\tilde{C}^{\circ}.\]
Let $\bGr^{BD}_{\calG,\mu,w}$ denote its flat closure in
$\Gr^{BD}_\calG\times_C\tilde{C}$, and $\bGr^{Conv}_{\calG,\mu,w}$
denote its flat closure in
$\Gr^{Conv}_\calG\times_C\tilde{C}$. Then
$\bGr^{Conv}_{\calG,\mu,w}$ maps to $\bGr^{BD}_{\calG,\mu,w}$ via
$m$. In addition, we have
\begin{lem}\label{fibration}
$\bGr^{Conv}_{\calG,\mu,w}$ is a fibration (via $\pi$) over
$\bGr_{\calG,\mu}$ with fibers isomorphic to $\Fl_w^Y$.
\end{lem}
\begin{proof}
From the second definition of $\Gr_{\calG}^{Conv}$, it is clear that
\[\bGr^{Conv}_{\calG,\mu,w}\cong \bGr_{\calG,\mu}\times_{\Gr_\calG}(\Gr_{\calG.\underline{0}}\times^{L^+\calG_{\calO_0}}\Fl^Y_w).\]
\end{proof}
As by definition $(\bGr_{\calG,\mu})_0$ is a $L^+\calG_{\calO_0}$-stable closed subscheme of $\Fl^Y$ and $(\bGr^{BD}_{\calG,\mu,\nu})_{\tilde{0}}$ is the image of $(\bGr^{Conv}_{\calG,\mu,w})_0$ under $m$, we obtain 
\begin{cor}\label{BDstable}
The scheme $(\bGr^{BD}_{\calG,\mu,\nu})_{\tilde{0}}$ is a $L^+\calG_{\calO_0}$-stable closed subscheme of $\Fl^Y$.
\end{cor}

\begin{prop}\label{key prop}
Let $\nu\in\xcoch(T)$ be a sufficiently dominant coweight. Then the
variety $\bGr^{BD}_{\calG,\mu,\nu}$ is normal and the fiber
$(\bGr^{BD}_{\calG,\mu,\nu})_{\tilde{0}}$ over $\tilde{0}\in
\tilde{C}$ is reduced.
\end{prop}
\begin{proof} The key observation
\begin{lem}Suppose that $\nu$ is sufficiently large. Then the fiber
$(\bGr^{BD}_{\calG,\mu,\nu})_{\tilde{0}}$ is irreducible and
generically reduced.
\end{lem}
Assuming the lemma, then the proposition follows from Hironaka's
lemma (cf. EGA IV.5.12.8). Namely, let $V$ denote the underlying
reduced subscheme of $(\bGr^{BD}_{\calG,\mu,\nu})_{\tilde{0}}$. Then
$V$ is irreducible, and therefore by Corollary \ref{BDstable}, is a Schubert variety of $\Fl^Y$, which is normal by Theorem
\ref{NSV}. Therefore, the proposition follows.

So it remains to prove the lemma. Let us first prove that
$(\bGr^{BD}_{\calG,\mu,\nu})_{\tilde{0}}$ is irreducible. Clearly,
$\bGr^{Conv}_{\calG,\mu,\nu}$ maps surjectively onto
$\bGr^{BD}_{\calG,\mu,\nu}$. Therefore,
$(\bGr^{Conv}_{\calG,\mu,\nu})_{\tilde{0}}$ dominates
$(\bGr^{BD}_{\calG,\mu,\nu})_{\tilde{0}}$. We know that
$(\bGr^{Conv}_{\calG,\mu,\nu})_{\tilde{0}}$ is a fibration over
$(\bGr_{\calG,\mu})_{\tilde{0}}$, with fibers isomorphic to
$\Fl_\nu^Y$. Therefore by Theorem \ref{top fiber}, the underlying
reduced subschemes of irreducible components of
$(\bGr^{Conv}_{\calG,\mu,\nu})_{\tilde{0}}$ are just
\[
\Fl^Y_\la\tilde{\times}\Fl^Y_\nu,\quad \la\in\Lambda.
\]
Here and in what follows we use the following notation: Let $S_1,S_2$ are two
subschemes of $\Fl^Y$ and assume that $S_2$ is $L^+\calG_{\calO_0}$
stable, then we denote
\begin{equation}\label{twp}S_1\tilde{\times}S_2:=\widetilde{S_1}\times^{L^+\calG_{\calO_0}}S_2,\end{equation}
where $\widetilde{S_1}$ is the preimage of $S_1$ under
$L\calG_{\calO_0}\to\Fl^Y$.

Therefore, the underlying reduced subscheme of each irreducible
component of
\[(\bGr^{BD}_{\calG,\mu,\nu})_{\tilde{0}}\subset\Fl^Y\]
is contained in one of $m(\Fl^Y_\la\tilde{\times}\Fl^Y_\nu),
\la\in\Lambda$. Observe that if $\la\in\Lambda$ is not dominant, for $\nu$ sufficiently dominant so that $\la+\nu$ is dominant,
we have
\[\ell(t_{\nu+\la})=(2\rho,\nu+\la)<(2\rho,\nu)+(2\rho,\mu)= \ell(t_\nu)+\ell(t_\la)\]
by Lemma \ref{length}. However, by
flatness, all the irreducible components of
$(\bGr^{BD}_{\calG,\mu,\nu})_{\tilde{0}}$ have dimension
$\ell(t_\mu)+\ell(t_\nu)$. This implies that
$(\bGr^{BD}_{\calG,\mu,\nu})_{\tilde{0}}$ has only one irreducible
component, whose underlying reduced subscheme is
$m(\Fl^Y_\mu\tilde{\times}\Fl^Y_\nu)=\Fl^Y_{\mu+\nu}$.

Next, we prove that $(\bGr^{BD}_{\calG,\mu,\nu})_{\tilde{0}}$ is
generically reduced. In \S \ref{affine chart}, we have constructed the
affine open chart $c_\mu: U_\mu\subset\bGr_{\calG,\mu}$ that satisfies:
\begin{enumerate}
\item $s_\mu(\tilde{C})\subset U_\mu$;
\item $U_\mu$ is an affine space over $\tilde{C}$ and therefore smooth over
$\tilde{C}$;
\item $(U_\mu)_{\tilde{0}}=C(\mu)\subset\Fl^Y$ is the Schubert cell containing
$t_\mu$, i.e. the $L^+\calG_{\mathbf{a}}$-orbit containing $t_\mu$.
\end{enumerate}
Let us restrict $\bGr^{Conv}_{\calG,\mu,\nu}$ over $U_\mu$. Then
clearly, the fiber over $\tilde{0}$ of this family is
$(U_\mu)_{\tilde{0}}\tilde{\times}\Fl_\nu^Y$ and therefore is
irreducible and reduced. Let $\xi$ be the generic point of
$(U_\mu)_{\tilde{0}}\tilde{\times}\Fl_\nu^Y$. By the above argument,
$\eta=m(\xi)$ is the generic point of
$(\bGr^{BD}_{\calG,\mu,\nu})_{\tilde{0}}$. Let $A$ denote the local
ring of $\bGr^{BD}_{\calG,\mu,\nu}$ at $\eta$ and $B$ denote the
local ring of $\bGr^{Conv}_{\calG,\mu,\nu}$ at $\xi$. Both are
discrete valuation rings, flat over $\tilde{C}$, and there is an
injective map $A\to B$ which is an isomorphism over
$\tilde{C}^{\circ}$. Since $\bGr^{Conv}_{\calG,\mu,\nu}$ is
proper over $\tilde{C}$, we obtain a morphism $\Spec
A\to\bGr^{Conv}_{\calG,\mu,\nu}$ which must factors through $\Spec
A\to\Spec B$. That is, $A\to B$ is split injective. Therefore
$A/uA\subset B/uB$ is a subfield. That is
$\bGr^{BD}_{\calG,\mu,\nu}$ is generically reduced.
\end{proof}

In fact, we proved that the fiber
$(\bGr^{BD}_{\calG,\mu,\nu})_{\tilde{0}}$ is isomorphic to
$\Fl^Y_{\mu+\nu}$.

If $\nu\in \xcoch(T_\s)_\Ga\subset \xcoch(T)_\Ga$ so that $\nu\in
W_\aff$, then $z(\bGr_{\calG,\mu})\subset
\bGr^{Conv}_{\calG,\mu,\nu}$ (resp.
$m\circ z(\bGr_{\calG,\mu})\subset\bGr^{BD}_{\calG,\mu,\nu}$) is naturally
a closed subscheme.

By Corollary \ref{change char}, we just need to prove Theorem \ref{sch fiber} for one prime.
Therefore, we will assume $\on{char}k=p>2$. Recall the notation of Frobenius splitting (cf. \cite{MR,BK}) for varieties in characterisic $p>0$.

\begin{thm}\label{aux}
Assume that $\nu\in \xcoch(T_\s)$ is a coweight dominant enough so
that Proposition \ref{key prop} holds. Then $\bGr^{BD}_{\calG,\mu,\nu}$
is Frobenius split, compatibly with $\bGr_{\calG,\mu}$ and
$(\bGr^{BD}_{\calG,\mu,\nu})_{\tilde{0}}$.
\end{thm}

\begin{cor}Theorem \ref{sch fiber} holds. That is, the scheme $(\bGr_{\calG,\mu})_{\tilde{0}}$ is reduced.
\end{cor}
\begin{proof}This is because that
\[(\bGr_{\calG,\mu})_{\tilde{0}}=\bGr_{\calG,\mu}\cap(\bGr^{BD}_{\calG,\mu,\nu})_{\tilde{0}},\]
and therefore is Frobenius split. In particular, it is reduced.
\end{proof}

The remaining goal of this subsection is to reduce Theorem \ref{aux}
to the following Theorem \ref{Fsplit} via Proposition \ref{aux2}. Theorem
\ref{Fsplit} itself will be proven in later subsections. First, it is enough to prove Theorem \ref{aux} for the case
$\calG_{\calO_0}\cong\calG_{\mathbf{a}}$ is the Iwahori group scheme.
To see this, assume that we have
$\calG_1\to\calG_2$ where $(\calG_1)_{\calO_0}$ is Iwahori and
$(\calG_2)_{\calO_0}$ is a general parahoric group scheme. Then the
natural projection
$\bGr_{\calG_1,\mu,\nu}^{BD}\to\bGr_{\calG_2,\mu,\nu}^{BD}$ is proper
birational and therefore the push-forward of the structure sheaf is
the structure sheaf by the normality. Furthermore, under the
projection, the scheme-theoretical image of $\bGr_{\calG_1,\mu}$
(resp. $(\bGr^{BD}_{\calG_1,\mu,\nu})_{\tilde{0}}$) is
$\bGr_{\calG_2,\mu}$ (resp.
$(\bGr^{BD}_{\calG_2,\mu,\nu})_{\tilde{0}}$). Therefore, from now on we assume that $\calG_{\calO_0}=\calG_{\mathbf{a}}$ and write
$I=L^+\calG_{\mathbf{a}}$.

Since $\bGr^{BD}_{\calG,\mu,\nu}$ is normal, we just need to find an
open subscheme of $U\subset\bGr^{BD}_{\calG,\mu,\nu}$, whose
complement has codimension at least two, such that $U$ is Frobenius
split, compatibly with $U\cap
(\bGr^{BD}_{\calG,\mu,\nu})_{\tilde{0}}$ and $U\cap
\bGr_{\calG,\mu}$ (\cite[Lemma 1.1.7 (iii)]{BK}). Therefore, we can throw away some bad loci of
$\bGr^{BD}_{\calG,\mu,\nu}$ which is hard to control. In particular,
we can throw away
$(\bGr_{\calG,\mu})_{\tilde{0}}\subset\bGr_{\calG,\mu}\subset\bGr^{BD}_{\calG,\mu,\nu}$
which is of our main interests!

More precisely, we have
\begin{prop}\label{aux2}
There is an open subscheme $U$ of $\bGr^{Conv}_{\calG,\mu,\nu}$,
such that
\begin{enumerate}
\item $m:\bGr^{Conv}_{\calG,\mu,\nu}\to \bGr^{BD}_{\calG,\mu,\nu}$ maps $U$ isomorphically onto an open subscheme $m(U)$ of $\bGr^{BD}_{\calG,\mu,\nu}$, and the complement of $p(U)$ in $\bGr^{BD}_{\calG,\mu,\nu}$ has codimension two;
\item  $U$ is Frobenius split, compatible with $U\cap(\bGr^{Conv}_{\calG,\mu,\nu})_{\tilde{0}}$ and $U\cap z(\bGr_{\calG,\mu})$.
\end{enumerate}
\end{prop}
It is clear that Theorem \ref{aux} will follow from this proposition.
\begin{proof}
Let us first construct this open subscheme $U$. Recall that we
constructed the section $s_\mu:\tilde{C}\to \wGr_\calG$ and
$\bGr_{\calG,\mu}$ is the minimal irreducible closed subvariety of
$\wGr_\calG$ that is invariant under $\wJG$ and contains
$s_\mu(\tilde{C})$. Let $\Gr_{\calG,\mu}$ denote the $\wJG$-orbit
through $s_\mu$. Then $\Gr_{\calG,\mu}$ is an open subscheme of
$\bGr_{\calG,\mu}$, which is smooth over $\tilde{C}$. In fact $\Gr_{\calG,\mu}$ is open in $\bGr_{\calG,\mu}$ since each point in $\Gr_{\calG,\mu}$ can be translated to an element in $s_\mu(\tilde{C})$, which is contained in the open subset $U_\mu$ of $\bGr_{\calG,\mu}$.
latter is the closure of the former. Therefore $\Gr_{\calG,\mu}$ is
flat over $\tilde{C}$. Observe that under the isomorphism
$\Gr_\calG\times_C\tilde{C}^{\circ}\cong\Gr_H\times\tilde{C}^{\circ}$,
\[\Gr_{\calG,\mu}|_{\tilde{C}^{\circ}}\cong \Gr_\mu\times\tilde{C}^{\circ},\]
where $\Gr_\mu$ denotes the $L^+H$-orbit in $\Gr_H$ through $t_\mu$,
which is smooth. On the other hand
$(\Gr_{\calG,\mu})_{\tilde{0}}=(U_\mu)_{\tilde{0}}$ is the Schubert
cell $C(\mu)$ in $\Fl$ containing $\mu$, which is irreducible and
smooth. Therefore, $\Gr_{\calG,\mu}$ is smooth over $\tilde{C}$.

Let $U_1$ be the preimage of $\Gr_{\calG,\mu}$ under
$\pi:\bGr^{Conv}_{\calG,\mu,\nu}\to\bGr_{\calG,\mu}$. Then $U_1$ is
a fibration over $\Gr_{\calG,\mu}$ with fibers $\Fl_\nu$. As a
scheme over $\tilde{C}$, the fiber of $U_1$ over $\tilde{0}$ is
\[C(\mu)\tilde{\times}\Fl_\nu.\]
We define $U$ to be the open subscheme of $U_1$ which coincides with
$U_1$ over $\tilde{C}^{\circ}$, and which is given by
\[C(\mu)\tilde{\times}C(\nu)\subset C(\mu)\tilde{\times}\Fl_\nu\]
over $\tilde{0}$.

We claim that $m:U\to m(U)$ is an isomorphism and the complement of
$m(U)$ in $\bGr^{BD}_{\calG,\mu,\nu}$ has codimension two. Over
$\tilde{C}^{\circ}$, $m$ is an isomorphism. Over $\tilde{0}$, the
morphism
\[m:U_{\tilde{0}}\to (\bGr^{BD}_{\calG,\mu,\nu})_{\tilde{0}}\]
is the same as
\[m: C(\mu)\tilde{\times}C(\nu)\to \Fl_{\mu+\nu}.\]
It is well-known (e.g \cite[IV]{Ma}) that over $\tilde{0}$ $m$ induces an isomorphism from
$C(\mu)\tilde{\times}C(\nu)$ onto $C(\mu+\nu)$, and the preimage of $C(\mu+\nu)$ is $C(\mu)\tilde{\times}C(\mu)$. Therefore, $m:U\to m(U)$ is a homeomorphism and $m^{-1}m(U)=U$. Therefore, $m:U\to m(U)$ is
a proper, birational homeomorphism with $m(U)$ normal, which must be
an isomorphism. Note that
$(\bGr_{\calG,\mu})_{\tilde{0}}\subset\bGr_{\calG,\mu}\subset\bGr_{\calG,\mu,\nu}^{BD}$
is not contained in $m(U)$.

To see that the complement of $m(U)$ has codimension two, first
observe that over $\tilde{C}^{\circ}$,
\[\bGr^{BD}_{\calG,\mu,\nu}|_{\tilde{C}^{\circ}}-m(U)|_{\tilde{C}^{\circ}}\cong (\bGr_\mu-\Gr_\mu)\times \Fl_\nu\times\tilde{C}^{\circ},\]
which has codimension two, since $\bGr_\mu-\Gr_\mu$ has codimension
two in $\bGr_\mu$. Over $\tilde{0}$,
\[(\bGr^{BD}_{\calG,\mu,\nu})_{\tilde{0}}-m(U)_{\tilde{0}}\cong \Fl_{\mu+\nu}-C(\mu+\nu),\]
which has codimension at least one. This proves that the complement
of $m(U)$ in $\bGr^{BD}_{\calG,\mu,\nu}$ has codimension two.

\medskip

Next we turn to the second part of the proposition. Recall that
$U_1$ is the preimage of $\Gr_{\calG,\mu}$ under
$\pi:\bGr^{Conv}_{\calG,\mu,\nu}\to\bGr_{\calG,\mu}$. From the
construction of $U$, we know that $U\subset
U_1\subset\bGr^{Conv}_{\calG,\mu,\nu}$. Therefore, it is enough to
show that the same statement of Proposition \ref{aux2} (2) holds for
$U_1$. Recall that
\[U_1\cong (\Gr_{\calG,\mu}\times_{\Gr_\calG}\Gr_{\calG,\underline{0}})\times^{I}\Fl_\nu,\]
where $\Gr_{\calG,\underline{0}}$ is the $I$-torsor over $\Gr_\calG$
as in \eqref{B-torsor}. To simplify
the notation, for any $I$-variety $V$, we denote 
\[\Gr_{\calG,\mu}\tilde{\times}V:=(\Gr_{\calG,\mu}\times_{\Gr_\calG}\Gr_{\calG,\underline{0}})\times^{I}V.\]
Now, let $\ast\in\Fl_\nu$ be the base point
(recall that $\nu\in \xch(T_\s)$, so that $\ast$, the Schubert
variety corresponding to the identity element in the affine Weyl
group, is contained in $\Fl_\nu$). Then the closed embedding
$z:\Gr_{\calG,\mu}\to U_1$ corresponds to
\[\Gr_{\calG,\mu}\tilde{\times}\ast\to\Gr_{\calG,\mu}\tilde{\times}\Fl_\nu.\]
Now the assertion follows from the following more general statement.
\end{proof}

\begin{thm}\label{Fsplit}
For any $w\in \widetilde{W}$, there is a Frobenius splitting of
$\Gr_{\calG,\mu}\tilde{\times}\Fl_w$, compatibly with
\[(\Gr_{\calG,\mu}\tilde{\times}\Fl_w)_{\tilde{0}}\cong(\Gr_{\calG,\mu})_{\tilde{0}}\tilde{\times}\Fl_w\cong
C(\mu)\tilde{\times}\Fl_w.\] In addition, for any $v\leq w$ in
$\widetilde{W}$,
$\Gr_{\calG,\mu}\tilde{\times}\Fl_v\subset\Gr_{\calG,\mu}\tilde{\times}\Fl_w$
is also compatible with this splitting.
\end{thm}

The remaining goal of this section is to prove this theorem.

\subsection{Special parahorics}\label{sp}
We continue assume that $G$ and $\calG$ are as given in \S \ref{grpsch} but we are particularly interested in the case when $\calG=\calG^s$ is the group scheme over $C$ such that
$\calG^s_{\calO_0}$ is a special parahoric group scheme of $G$. In
this case, we can easily deduce Theorem \ref{sch fiber} (assuming
Theorem \ref{top fiber}) directly from Hironaka's lemma (without
going into the argument presented in the previous subsection). This
will in turn help us prove a special case of Theorem \ref{Fsplit},
namely, the case when $w=1$ (see Corollary \ref{aux3}). Let us
remark that if $G$ is split, Proposition \ref{special parahoric
case} directly follows from Frobenius splitting of Schubert varieties, and those who are only interested in split groups
can go to the paragraph after this proposition directly.

So let $v\in A(G,S)$ be a special point in the apartment associated
to $(G,S)$, and let $\calG_v$ be the corresponding special parahoric
group scheme over $\calO$. Let $\Fl_v=LG/L^+\calG_v$ be the partial
affine flag variety. To emphasize that it is the affine flag variety
associated to a special parahoric, we sometimes also denote it by
$\Fl^s$. As before, for each $\mu\in\xcoch(T)_\Ga$, let us use
$t_\mu$ to denote its lifting to $T(F)$ under the Kottwitz
homomorphism $T(F)\to\xcoch(T)_\Ga$. It gives a point in $\Fl^s$,
still denoted by $t_\mu$. Then the Schubert variety $\Fl^s_\mu$ is
the closure of the $L^+\calG_v$-orbit in $\Fl_v$ passing through
$t_\mu$. We have the following results special for Schubert
varieties in $\Fl^s$, which generalize the corresponding results
for $\Gr_H$ (see also \cite[Corollary 2.10]{Ri} for more detailed discussion).

\begin{lem}The Schubert varieties are parameterized by
$\xcoch(T)_\Ga^+$. For $\mu\in\xcoch(T)_\Ga^+$, the dimension of
$\Fl^s_\mu$ is $(\mu,2\rho)$. Let
$\mathring{\Fl}^s_\mu\subset\Fl^s_\mu$ be the unique open
$L^+\calG_v$-orbit in $\Fl^s_\mu$. Then
$\Fl^s_\mu-\mathring{\Fl}^s_\mu$ has codimension at least two.
\end{lem}
\begin{proof}Observe that the natural map
$\xcoch(T)_\Ga^+\subset\xcoch(T)_\Ga\to{W_0}\setminus\widetilde{W}/{W_0}$
is a bijection. The first claim follows. Let $I\subset L^+\calG_v$
be the Iwahori subgroup of $LG$ corresponding to the alcove $\mathbf{a}$
(recall that $v$ is contained in the closure of $C$). Then the
$I$-orbits in $\Fl^s$ are parameterized by minimal length
representatives in $\widetilde{W}/{W_0}$. Let
$\la\in\La={W_0}\mu\subset\xcoch(T)_\Ga$. By Lemma \ref{length} and \ref{length2}, if
$w\in\widetilde{W}$ is a minimal length representative for the coset
$t_\la {W_0}$, then
\[\dim IwL^+\calG_v/L^+\calG_v\leq (\mu,2\rho)\]
and if $\la\in\xcoch(T)_\Ga^+$, the equality holds. Therefore,
$\dim\Fl^s_\mu=(\mu,2\rho)$. To prove the last claim, observe that
if $\Fl^s_\la\subsetneq\Fl^s_\mu$, then $\mu-\la\in\xcoch(T_\s)_\Ga$
and therefore $(\mu-\la,2\rho)$ is an even integer.
\end{proof}

Recall that in \cite[\S 4.6]{BD}, Beilinson and Drinfeld proved that
$\bGr_\mu$ is Gorenstein, i.e., the dualizing sheaf
$\omega_{\bGr_\mu}$ is indeed a line bundle (see Equation (241) in \emph{loc. cit.}). It is natural to ask
whether the same result hold for $\Fl^s_\mu$. However, the
situation is more complicated in the ramified case due to the fact
that not all special points in the building of $G$ are conjugate
under $G_\ad(F)$. More precisely, if $G_\der$ is the
odd ramified special unitary group $\SU_{2n+1}$ (see \S \ref{lines on local models} for
the definition), then there are two types of special parahoric group
schemes (see Remark \ref{special parahoric for unitary} (ii)).

Let us begin with the following lemma. Let $v$ be any point in the
apartment $A(G,S)$ and let $\calG_v$ be the corresponding parahoric
group scheme for $G$. For simplicity, we write $K=L^+\calG_v$. Then
$K$ acts on $\Lie G$ by the adjoint representation. Let
$\mu\in\xcoch(T)_\Ga$. Let \[P=K\cap\Ad_{t_\mu}K,\]
considered as a proalgebraic group over $k$.
Then $\Lie K$
and $\Ad_{t_\mu}\Lie K$ are $P$-modules.

\begin{lem}\label{P-mod}
As $P$-modules,
\begin{equation}\label{isom}\det\frac{\Lie
K}{\Lie K\cap\Ad_{t_\mu}\Lie K}\cong (\det\frac{\Ad_{t_\mu}\Lie
K}{\Lie K\cap\Ad_{t_\mu}\Lie K})^{-1}.
\end{equation}
\end{lem}
\begin{proof}Recall that we denote by $S$ the chosen maximal split $F$-torus
of $G$. Its (connected) N\'{e}ron model $\calS$ maps naturally into
$\calG_{v}$ since $v\in A(G,S)$ (\cite[Sect. 5.2]{BT2}), and $L^+\calS$ maps to $P$. The special fiber $S_k$ of $\calS$ can be regarded as the ``constant" maps from $\calO$ to $\calS$, and therefore can be regarded a subgroup of $L^+\calS$. Then $S_k\subset P$ is a maximal torus of $P$. Therefore,
$\xch(P)\subset \xch(S_k)$. Therefore, it is enough to prove
\eqref{isom} as $S_k$-modules.

In \S \ref{root subgroup}, in particular Remark \ref{local root
group} (see also \cite[9.a,9.b]{PR1}), we have attached to each affine root
$\al$ of $(G,S)$ a 1-dimensional unipotent subgroup
$U_\al\cong\bbG_a\subset LG$. Let
$\fraku_\al$ be the Lie algebra of $U_\alpha$. By definition
\[\Lie K=\Lie \calT^{\flat,0}\oplus\prod_{\al(v)\geq 0}\fraku_{\al},\]
where $\calT^{\flat,0}$ is the connected N\'{e}ron model of $T$.
Then clearly, as $S_k$-modules (we fix an embedding $S_k\to
L^+\calS$)
\[\frac{\Lie K}{\Lie K\cap\Ad_{t_\mu}\Lie K}\cong\bigoplus_{\al(v)\geq 0,\al(v-\mu)<0} \fraku_{\al},\quad \frac{\Ad_{t_\mu}\Lie K}{\Lie K\cap\Ad_{t_\mu}\Lie K}\cong\bigoplus_{\al(v)<0,\al(v-\mu)\geq 0} \fraku_{\al}.\]
By identifying $A(G,S)$ with $\xcoch(S)_\bbR$ using the point $v$,
we can write affine roots of $G$ by $\al=a+m, a\in\Phi(G,S)$, where
$a$ is the vector part of $\al$ and $m=v(\al)$. Therefore,
\[\{\al(v)\geq 0,\al(v-\mu)<0\}=\{a+m|a\in\Phi(G,S)^+, 0\leq
m<(\mu,a)\}\] and \[\{\al(v)< 0,\al(v-\la)\geq
0\}=\{a+m|a\in\Phi(G,S)^-, (\mu,a)\leq m<0\}.\] Since $S_k$ acts on
$\fraku_{a+m}$ via the weight $a$, the lemma follows.
\end{proof}

\quash{
\begin{rmk}If $G$ is split and $v$ is a hyperspecial vertex, under some mild restriction of the characteristic of $k$, one can
even show that as $P$-modules,
\begin{equation}\label{P-module}\frac{\Lie
K}{\Lie K\cap\Ad_{t_\mu}\Lie K}\cong(\frac{\Ad_{t_\mu}\Lie K}{\Lie
K\cap\Ad_{t_\mu}\Lie K})^*.\end{equation} This is because that there
is an invariant non-degenerate symmetric bilinear form on $\Lie LG$
which induces a $P$-invariant non-degerate bilinear form
\[\frac{\Lie K}{\Lie K\cap\Ad_{t_\mu}\Lie K}\otimes
\frac{\Ad_{t_\mu}\Lie K}{\Lie K\cap\Ad_{t_\mu}\Lie K}\to k.\] I do
not know whether \eqref{P-module} is true if $G$ is ramified and if
$v$ is a special vertex, but the argument fails.
\end{rmk}}

Now we should specify the special vertex. Recall that we assume that $G_\der$ is simple and simply-connected. If $G_\der\neq\SU_{2n+1}$,
we can choose arbitrary special vertex in the building of $G$ since
they are conjugate under $G_\ad(F)$. If $G_\der=\SU_{2n+1}$, we
choose the special vertex so that the corresponding parahoric group
has reductive quotient $\on{Sp}_{2n}$ (see Remark \S \ref{special
parahoric for unitary}).

\begin{thm}\label{Goren}
Let $G$ as in \S \ref{grpsch}. With the choice of the special vertex $v$ as above, the Schubert
variety $\Fl^s_\mu$ is Gorenstein for all $\mu$.
\end{thm}
\begin{proof}As above, we denote by $\calG_v$ the
parahoric group of $G$ corresponding to $v$ and $K=L^+\calG_v$. Recall that $\Fl^s_\mu$ is Cohen-Macaulay, the dualizing sheaf $\omega_{\Fl^s_\mu}$ exists. We
need to show that it is indeed
a line bundle. Let $j:\mathring{\Fl}^s_\mu\to\Fl^s_\mu$ be the open $K$-orbit in
$\Fl^s_\mu$. Then we have shown that
$\Fl^s_\mu-\mathring{\Fl}^s_\mu$ has codimension at least two. As $\Fl^s_\mu$ is normal, $\omega_{\Fl^s_\mu}=j_*(\omega_{\mathring{\Fl}^s_\mu})$.  Let
$\calL_{2c}$ be the pullback to $\Fl^s$ of the determinant line
bundle $\calL_{\det}$ of $\Gr_{\GL(\Lie\calG_{v})}$ along
$i:\Fl^s\to\Gr_{\GL(\Lie\calG_{v})}$. We first prove that there is an isomorphism of line bundles
$\omega_{\mathring{\Fl}^s_\mu}^{-2}\cong
\calL_{2c}|_{\mathring{\Fl}^s_\mu}$ on $\mathring{\Fl}^s_\mu$.

Indeed, observe that both sheaves are $K$-equivariant. The
$K$-equivariant structure of $\omega_{\mathring{\Fl}^s_\mu}^{-2}$ is
induced from the action of $K$ on $\mathring{\Fl}^s_\mu$. On the
other hand, a central extension of $LG$ acts on $\calL_{2c}$, and a
splitting of this central extension over $K$ defines a
$K$-equivariant structure on $\calL_{2c}$. To fix this
$K$-equivariant structure uniquely, we will require that the action of
$K$ on the fiber of $\calL_{2c}$ over $\ast\in\Fl^s$ is trivial.
Then the $K$-equivariant structure on $\calL_{2c}$ is given as
follows (for simplicity, we only describe it at the level of
$k$-points, but the generalization to $R$-points is clear, for
example see cf. \cite[\S 2.2.2-2.2.3]{FZ}): recall that for $x\in\Fl^s$, $i(x)$ is a
lattice in $\Lie G$ and $\calL_{2c}|_x$ is the $k$-line
\[\calL_{2c}|_x=\det(i(x)|\Lie K):=\det \frac{\Lie K}{\Lie K\cap i(x)}\otimes \det(\frac{i(x)}{\Lie K\cap i(x)})^{-1}\]
Then for $g\in K$, $\calL_{2c}|_x\to\calL_{2c}|_{gx}$ is given by
\[\det(g):\det(i(x)|\Lie K)\cong\det(i(gx)|g\Lie K)=\det(i(gx)|\Lie
K).\]

Now it is enough to prove that there is an isomorphism
$\calL_{2c}|_{t_\mu}\cong
\omega_{\mathring{\Fl}^s_\mu}^{-2}|_{t_\mu}$ as $1$-dimensional
representations of $P=t_\mu K t_\mu^{-1}\cap K$, the stabilizer of
$t_\mu\in\Fl^s$ in $K$. As the tangent space of $\mathring{\Fl}^s_\mu$ at $t_\mu$ as a $P$-module is isomorphic to $\frac{\Lie K}{\Lie
K\cap\Ad_{t_\mu}\Lie K}$,
\[\omega_{\mathring{\Fl}^s_\mu}^{-2}|_{t_\mu}\cong (\det\frac{\Lie K}{\Lie
K\cap\Ad_{t_\mu}\Lie K})^2\] as $P$-modules. On the other hand, it
follows from the construction of the determinant line bundle that
\[\calL_{2c}|_{t_\mu}\cong \det\frac{\Lie K}{\Lie
K\cap\Ad_{t_\mu}\Lie K}\otimes (\det\frac{\Ad_{t_\mu}\Lie K}{\Lie
K\cap\Ad_{t_\mu}\Lie K})^{-1}\] as $P$-modules. Therefore, the
assertion follows from Lemma \ref{P-mod}.

Next, we prove that there is a $K$-equivariant line bundle $\calL_c$ on
$(\Fl^s)_{\rm red}$ such that $\calL_c^2\cong\calL_{2c}$. Indeed, for
any $g\in G(F)$ acting on $\Fl^s$ by left translation, we have
$g^*\calL_{2c}\cong\calL_{2c}$. Therefore, it is enough to construct
$\calL_c$ in the neutral connected component of $(\Fl^s)_{\rm red}$,
which is isomorphic to $\Fl^s_\s$, the corresponding affine flag
variety for $G_\der$ by \cite[\S 6]{PR1}. Since $v$ is a special
vertex, $\Pic(\Fl^s_\s)\cong\bbZ\calL(\epsilon_i)$, where $i\in\bold
S$ is a special vertex in the local Dynkin diagram of $G$
corresponding to $v$. By checking \cite[\S 4, \S
6]{Kac}, we see that for our choice of $v$,  we have$a_i^\vee=1$. (For $\SU_{2n+1}$,
there is another special vertex $i'\in\bold S$ such that
$a_{i'}^\vee=2$, and the reductive quotient of the corresponding
parahoric group is $\on{SO}_{2n+1}$, see the following remark and
Remark \ref{special parahoric for unitary}). Therefore, the central
charge of $\calL(\epsilon_i)$ is $1$, whereas the central charge of
$\calL_{2c}$ is $2h^\vee$ by \eqref{cc} and Lemma \ref{adj}. Therefore,
$\calL_c=\calL(h^\vee\epsilon_i)$.

As $\xch(P)\subset\xch(S_k)$ is torsion free, $\calL_c|_{t_\mu}\cong \det\frac{\Lie K}{\Lie
K\cap\Ad_{t_\mu}\Lie K}$ as $P$-modules. Therefore, we have
$\omega_{\mathring{\Fl}^s_\mu}^{-1}\cong\calL_c|_{\mathring{\Fl}^s_\mu}$,
which in turn implies that $\omega_{\Fl^s_\mu}^{-1}=j_*(\omega_{\mathring{\Fl}^s_\mu})\cong j_*(\calL_c|_{\mathring{\Fl}^s_\mu})=\calL_c$.
\end{proof}

Observe the above proof implies that no matter what special vertex
we choose, $\omega_{\Fl^s_\mu}^{-2}$ is always a line bundle, where following \S \ref{Frob}, we denote $j_*((\omega_{\Fl^s_\mu}|_{\mathring{\Fl}^s_\mu})^n)$ by $\omega_{\Fl^s_\mu}^n$.

The
following corollary is what we need in the sequel.
\begin{cor}\label{vanishing H1}
For any special vertex $v$ of $G$,
$H^1(\Fl^s_\mu,\omega_{\Fl^s_\mu}^{-n})=0$ for all positive even integers $n$.
\end{cor}

\begin{rmk}\label{non-Goren}
In the case $G_\der=\SU_{2n+1}$, if we take the special vertex to be
$v_0$, the one defined by the pinning \eqref{pinned isom} so that
$\calG_{v_0}$ is of the form \eqref{constr of special}, then the
reductive quotient is $\on{SO}_{2n+1}$ and the corresponding
$a_i^\vee=2$. Since the dual Coxeter number of $\SL_{2n+1}$ is
$2n+1$, this means that on the partial flag variety $\Fl^s$
corresponding to this special vertex, $\calL_{2c}$ does NOT have a
square root. Let $I$ be the Iwahori group of $G_\der$
corresponding to the chosen alcove $\mathbf{a}$, $i\in\bold S$ given by $v_0$.
Let $\bbP^1_i=P_i/I$ be the rational line in $\Fl_\s=LG_\der/I$ as
constructed in \S \ref{loop grp and flag var}. It projects to a
rational curve in $\Fl^s$ under $LG_\der/I\to
LG_\der/L^+\calG_{v_0}$ (an explicit description of this rational
line is given in \eqref{rational line}). Then the restriction of
$\calL_{2c}$ to this rational line has degree $2n+1$. Since this
line is contained in any Schubert variety $\Fl^s_\mu$, this
means that $\omega_{\Fl^s_\mu}^{-1}$ is NOT a line bundle, i.e.
$\Fl^s_\mu$ is not Gorenstein.
\end{rmk}

\medskip

Now we turn to the global Schubert varieties. Let
$\calG^s=((\Res_{\tilde{C}/C}(H\times\tilde{C}))^\Ga)^0$ be the
Bruhat-Tits group scheme over $C$ as constructed in \S \ref{grpsch}.
Therefore $\calG^s_{\calO_0}\cong\calG_{v_0}$ is the special
parahoric group scheme for $\calG_{F_0}$ as in \eqref{constr of
special}.

\begin{prop}Assume Theorem \ref{top fiber}. Then Theorem \ref{sch
fiber} holds for $\calG^s$.
\end{prop}
\begin{proof}By Theorem \ref{top fiber}, the support of
$(\bGr_{\calG^s,\mu})_{\tilde{0}}$ is a single Schubert variety.
This is because, when $\calG^s_{\calO_0}=\calG_{v_0}$ is a special
parahoric group scheme, $W^Y={W_0}$ and
${W_0}\setminus\Adm^Y(\mu)/{W_0}$ consists of only one extremal
element in the Bruhat order, namely $t_\mu$ under the projection
$\Adm^Y(\mu)\to{W_0}\setminus\Adm^Y(\mu)/{W_0}$. This proves that
the special fiber of $\bGr_{\calG^s,\mu}$ is irreducible. On the
other hand, we have the affine chart $U_\mu$ which is an affine space over $\tilde{C}$ (see \S \ref{affine
chart}) of $\bGr_{\calG^s,\mu}$ and $(U_\mu)_{\tilde{0}}$ is open in
$(\bGr_{\calG^s,\mu})_{\tilde{0}}$. Therefore, the special fiber of
$\bGr_{\calG^s,\mu}$ is generically reduced. By Hironaka's lemma
again, $\bGr_{\calG^s,\mu}$ is normal over $\tilde{C}$, with special
fiber reduced, indeed isomorphic to $\Fl^s_\mu$.
\end{proof}

\begin{cor}
The global Schubert variety $\bGr_{\calG^s,\mu}$ is normal and
Cohen-Macaulay.
\end{cor}
\begin{proof}The normality follows from Hironaka's lemma. Since $\Fl_\mu^s$ is Cohen-Macaulay, the assertion follows.
\end{proof}

We refer to \S \ref{Frob} for a brief discussion of some facts about Frobenius
splittings.

\begin{prop}\label{special parahoric case}
The variety $\bGr_{\calG^s,\mu}$ is Frobenius split, compatibly with
$(\bGr_{\calG^s,\mu})_{\tilde{0}}$.
\end{prop}
\begin{proof}For simplicity, let us denote $\bGr_{\calG^s,\mu}$ by
$X$. Then $f:X\to\tilde{C}$ is flat and is fiberwise normal and
Cohen-Macaulay (since each $X_y$ is a Schubert variety). Let
$\omega_{X/\tilde{C}}$ be the relative dualizing sheaf on $X$. We know
that $f_*\omega_{X/\tilde{C}}^{1-p}$ is a vector bundle on
$\tilde{C}$ by Corollary \ref{vanishing H1}.

By the construction of \S \ref{Gm action}, the sheaf
$f_*\omega_{X/\tilde{C}}^{1-p}$ is $\bbG_m$-equivariant, and
therefore, we can choose a $\bbG_m$-equivariant isomorphism
\begin{equation}\label{tl}
f_*\omega_{X/\tilde{C}}^{1-p}\cong
H^0(X_{\tilde{0}},\omega_{X_{\tilde{0}}}^{1-p})\otimes\calO_{\tilde{C}},
\end{equation}
where the $\bbG_m$ action on $H^0(X_{\tilde{0}},\omega_{X_{\tilde{0}}}^{1-p})$ comes from the $\bbG_m$-equivariant structure on $\omega_{X_{\tilde{0}}}^{1-p}$.
Let $\sigma\in H^0(X_{\tilde{0}},\omega^{1-p}_{X_{\tilde{0}}})$ be a
$\bbG_m$-invariant section which splits $X_0$ (i.e. $\sigma$ is a
splitting of the natural map $\calO_{X_{\tilde{0}}}\to
F_*\calO_{X_{\tilde{0}}}$, when regarded as a morphism from
$F_*\calO_{X_{\tilde{0}}}\to\calO_{X_{\tilde{0}}}$ via
\eqref{dual}). Such a section always exists by Lemma \ref{Gm split} below. Let
$\sigma\otimes 1$ be a section of $f_*\omega_{X/\tilde{C}}^{1-p}$
via the isomorphism \eqref{tl}. We claim that $\sigma\otimes 1$,
regarded as a morphism
$(F_{X/\tilde{C}})_*\calO_X\to\calO_{X^{(p)}}$ via \eqref{rel dual},
will map $1$ to $1$. In fact, $(\sigma\otimes 1)(1)$ is a
$\bbG_m$-invariant non-zero function on $\bGr_{\calG^s,\mu}$ since
its restriction to $X_{\tilde{0}}$ is non-zero by \eqref{base change
for F-split}. But since all regular functions on $\bGr_{\calG^s,\mu}$
come from $\tilde{C}$, $(\sigma\otimes 1)(1)$ is a
$\bbG_m$-invariant non-zero function on $\tilde{C}$, which must be a
constant. But its restriction to $X_{\tilde{0}}$ is $1$, the claim
follows.

Now, let $(\sigma\otimes 1)\otimes (\frac{u}{du})^{p-1}\in
f_*\omega_{X/\tilde{C}}^{1-p}\otimes\omega_{\tilde{C}}^{1-p}\cong
f_*\omega_X^{1-p}$. By the formula \eqref{Cartier isom} (applied to $\tilde{C}$) and the
commutative diagram \eqref{comm diag}, the proposition follows.
\end{proof}

\begin{lem}\label{Gm split}
Let $X$ be an algebraic variety an algebraically closed field $k$ of positive characteristic with a $\bbG_m$-action. Let $\tau: F_*\calO_X\to \calO_X$ be a splitting map of the inclusion $\calO_X\to F_*\calO_X$. Decompose $\tau=\sum_j \tau_j$ according to the
$\bbG_m$-weights, then $\tau_0$ is also a splitting map.
\end{lem}
\begin{proof}
By definition $1=\tau(1)=\sum_j \tau_j(1)$, where $\tau_j(1)$ is a function on $X$ of weight $j$ under the action of $\bbG_m$. Comparing the weights of both sides, we find that $\tau_0(1)=1$ and $\tau_j(1)=0$ for $j\neq 0$.
\end{proof}

Let $\calG$ be the group scheme with $\calG_{\calO_0}=\calG_{\mathbf{a}}$.  Let
$I=L^+\calG_{\mathbf{a}}$. Observe that the natural projection
$\Gr_{\calG}\to\Gr_{\calG^s}$ induces an isomorphism from
$\Gr_{\calG,\mu}$ to its image in $\bGr_{\calG^s,\mu}$. To see this,
observe that $\Gr_{\calG,\mu}$ is covered by $U_\mu$ and
$\Gr_{\calG,\mu}|_{\tilde{C}^{\circ}}$, both of which map
isomorphically to their images in $\bGr_{\calG^s,\mu}$. We thus
regard $\Gr_{\calG,\mu}$ as an open subscheme of
$\bGr_{\calG^s,\mu}$ under this map. The boundary
$\bGr_{\calG^s,\mu}-\Gr_{\calG,\mu}$ has codimension at least two.
Therefore, we have proven
\begin{cor}\label{aux3}
$\Gr_{\calG,\mu}$ is Frobenius split, compatibly with
$(\Gr_{\calG,\mu})_{\tilde{0}}$.
\end{cor}

\begin{cor}\label{reg functions}
The pullback along $f:\Gr_{\calG,\mu}\to\tilde{C}$ gives an isomorphism $f^*:H^0(\tilde{C},\calO_{\tilde{C}})\stackrel{\sim}{\to} H^0(\Gr_{\calG,\mu},\calO_{\Gr_{\calG,\mu}})$.
\end{cor}

\subsection{Proof of Theorem \ref{Fsplit}}
The goal of this subsection is to prove Theorem \ref{Fsplit}.
Without loss of generality, we can assume that $w\in W_\aff$. Let
$s_i (i\in \bold S)$ be the simple reflections (determined by the
alcove $\mathbf{a}$). Let us recall that for
$\tilde{w}=(s_{i_1},s_{i_2},\ldots,s_{i_m})$ a sequence of
simple reflections corresponding to affine simple roots with $w=s_{i_1}\cdots s_{i_m}$, the
Bott-Samelson-Demazure-Hasen variety is defined as
\[D_{\tilde{w}}=L^+P_{i_1}\times^{I}L^+P_{i_2}\times^{I}\cdots \times^{I}L^+P_{i_n}/I,\]
where $P_i$ is the parahoric group corresponding to $i$ (so that
$L^+P_i/I\cong \bbP^1$). This is a smooth variety which is an iterated
fibration by $\bbP^1$. For any subset
$\{j_1,\ldots,j_n\}\subset\{1,\ldots,m\}$, let
$\tilde{v}=(s_{i_{j_1}},\ldots,s_{i_{j_n}})$ be the corresponding subsequence of
$\tilde{w}$, let $H_{i_p}, p=1,2,\ldots,m$ be defined as
\[H_{i_p}=\left\{\begin{array}{ll}I &\mbox{ if } p\not\in \{j_1,\ldots, j_n\}\\ L^+P_{i_p} &\mbox{ if } p\in \{j_1,\ldots, j_n\}.\end{array}\right.\]
Then there is a closed embedding
$\sigma_{\tilde{v},\tilde{w}}:D_{\tilde{v}}\to D_{\tilde{w}}$ given
by
\begin{multline}
D_{\tilde{v}}=L^+P_{j_1}\times^IL^+P_{j_2}\times^I\cdots\times^IL^+P_{j_n}/I\cong H_{i_1}\times^I\cdots\times^IH_{i_m}/I\\
\hookrightarrow L^+P_{i_1}\times^IL^+P_{i_2}\times^I\cdots
\times^IL^+P_{i_m}/I=D_{\tilde{w}}.
\end{multline}
In particular, let $\tilde{w}[j]$ denote the subsequence of
$\tilde{w}$ obtained by deleting $s_{i_j}$. Then
\[\sigma_{\tilde{w}[j],\tilde{w}}:D_{\tilde{w}[j]}\hookrightarrow D_{\tilde{w}}\]
is a divisor. This way, we obtain $m$ divisors of
$D_{\tilde{w}}$. If $\tilde{v}_1,\tilde{v}_2$ are two subsequences
of $\tilde{w}$, then the scheme-theoretical intersection
$D_{\tilde{v}_1}\cap D_{\tilde{v}_2}$ inside $D_{\tilde{w}}$ is
$D_{\tilde{v}_1\cap\tilde{v}_2}$.

For $w\in W_\aff$, let $m=\ell(w)$, let us fix a reduced
expression of $w=s_{i_1}\cdots s_{i_{m}}$ and let
$\tilde{w}=(s_{i_1},s_{i_2},\ldots,s_{i_m})$. Let $D_{\tilde{w}}$ be
the corresponding BSDH variety so that $D_{\tilde{w}}$ is smooth and
$\pi_{\tilde{w}}:D_{\tilde{w}}\to \Fl_w$ is birational. By twisting
by the $I$-torsor
$\Gr_{\calG,\mu}\times_{\Gr_\calG}\Gr_{\calG,\underline{0}}$, we
have
$\Gr_{\calG,\mu}\tilde{\times}D_{\tilde{w}}\to\Gr_{\calG,\mu}\tilde{\times}\Fl_w$,
still denoted by $\pi_{\tilde{w}}$.

By the standard argument, to prove Theorem \ref{Fsplit}, it is
enough to prove that:
\begin{prop}\label{splitting}
The variety $\Gr_{\calG,\mu}\tilde{\times}D_{\tilde{w}}$ is Frobenius split,
compatibly with all $\Gr_{\calG,\mu}\tilde{\times}D_{\tilde{w}[j]}$ for all $j$, 
and with $(\Gr_{\calG,\mu})_{\tilde{0}}\tilde{\times}D_{\tilde{w}}$.
\end{prop}
Let $\omega_{D_{\tilde{w}}}$ be the canonical sheaf of
$D_{\tilde{w}}$. It is known that there is an isomorphism (for
example, see \cite[Proposition 3.19]{Go1} for the $\SL_n$ case,
\cite[proof of Proposition 9.6]{PR1}, \cite[Ch. 8.18]{Ma} for the
general case)
\begin{equation}\label{canonical line on BSDH}
\omega_{D_{\tilde{w}}}^{-1}\cong\calO(\sum_{j=1}^mD_{\tilde{w}[j]})\otimes\pi_{\tilde{w}}^*\calL(\sum_{i\in\bold
S}\epsilon_i),
\end{equation}
where $\calL(\sum_{i\in\bold S}\epsilon_i)$ is the line bundle on
$\Fl_{\s}$ as defined in \S \ref{loop grp and flag var}. (Recall that
since we assume that $w\in W_\aff$,
$\Fl_w\subset\Fl_\s=(\Fl)^0_{\rm red}$ by \cite[\S 6]{PR1}). If we endow  $\calL(\sum_{i\in\bold S}\epsilon_i)$ with the $I$-equivariant structure such that  $I$ acts on its fiber over $\ast\in\Fl_{\s}$ trivially, then the isomorphism \eqref{canonical line on BSDH}  is $I$-invariant. This observation allows us to formulate a
relative version of this isomorphism. 

Let us denote by $\calL_{2c}$ the line bundle on $\Fl$ which is the
pullback of $\calL_{\det}$ along $\Fl\to\Gr_{\GL(\Lie I)}$ as
in \S \ref{sub;lines on globaff} (as before, by abuse of notation, $\Lie I$ is considered as an $\calO$-module). We endow it with the $I$-equivariant structure so that $I$ acts its fiber over $\ast\in\Fl$ trivially. By twisting by the $I$-torsor $\Gr_{\calG,\mu}\times_{\Gr_\calG}\Gr_{\calG,\underline{0}}$, we obtain
a line bundle on $\Gr_{\calG,\mu}\tilde{\times}\Fl$, still denoted by $\calL_{2c}$. In addition, to
simply the notation, let us denote the projection
$\Gr_{\calG,\mu}\tilde{\times}D_{\tilde{w}}\to\Gr_{\calG,\mu}$ by
$f:X\to V$. Then by the same proof as \eqref{canonical line on BSDH}
(i.e. induction on the length of $w$), we have
\begin{equation}\label{canonical line on BSDH, rel}
\omega_{X/V}^{-2}\cong\calO(2\sum_{j=1}^m\Gr_{\calG,\mu}\tilde{\times}D_{\tilde{w}[j]})\otimes\pi_{\tilde{w}}^*\calL_{2c}.
\end{equation}
\quash{\begin{rmk}Observe that in the above isomorphism we use
$\omega_{X/V}^{-2}$ rather than $\omega_{X/V}^{-1}$. The reason is
as follows. Clearly, $\omega_{D_{\tilde{w}}}^{-1}$ is
$I$-equivariant. If $G=G_\der$ is simply-connected, both lines
$\calO(\sum_{j=1}^mD_{\tilde{w}[j]})$ and $\calL(\sum_{i\in\bold
S}\epsilon_i)$ are also $B$-equivariant and \eqref{canonical line on
BSDH} is naturally upgraded to a $B$-equivariant isomorphism so by
twisting the $B$-torsor
$\Gr_{\calG,\mu}\times_{\Gr_\calG}\Gr_{\calG,\underline{0}}$ we can
obtain a relative version of \eqref{canonical line on BSDH} whose
square gives rise to \eqref{canonical line on BSDH, rel}. However,
if $G$ is not simply-connected, the $B$-equivariant structure on
$\calO(\sum_{j=1}^mD_{\tilde{w}[j]})$ or $\calL(\sum_{i\in\bold
S}\epsilon_i)$ may not exist nor be unique. Therefore the natural
relative version for \eqref{canonical line on BSDH} is
\eqref{canonical line on BSDH, rel}.
\end{rmk}}

We will later prove the following lemma.
\begin{lem}\label{section}
There is a section $\sigma_0$ of $\calL_{2c}$ whose divisor
$\on{div}(\sigma_0)\subset \Gr_{\calG,\mu}\tilde{\times}\Fl_w$ does
not intersect
$z(\Gr_{\calG,\mu})=\Gr_{\calG,\mu}\tilde{\times}\ast$.
\end{lem}
Let us remark that the line bundle $\calL(\sum_{i\in\bold
S}\epsilon_i)$ is very ample on $\Fl_w$, and therefore there exists
a section of $\calL(\sum_{i\in\bold S}\epsilon_i)$ that does not
pass through $\ast$. However, $\calL_{2c}$ is twisted by the
$I$-torsor
$\Gr_{\calG,\mu}\times_{\Gr_\calG}\Gr_{\calG,\underline{0}}$, and it
is not ample. Therefore, some detailed analysis of this line bundle
is needed.

Let us first assume this lemma, and let $\sigma$  be a section of
$\omega_{X/V}^{-2}$ whose divisor is of the form
\begin{equation}\label{splitting section}
\on{div}(\sigma)=2\sum_{j=1}^m\Gr_{\calG,\mu}\tilde{\times}D_{\tilde{w}[j]}+\on{div}(\pi_{\tilde{w}}^*\sigma_0).
\end{equation}
We claim that
\begin{lem}
A non-zero scalar multiple of the section
$\sigma^{\frac{p-1}{2}}\in\omega_{X/V}^{1-p}$ (recall that we assume
that $p>2$), when regarded as a morphism
$(F_{X/V})_*\calO_X\to\calO_{X^{(p)}}$ via \eqref{rel dual}, will
send $1$ to $1$.
\end{lem}
\begin{proof}
Let
\[h=\sigma^{\frac{p-1}{2}}(1)\in\Gamma(X^{(p)},\calO_{X^{(p)}})\] be
the function as in the lemma. 
By Corollary \ref{reg functions}, we have $\Gamma(X^{(p)},\calO_{X^{(p)}})=\Ga(\tilde{C},\calO_{\tilde{C}})$ and so $h$ is obtained by pullback from a function on $\tilde C$ which then has to be a constant. To see $h$ is
non-where vanishing, let $x\in \Gr_{\calG,\mu}$ be a point, and it
is enough to show the restriction of $h$ to
$(D_{\tilde{w}})_x:=\Gr_{\calG,\mu}\tilde{\times}D_{\tilde{w}}|_x\cong
D_{\tilde{w}}$ is not zero. This is because the restriction of
$\sigma$ to $\Gr_{\calG,\mu}\tilde{\times}D_{\tilde{w}}|_x$ gives a
divisor of the form $2\sum_{j=1}^mD_{\tilde{w}[j]}+D$ for some $D$
which does not pass through $\ast$. Therefore, by (a slight
variant of) \cite[Proposition 8]{MR},
$\sigma^{\frac{p-1}{2}}|_{(D_{\tilde{w}})_x}$, when regarded as a
morphism from $F_*\calO_{D_{\tilde{w}}}$ to
$\calO_{(D_{\tilde{w}})_x}$ via \eqref{dual}, will send $1$ to a
non-zero constant function on $(D_{\tilde{w}})_x$. Therefore, by
\eqref{base change for F-split},
\[h|_{(D_{\tilde{w}})_x}=\sigma^{p-1}|_{(D_{\tilde{w}})_x}(1)\]
is a non-zero constant. This finishes the proof of the lemma.
\end{proof}

Now let $\tau\in\omega_{\Gr_{\calG,\mu}}^{1-p}$ be a section which
gives rise to a Frobenius splitting of $\Gr_{\calG,\mu}$, compatible
with $(\Gr_{\calG,\mu})_{\tilde{0}}$ by Corollary \ref{aux3}.
Consider $\sigma^{\frac{p-1}{2}}\otimes f^*\tau\in\omega^{1-p}_X$. By
\eqref{comm diag}, it gives a splitting of
$\Gr_{\calG,\mu}\tilde{\times}D_{\tilde{w}}$, compatible with
$(\Gr_{\calG,\mu})_{\tilde{0}}\tilde{\times}D_{\tilde{w}}$. Again,
by (a slight variant of) \cite[Proposition 8]{MR}, this splitting is
also compatible with all
$\Gr_{\calG,\mu}\tilde{\times}D_{\tilde{w}[j]}$. This finishes the
proof of Theorem \ref{Fsplit}.

\medskip

It remains to prove Lemma \ref{section}. Let us consider the surjective map
\[V_w=\Ga(\Fl_w,\calL_{2c})\to \Ga(\ast,\calL_{2c})=V_1\cong k.\]
By twisting with the $I$-torsor
$\Gr_{\calG,\mu}\times_{\Gr_\calG}\Gr_{\calG,\underline{0}}$, we
obtain a surjective morphism of vector bundles
$\calV_w\to\calV_1\cong\calO_{\Gr_{\calG,\mu}}$ over
$\Gr_{\calG,\mu}$. Clearly, $\calV_w$ is $\pi_*\calL_{2c}$, where $\pi: \Gr_{\calG,\mu}\tilde{\times}\Fl\to \Gr_{\calG,\mu}$ is the base change of $\pi:\Gr_{\calG}^{conv}\to \Gr_\calG$. Then to
prove Lemma \ref{section} is equivalent to prove that there is a
morphism $\calO_{\Gr_{\calG,\mu}}\to\calV_w$ (which determines the
section $\sigma_0$ of $\calL_{2c}$) such that the composition
$\calO_{\Gr_{\calG,\mu}}\to\calV_w\to\calV_1$ is an isomorphism.

To this goal, let us first observe that the $I$-torsor
$\Gr_{\calG,\underline{0}}\times_CC^{\circ}\to\Gr_\calG\times_CC^{\circ}$
has a canonical section. Namely, we associated an $R$-point
$(y,\calE,\beta)$ of $\Gr_\calG\times_CC^{\circ}$ an $R$-point
$(y,\calE,\beta,\ga)$ of
$\Gr_{\calG,\underline{0}}\times_CC^{\circ}$ as follows. Since
the graph $\Ga_y$ of $y:\Spec R\to C$ does not intersect with
$\{0\}\times \Spec R\subset C\times \Spec R$, we can define
\[\ga:\calE|_{\widehat{\{0\}\times \Spec
R}}\to\calE^0|_{\widehat{\{0\}\times \Spec R}}\] as the restriction
of $\beta: \calE|_{C_R-\Ga_y}\cong\calE^0|_{C_R-\Ga_y}$. By base
change, we get a canonical section (a canonical trivialization)
$\psi$ of the $I$-torsor
$W\times_{\Gr_\calG}\Gr_{\calG,\underline{0}}\to W$, where
$W=\Gr_{\calG,\mu}|_{\tilde{C}^{\circ}}\cong \Gr_\mu\times
\tilde{C}^{\circ}$. Therefore,
$\Gr_{\calG,\mu}\tilde{\times}\Fl_w|_W\cong W\times\Fl_w$
canonically, and over $W$, we have
\[\begin{CD}
V_w\otimes\calO_W@>>>V_1\otimes\calO_W\\
@V\cong VV@V\cong VV\\
\calV_w|_W@>>>\calV_1|_W.
\end{CD}\]
To complete the proof of the lemma, it is enough to show
\begin{enumerate}
\item the isomorphism $V_1\otimes\calO_W\to\calV_1|_W$ extends to an
isomorphism $V_1\otimes \calO_{\Gr_{\calG,\mu}}\to\calV_1$;

\item there is a splitting $V_1\to V_w$ (equivalently, a section of $\calL(2\sum_{i\in\bold S}\epsilon_i)$ whose divisor does not pass through
$\ast\in \Fl_w$), such that the induced map \[V_1\otimes\calO_W\to
V_w\otimes\calO_W\to\calV_w|_W\] extends to
$V_1\otimes\calO_{\Gr_{\calG,\mu}}\to\calV_w$.
\end{enumerate}

Let us first prove (1). Let us consider the general situation: Let $E\to B$ be a torsor under some group $K$, and $M$ be a space with the trivial $K$-action. Then there is a canonical isomorphism $t:E\times^KM\simeq E/K\times M=B\times M$. In addition, for any section $s:B\to E$, the induced isomorphism $E\times^K M\simeq (B\times K)\times^{ K} M\simeq B\times M$ coincides with $t$. Back to our situation, as the $I$-module $V_1$ is trivial, we can apply this general remark to conclude that $\calV_1$ is canonically trivialized, which restricts to its trivialization over $W$ induced from the canonical trivialization of the $I$-torsor over $W$.

To prove (2),
let us first complete the curve $\bar{C}=C\cup\{\infty\}\cong\bbP^1$ and $\bar{\tilde{C}}=\tilde{C}\cup\{\tilde{\infty}\}$. We extend $\calG$ to a group scheme over $\bar{C}$ so that $\calG_{\calO_\infty}$ is the pro-unipotent radical of the Iwahori opposite to $\calG_{\calO_0}$. More precisely, the pinning of $H$ (\S \ref{grp data} and \S \ref{grpsch}) determines a unique Borel $B^-$ such that $B_H\cap B^-=T_H$. Let $U^-=[B^-,B^-]$. Let $\tilde{\calG}$ be the group scheme over $\bar{\tilde{C}}$ obtained by dilatation of $H\times\bar{\tilde{C}}$ along $B_H\times\{\tilde{0}\}$ and $U^-\times\{\tilde{\infty}\}$. Then $\calG$ is the neutral connected component of $(\Res_{\bar{\tilde{C}}/\bar{C}}\tilde{\calG})^\Ga$. This group scheme is the same as the group scheme $\calG(0,1)$ in \cite{HNY}. Let $$I^{u,-}=\Gamma(\bar{C}-\{0\},\calG),$$ considered as an ind-group over $k$. Then the Birkhoff decomposition (cf. \cite[Proposition 1(4)]{HNY}) implies that $\Lie G= \Lie I\oplus \Lie I^{u,-}$ as $k$-vector spaces (this is the triangular decomposition in the Kac-Moody theory). For an $\calO$-lattice $L$ in $\Lie G$, consider the the complex of $k$-vector spaces $$L\oplus \Lie I^{u,-}\to \Lie G.$$ As $L$ varies, its determinant defines a section of $\calL_{\det}$ (over the neutral connected component of $\Gr_{\GL(\Lie I)}$), whose pullback defines a section $\sigma^0$ of $\calL_{2c}$ vanishing away from $\ast\in\Fl$. This gives us a
splitting $V_1\to V_w$ which we claim is the desired splitting
satisfying (2).

The prove this claim, we need two more ingredients.
Let
$\Bun_\calG$ be the moduli stack of $\calG$-bundles on $\bar{C}$.
Let us express $\Fl$ as the ind-scheme representing $(\calE,\beta)$,
where $\calE$ is a $\calG$-torsor on $\bar{C}$ and $\beta$ a
trivialization of $\calE$ away from $0\in \bar{C}$. 
Let
$\omega_{\Bun_\calG}^{-1}$ be the anti-canonical bundle of
$\Bun_\calG$. Its fiber over a $\calG$-torsor $\calE$ is the inverse of the determinant of the cohomology $\det R\Gamma(\bbP^1,\ad\calE)^{-1}$. Therefore
$\omega_{\Bun_\calG}^{-1}$ is isomorphic to the pullback along $\Bun_\calG\to\Bun_{\GL(\Lie \calG)}$ of the inverse of the determinant of cohomology line bundle.
As is well-known (e.g. \cite{F}), the pullback of the latter line bundle on $\Bun_{\GL(\calV)}$ to $\Gr_{\GL(\calV)}$ is the determinant line bundle $\calL_{\det}$ we introduced in \S \ref{lines on globaff}.
Therefore, we have
$\calL_{2c}\cong h^*\omega_{\Bun_\calG}^{-1}$. 

The following lemma is first ingredient we need.
\begin{lem}The section $\sigma^0$ of $\calL_{2c}$ descends to a section $\Theta\in\omega_{\Bun_\calG}^{-1}$.
\end{lem}
\begin{proof}
Clearly, the adjoint action of $I^{u,-}$ preserves the determinant of $L\oplus \Lie I^{u,-}\to \Lie G$ up to a scalar. As $I^{u,-}$ has no non-trivial characters,  the left action of $I^{u,-}$ on $\Fl$ preserves $\sigma^0$. As $\Bun_\calG$ is the quotient of $\Fl$ by $I^{u,-}$ (cf. \cite[Proposition 1]{HNY}), $\sigma^0$ descends. 
\end{proof}
By \cite[Corollary 1.2]{HNY}, we can translate $\Theta$ to sections of $\omega_{\Bun_\calG}^{-1}$ over other connected components of $\Bun_\calG$, still denoted by $\Theta$. 

Next, consider
the following morphisms
\[\Bun_\calG\stackrel{h_1}{\leftarrow}\Gr_\calG\stackrel{\pi}{\leftarrow}\Gr_\calG^{Conv}\stackrel{m}{\to}\Gr_{\calG}^{BD}\stackrel{h_2}{\to}\Bun_{\calG}.\]
The second ingredient we need is as follows.
\begin{lem}
Over $\Gr_\calG\tilde{\times}\Fl_w\subset\Gr_\calG^{Conv}$, there is
an isomorphism
\[\calL_{2c}\cong m^*h_2^*\omega_{\Bun_\calG}^{-1}\otimes
\pi^*h_1^*\omega_{\Bun_\calG}.\]
\end{lem}
\begin{proof}Since $\Gr_\calG\tilde{\times}\Fl_w$ is proper over $\Gr_\calG$, by the see-saw principle,
it is enough to show that: (i) for each $x\in\Gr_{\calG}$, the
restrictions of $m^*h_2^*\omega_{\Bun_\calG}^{-1}\otimes
\pi^*h_1^*\omega_{\Bun_\calG}$ and $\calL_{2c}$ to $\Fl_w\subset
\pi^{-1}(x)$ are isomorphic; and (ii) when restricting both line
bundles via the section $z:\Gr_\calG\to\Gr_\calG^{Conv}$, they are isomorphic.

Indeed, recall that over $C^{\circ}$,
$\Gr_{\calG}^{Conv}|_{C^{\circ}}\cong\Gr_\calG^{BD}|_{C^{\circ}}\cong\Gr_\calG|_{C^{\circ}}\times\Fl$,
and over $0\in C(k)$, the morphisms
$(\Gr_\calG)_0\stackrel{\pi}{\leftarrow}(\Gr_\calG^{Conv})_0\stackrel{m}{\to}(\Gr_{\calG}^{BD})_0$
identify with
$\Fl\stackrel{\pi}{\leftarrow}\Fl\tilde{\times}\Fl\stackrel{m}{\to}\Fl$.
Under these isomorphisms
\[h_2^*\omega_{\Bun_\calG}^{-1}|_{\Gr_{\calG}^{BD}|_{C^{\circ}}}\cong h_1^*\omega_{\Bun_\calG}^{-1}|_{\Gr_\calG|_{C^{\circ}}}\otimes h^*\omega_{\Bun_\calG}^{-1},\quad h_2^*\omega_{\Bun_\calG}^{-1}|_{(\Gr_\calG^{BD})_0}\cong h^*\omega_{\Bun_\calG}^{-1}.\]
Therefore, for all $x\in\Gr_\calG$, the restriction of
$m^*h_2^*\omega_{\Bun_\calG}^{-1}\otimes
\pi^*h_1^*\omega_{\Bun_\calG}$ to $\Fl_w\subset \pi^{-1}(x)$ is
isomorphic to $\calL_{2c}$. The first fact is established. For the
second fact, one can easily see that when restricting both line
bundles via $z:\Gr_\calG\to\Gr_\calG^{Conv}$, they are isomorphic to
the trivial bundle.
\end{proof}

Finally, we prove that $\sigma^0$ gives the desired splitting
satisfying (2).
Indeed, since the $I$-torsor
$\Gr_{\calG,\underline{0}}\times_CC^{\circ}\to\Gr_\calG\times_CC^{\circ}$
has a canonical section, we can spread out $\sigma^0$ as a section
of $\calL_{2c}$ over $\Gr_{\calG,\mu}\tilde{\times}\Fl_w|_W$, still
denoted by $\sigma^0$. This induces a map
$V_1\otimes\calO_W\to V_w\otimes\calO_w$. Then to prove (2), it is
equivalent to show that $\sigma^0$ indeed extends to a section of
$\calL_{2c}$ over the whole $\Gr_{\calG,\mu}\tilde{\times}\Fl_w$.
Otherwise, let $n>0$ be the smallest integer such that $u^n\sigma^0$
would extend (recall that we use $u$ to denote the global coordinate
on $\tilde{C}$ so that $u=0$ defines the divisor
$(\Gr_{\calG,\mu})_{\tilde{0}}\tilde{\times}\Fl_w$ inside
$\Gr_{\calG,\mu}\tilde{\times}\Fl_w$). Then
$u^n\sigma^0|_{(\Gr_{\calG,\mu})_{\tilde{0}}\tilde{\times}\Fl_w}$
would not be zero. Observe that by construction, over
$\Gr_{\calG,\mu}\tilde{\times}\Fl_w|_W$, we have
\[\pi^*h_1^*\Theta\otimes\sigma^0=m^*h_2^*\Theta,\]
as sections in
$m^*h_2^*\omega_{\Bun_\calG}^{-1}|_{\Gr_{\calG,\mu}\tilde{\times}\Fl_w|_W}$.
Then as sections in $m^*h_2^*\omega_{\Bun_\calG}^{-1}$ over the
whole $\Gr_{\calG,\mu}\tilde{\times}\Fl_w$, we would have
\[\pi^*h_1^*\Theta\otimes u^n\sigma^0=u^nm^*h_2^*\Theta.\]
When restricting this equation to
$(\Gr_{\calG,\mu})_{\tilde{0}}\tilde{\times}\Fl_w$, the right hand
side is zero. However, the left have side is not since
$\pi^*h_1^*\Theta|_{(\Gr_{\calG,\mu})_{\tilde{0}}\tilde{\times}\Fl_w}\neq
0$. This is a contradiction!

\section{Proofs II: the nearby cycles}\label{PfII}
\subsection{The strategy}
In this section, we prove Theorem
\ref{top fiber}. As mentioned in the introduction, a direct proof
would be to write down a moduli problem $\calM_\mu$ over
$\tilde{C}$, which is a closed subscheme of $\wGr_\calG$, such that:
(i)
$\calM_\mu|_{\tilde{C}^{\circ}}\cong\bGr_{\calG,\mu}|_{\tilde{C}^{\circ}}$;
and (ii)
$(\calM_{\mu})_{\tilde{0}}(k)=\bigcup_{w\in\Adm^Y(\mu)}\Fl_w^Y(k)$.
Then by Lemma \ref{easy}, Theorem \ref{top fiber} would follow.
Unfortunately, so far, such a moduli functor is not available for
general group $G$ and general coweight $\mu$. In certain cases, such
a moduli problem is available. We refer to \cite{PRS} for a survey
of the known results.

The proof presented here is indirect. Let $(S,s,\eta)$ be a
Henselian trait, i.e. $S$ is the spectrum of a discrete valuation
ring, $s$ is the closed point of $S$ and $\eta$ is the generic point
of $S$. Assume that the residue field $k(s)$ of $s$ is algebraically
closed and let $\ell$ be a prime different from $\on{char}k(s)$.
Recall that if $p:\frakX\to S$ is a morphism, where $\frakX$ is a
scheme, (separated) and of finite type over $S$ there is the
so-called nearby cycle functor
\[\Psi_V: \on{D}_c^b(V_\eta,\overline{\bbQ}_\ell)\to \on{D}_c^b(V_s\times_s\eta,\overline{\bbQ}_\ell),\]
which restricts to an exact functor (\cite[Sect.
4]{Il} ) between the categories of perverse sheaves
\[\Psi_V: \on{Perv}(V_\eta,\overline{\bbQ}_\ell)\to \on{Perv}(V_s\times_s\eta,\overline{\bbQ}_\ell).\]
For $V$ a variety over a field whose characteristic prime to $\ell$,
the intersection cohomology sheaf is the Goresky-MacPherson
extension to $V$ of the (shifted) constant sheaf $\bbQ_\ell[\dim V]$
on the smooth locus of $V$. We will use the following lemma.
\begin{lem}\label{support}
Let $f:V\to S$ be a proper flat morphism. Let $\IC$ be the
intersection cohomology sheaf of $V_\eta:=V\times_S\eta$ and let
$\Psi_V(\IC)$ be the nearby cycle of $\IC$. Then the support of
$\Psi_V(\IC)$ is $V_s$.
\end{lem}
\begin{proof}Let $x\in V$ be a point in the special fiber $V_s$ and $\bar{x}$ be a geometric point over $x$.
Then by definition $\Psi_V(\IC)_{\bar{x}}\cong
H^*((V_{(\bar{x})})_{\bar{\eta}},\IC|_{(V_{(\bar{x})})_{\bar{\eta}}})$,
where $V_{(\bar{x})}$ is the strict Henselization of $V$ at
$\bar{x}$, and $(V_{(\bar{x})})_{\bar{\eta}}$ is its fiber over
$\bar{\eta}$, a geometric point over $\eta$. Let $x$ be a generic
point of $V_s$, then $(V_{(\bar{x})})_{\bar{\eta}}$ is the union of
finite many points and
$\IC|_{(V_{(\bar{x})})_{\bar{\eta}}}\cong\bbQ_\ell[\dim V]^m$ for
some $m>0$. The lemma follows.
\end{proof}

Now, let $\ell$ be a prime different from $p$. Let $\IC_\mu$ be the
intersection cohomology sheaf of
$\bGr_{\calG,\mu}|_{\tilde{C}^{\circ}}$. Then the nearby cycle
$\Psi_{\bGr_{\calG,\mu}}(\IC_\mu)$ is a perverse sheaf on $\Fl^Y$
whose support is $(\bGr_{\calG,\mu})_{\tilde{0}}$. Therefore, to
prove the theorem, it is enough to determine the support of
$\Psi_{\bGr_{\calG,\mu}}(\IC_\mu)$. In fact, we will give a
filtration of $\Psi_{\bGr_{\calG,\mu}}(\IC_\mu)$ and describe the
support of each associated graded piece.

When the group $G$ is split, such a description can be deduced from
\cite[Theorem 4]{AB} directly. In the non-split case, we will mostly
follow their strategy but with the following difference. We will not
make use of the results in \cite[Appendix]{B} and therefore we will
not generalize the full version of \cite[Theorem 4]{AB} to the
ramified case (but see Remark \ref{strengthen}). In particular, we
will not perform any categorical arguments as in \emph{loc. cit.}.

\subsection{Central sheaves}\label{central sheaves} 

Let us set $K^Y=L^+\calG_{\sigma_Y}$, and let $\on{P}_{K^Y}(\Fl^Y)$
denote the category of $K^Y$-equivariant perverse sheaves on
$\Fl^Y$. Recall that this category is defined as the direct limit of
categories of $K^Y$-equivariant of perverse sheaves supported on the
$K^Y$-stable finite dimensional subvarieties of $\Fl^Y$ (see
\cite[Appendix]{G} for details).

\begin{lem}\label{equivariance} The sheaf $\Psi_{\bGr_{\calG,\mu}}(\IC_\mu)$ naturally belongs to
$\on{P}_{K^Y}(\Fl^Y)$.
\end{lem}
\begin{proof}Let $\calL^+_n\calG$ be the $n$th jet group of $\calG$, i.e. the group scheme over
$C$, whose $R$-points classify $(y,\beta)$ where $y\in C(R)$ and
$\beta\in\calG(\Ga_{y,n})$, where $\Ga_{y,n}$ is the $n$th nilpotent
thickening of $\Ga_y$. It is clear that $\calL^+_n\calG$ is smooth
over $C$ and the action of $\wJG$ on $\bGr_{\calG,\mu}$ factors
through some $\calL^+_n\calG\times_C\tilde{C}$ for $n$ sufficiently
large.

Let
$m:\calL^+_n\calG\times_{\tilde{C}}\bGr_{\calG,\mu}\to\bGr_{\calG,\mu}$
be the multiplication and $p$ be the natural projection. Then there
is a canonical isomorphism $m^*\IC_\mu\cong p^*\IC_\mu$ as sheaves
on
$\calL^+_n\calG\times_{\tilde{C}}\bGr_{\calG,\mu}|_{\tilde{C}^{\circ}}$.
By taking nearby cycles, we have a canonical isomorphism
\[\Psi_{\calL^+_n\calG\times_{\tilde{C}}\Gr_{\calG,\mu}}(m^*\IC_\mu)\cong\Psi_{\calL^+_n\calG\times_{\tilde{C}}\Gr_{\calG,\mu}}(p^*\IC_\mu).\]
Since both $m$ and $p$ are smooth morphisms and taking nearby cycle
commutes with smooth base change, we have
\begin{equation}\label{equivar}
m^*\Psi_{\bGr_{\calG,\mu}}(\IC_\mu)\cong
p^*\Psi_{\bGr_{\calG,\mu}}(\IC_\mu).\end{equation} The isomorphism
$m^*\IC_\mu\cong p^*\IC_\mu$ satisfies the cocycle condition under
the pullback to
$\calL^+_n\calG\times_{\tilde{C}}\calL^+_n\calG\times_{\tilde{C}}\bGr_{\calG,\mu}|_{\tilde{C}^{\circ}}$.
This implies the cocycle condition for the isomorphism
\eqref{equivar}. The lemma follows.
\end{proof}

Let us define
\begin{equation}\label{Zla}
\calZ_\mu=\Psi_{\bGr_{\calG,\mu}}(\IC_\mu)
\end{equation}
as a $K^Y$-equivariant perverse sheaf on $\Fl^Y$.

Let $D(\Fl^Y)$ be the derived category of constructible sheaves on
$\Fl^Y$ and $D_{K^Y}(\Fl^Y)$ be the $K^Y$-equivariant derived
category on $\Fl^Y$. Recall that $D_{K^Y}(\Fl^Y)$ is a monoidal
category and there is a monoidal action (the ``convolution product")
of $D_{K^Y}(\Fl)$ on $D(\Fl^Y)$ (cf. \cite[Section 4]{MV}). Namely, we have the convolution
diagram
\[\Fl^Y\times\Fl^Y\stackrel{q}{\leftarrow}LG\times\Fl^Y\stackrel{p}{\to}LG\times^{K^Y}\Fl^Y=\Fl^Y\tilde{\times}\Fl^Y\stackrel{m}{\to}\Fl^Y\]
Let $\calF_1\in D(\Fl^Y)$, $\calF_2\in D_{K^Y}(\Fl^Y)$, and let $\calF_1\tilde{\times}\calF_2$ be the unique sheaf on
$LG\times^{K^Y}\Fl^Y$ such that
\begin{equation}\label{twp1}
p^*(\calF_1\tilde{\times}\calF_2)\cong
q^*(\calF_1\boxtimes\calF_2).
\end{equation}
Then
\begin{equation}\label{conv prod}
\calF_1\star\calF_2=m_!(\calF_1\tilde{\times}\calF_2),
\end{equation}
where $m_!$ is the derived pushforward functor with compact support. In general, if $\calF_1,\calF_2$ are perverse sheaves, it is not
necessarily the case that $\calF_1\star\calF_2$ is perverse.
However, we have

\begin{thm}(i)
Let $\calF$ be an arbitrary perverse sheaf on $\Fl^Y$. Then
$\calF\star\calZ_\mu$ is a perverse sheaf on $\Fl^Y$.

(ii) If $\calF\in \on{D}_{K^Y}(\Fl^Y)$, then there is a canonical
isomorphism $c_{\calF}:\calF\star\calZ_\mu\cong\calZ_\mu\star\calF$.
\end{thm}
\begin{rmk}
(i) The isomorphism $c_{\calF}$ is the composition of the
isomorphisms in Proposition \ref{Gait} below.

(ii) In the case when $G=H$ is a split group, this theorem is proved
by Gaitsgory (cf. \cite{G})\footnote{In fact, Part (ii) of the
theorem was proved in \cite{G} under the assumption that $\calF$ is
perverse. I am not sure whether the argument applies to the case
that $\calF$ is an arbitrary object in $D_{K^Y}(\Fl^Y)$.}. The
general case proved below follows his line of argument. Still, we take the
opportunity to spell out all the details for the following reasons.
First, the family $\wGr_\calG$ we use here is in fact different from
Gaitsgory's family which has no obvious generalization to the
ramified groups. On the other hand, this theorem for ramified groups
is used  in \cite{Z2} to establish the geometric Satake
correspondence for ramified groups. Second, the use of the
non-constant group schemes allows us to simplify Gaitsgory's
argument. Namely, we can treat (i) and (ii) in Proposition
\ref{Gait} below equal. This argument is generalized to a mixed
characteristic situation in \cite{PZ}. However, in \cite{G}, the
proof of part (i) of Proposition \ref{Gait} is considerably harder than
the proof of part (ii).

(iii) To simplify the notation, in the proof we only consider $Y=\mathbf{a}$
being an alcove. In this case, we denote by $I=K^{\mathbf{a}}$ the
corresponding Iwahori subgroup of $LG$, and denote $\Fl=\Fl^{\mathbf{a}}$.
However, the proof (with the only change by replacing $I$ by $K^Y$
and $\Fl$ by $\Fl^Y$) is valid in \emph{any} parahoric case.
\end{rmk}

\quash{(ii)The isomorphism $c_\calF$ in the above theorem satisfies
the following property: for $\calF_1,\calF_2\in D_I(\Fl)$, the
following hexagon commutes.
\[\begin{CD}(\calF_1\star\calF_2)\star\calZ_\mu@>\cong>>\calF_1\star(\calF_2\star\calZ_\mu)@>\cong>>\calF_1\star(\calZ_\mu\star\calF_2)\\
@V\cong VV@.@V\cong VV\\
\calZ_\mu\star(\calF_1\star\calF_2)@>\cong>>\calZ_\mu\star
\calF_1)\star\calF_2@>\cong
>>(\calF_1\star\calZ_\mu)\star\calF_2
\end{CD}\]
For a proof of this fact when $G=H$, see ???. The general case is the same. The means that we have a constructed a central functor (\bf See Bezrukavnikov \rm)
\[\calZ:\on{P}_{L^+H}(\Gr_H)\to \on{P}_{B}(\Fl)\subset
D_{B}(\Fl).\]}

\begin{proof}Recall the Beilinson-Drinfeld Grassmannian
$\Gr_{\calG}^{BD}$ as introduced in \eqref{BD Grass}. We have
\[\Gr_{\calG}^{BD}\times_C{\tilde{C}^{\circ}}\cong \Fl\times(\Gr_{\calG}\times_C\tilde{C}^{\circ}).\]
For $\calF\in D(\Fl)$, let
\[\calF\boxtimes\IC_\mu\subset D(\Fl\times(\Gr_{\calG}\times_C\tilde{C}^{\circ})),\]
which can be therefore regarded as a complex on
$\Gr_{\calG}^{BD}\times_C{\tilde{C}^{\circ}}$. Consider the
nearby cycle functor $\Psi_{\Gr_{\calG}^{BD}\times_C\tilde{C}}$.

\begin{prop}\label{Gait}
(i)If $\calF\in \on{D}(\Fl)$, there is a canonical isomorphism
\[\Psi_{\Gr_{\calG}^{BD}\times_C\tilde{C}}(\calF\boxtimes\IC_\mu)\cong\calF\star\calZ_\mu.\]

(ii) If $\calF\in D_I(\Fl)$, then there is a canonical isomorphism
\[\Psi_{\Gr_{\calG}^{BD}\times_C\tilde{C}}(\calF\boxtimes\IC_\mu)\cong\calZ_\mu\star\calF.\]
\end{prop}
It is clear that this proposition will imply the theorem. The
isomorphisms involved in the statement essentially come from the
fact that nearby cycles commute with the proper pushforward and the
smooth pullback. They will be constructed in the proof.

We first prove (ii). Let $\Gr_{\calG}^{Conv}$ be the convolution
Grassmannian as introduced in \eqref{Conv Grass}, which we recall is
a fibration over $\Gr_\calG$ with fibers isomorphic to $\Fl$. Regard
$\calF\boxtimes\IC_\mu$ as a complex of sheaves on
$\Gr_{\calG}^{Conv}\times_C\tilde{C}^{\circ}\cong\Fl\times(\Gr_\calG\times_C\tilde{C}^{\circ})$.
Since taking nearby cycles commutes with proper push-forward, it is
enough to prove that as complex of sheaves on
$\Fl\tilde{\times}\Fl$, there is a canonical isomorphism
\[\Psi_{\Gr_{\calG}^{Conv}\times_C\tilde{C}}(\calF\boxtimes\IC_\mu)\cong\calZ_\mu\tilde{\times}\calF,\]
where $\calZ_\mu\tilde{\times}\calF$ is the twisted product as
defined in \eqref{twp1}.

Recall the $I$-torsor $\Gr_{\calG,\underline{0}}$ over $\Gr_{\calG}$
defined in \eqref{B-torsor} and
$\Gr_{\calG}^{Conv}\cong\Gr_{\calG,\underline{0}}\times^I\Fl$. Let
$V\subset\Fl$ be the support of $\calF$, and
$I_n=L^+_n\calG_{\calO_0}$ (the $n$th jet group as defined in the
proof of Lemma \ref{equivariance}) be the finite dimensional
quotient of $I$ such that the action of $I$ on $V$ factors through
$I_n$. Let $\Gr_{\calG,0,n}$ be the $I_n$-torsor over $\Gr_\calG$
which classifies $(y,\calE,\beta,\ga)$ where $(y,\calE,\beta)$ is as
in the definition of $\Gr_\calG$ and $\ga$ is a trivialization of
$\calE$ on the $n$th infinitesimal neighborhood of $0\in C$. Then
$\IC_\mu\tilde{\times}\calF$ is supported on
\[(\tilde{C}\times_C\Gr_{\calG,\underline{0}})\times^IV\cong(\tilde{C}\times_C\Gr_{\calG,0,n})\times^{I_n}V\subset\Gr_\calG^{Conv}\times_C\tilde{C}.\]

Observe that over $\tilde{C}^{\circ}$, it makes sense to talk
about $\IC_\mu\tilde{\times}\calF$ (as defined via \eqref{twp1}),
which is canonically isomorphic to $\calF\boxtimes\IC_\mu$, we thus
need to show that
\begin{equation}\label{ser4}\Psi_{\Gr_{\calG}^{Conv}\times_C\tilde{C}}(\IC_\mu\tilde{\times}\calF)\cong\calZ_\mu\tilde{\times}\calF.\end{equation}
Let us denote the pullback of $\IC_\mu$ to
$\Gr_{\calG,0,n}\times_C\tilde{C}^{\circ}$ by
$\widetilde{\IC}_\mu$. Since $\Gr_{\calG,0,n}\to \Gr_{\calG}$ is
smooth,
$\Psi_{\Gr_{\calG,0,n}\times_C\tilde{C}}(\widetilde{\IC}_\mu)$ is
canonically isomorphic to the pullback of $\calZ_\mu$, and
\[\Psi_{(\Gr_{\calG,0,n}\times_C\tilde{C})\times V}(\widetilde{\IC}_\mu\boxtimes\calF)\cong\Psi_{\Gr_{\calG,0,n}\times_C\tilde{C}}(\widetilde{\IC}_\mu)\boxtimes\calF\]
is $I_n$-equivariant. We thus have \eqref{ser4}.

Next we prove (i). There is another convolution affine Grassmannian
$\Gr_\calG^{Conv'}$, which is an ind-scheme ind-proper over $C$ and represents the functor that associates to
every $k$-algebra $R$,
\begin{equation}
\Gr^{Conv'}_\calG(R)=\left\{(y,\calE,\calE',\beta,\beta') \left|\
\begin{split}&y\in C(R), \calE,\calE' \mbox{ are two }
\calG\mbox{-torsors}\mbox{ on } C_R,
\\&\beta:\calE|_{(C-\{0\})_R}\cong
\calE^0|_{(C-\{0\})_R} \mbox{ is a } \mbox{trivialization}, \\
&\mbox{and }\beta':\calE'|_{C_R-\Ga_y}\cong\calE|_{C_R-\Ga_y}
\end{split}\right.\right\}.
\end{equation}
Let us sketch the proof of the ind-representability of $\Gr^{Conv'}_{\calG}$. Let $\calL^+_n\calG$ be the $n$th jet group of $\calG$. As mentioned
before, $\calL^+_n\calG$ is smooth over $C$. Then one can present $\Gr_\calG$ as the inductive limit $\underrightarrow{\lim} Z_n$ where $Z_n$ is a $\calL^+\calG$-stable closed subscheme and the action of $\calL^+\calG$ on $Z_i$ factors through $\calL^+_n\calG$. Let us define the $\calL^+_n\calG$-torsor $\calP_n$ over $\Fl\times
C$ as follows. Its $R$-points are quadruples $(y,\calE,\beta,\ga)$,
where $y\in C(R)$, $(\calE,\beta)$ are as in the definition of $\Fl$
(and therefore $\beta$ is a trivialization of $\calE$ on
$C^{\circ}_R$), and $\ga$ is a trivialization of $\calE$ over
$\Ga_{y,n}$, the $n$th nilpotent thickening of the graph $\Ga_y$ of
$y$. Then it is not hard to see that $\Gr^{Conv'}_\calG=\underrightarrow{\lim} \calP_n\times^{\calL^+_n\calG}Z_n$ is an ind-scheme ind-proper over $C$.

Clearly, we have $m':\Gr^{Conv'}_\calG\to\Gr_\calG^{BD}$ by sending
$(y,\calE,\calE',\beta,\beta')$ to $(y,\calE',\beta'\circ\beta)$.
This is a morphism over $C$, which is an isomorphism over $C-\{0\}$,
and $m'_{0}$ again is the local convolution diagram
\[m:\Fl\tilde{\times}\Fl\to\Fl.\]
Again, regarding $\calF\boxtimes\IC_\mu$ as a sheaf on
$\Gr^{Conv'}_\calG|_{\tilde{C}^{\circ}}\cong\Fl\times(\Gr_\calG\times_C\tilde{C}^{\circ})$,
it is enough to prove that as sheaves on $\Fl\tilde{\times}\Fl$,
\[\Psi_{\Gr_{\calG}^{Conv'}\times_C\tilde{C}}(\calF\boxtimes\IC_\mu)\cong\calF\tilde{\times}\calZ_\mu.\]

Observe that the action of $\wJG$
on $\bGr_{\calG,\mu}$ factors through some
$\calL^+_n\calG\times_C\tilde{C}$ for $n$ sufficiently large.
Then we have the twisted product
\[(\calP_n\times_C\tilde{C})\times^{\calL^+_n\calG\times_C\tilde{C}}\bGr_{\calG,\mu}\subset\Gr_\calG^{Conv'}\times_C\tilde{C}.\]
Over the restriction of this ind-scheme to $\tilde{C}^{\circ}$, we can form the twisted product
$\calF[1]\tilde{\times}\IC_\mu$ as in \eqref{twp1}, which is
canonically isomorphic to $\calF\boxtimes\IC_\mu$. By the same
argument as in the proof of (ii) (i.e. by pulling back everything to
$\calP_n\times_{\tilde{C}}\bGr_{\calG,\mu}$), we have
\[\Psi_{(\calP_n\times_C\tilde{C})\times^{\calL^+_n\calG\times_C\tilde{C}}\bGr_{\calG,\mu}}(\calF[1]\tilde{\times}\IC_\mu)\cong\Psi_{\Fl\times\tilde{C}}(\calF[1])\tilde{\times}\Psi_{\bGr_{\calG,\mu}}(\IC_\mu)\cong\calF\tilde{\times}\calZ_\mu.\]
\end{proof}

\subsection{Wakimoto filtrations}\label{Waki fil}

Our goal to prove that the support of $\calZ_\mu$ is exactly the
Schubert varieties in $\Fl^Y$ labeled by the set
$W^Y\setminus\Adm^Y(\mu)/W^Y$, which will imply Theorem \ref{top
fiber} by Lemma \ref{support}. Clearly, it is enough to prove this
in the case $\calG_{\calO_0}$ is Iwahori.

Let us recall some standard objects in $\on{P}_I(\Fl)$. Recall that
$I$-orbits in $\Fl$ are labeled by elements $w\in\widetilde{W}$. For
any $w$, let $j_w:C(w)\to\Fl_w$ be the open embedding of the
Schubert cell to the Schubert variety. This is an affine embedding.
Let us denote
\[j_{w*}=(j_w)_*\overline{\bbQ}_\ell[\ell(w)],\quad\quad j_{w!}=(j_w)_!\overline{\bbQ}_\ell[\ell(w)].\]
Then it is well-known (e.g. \cite[Lemma 8]{AB}) that there are
canonical isomorphisms
\begin{equation}
\begin{gathered}\label{Iwahori-Hecke}
j_{w*}\star j_{w'*}\cong j_{ww'*},\quad j_{w!}\star j_{w'!}\cong j_{ww'!}, \quad \mbox{ if } \ell(ww')=\ell(w)+\ell(w'),\\
j_{w*}\star j_{w^{-1}!}\cong j_{w^{-1}!}\star j_{w*}\cong\delta_e.
\end{gathered}
\end{equation}
In addition, if $\ell(ww'w'')=\ell(w)+\ell(w')+\ell(w'')$, then the two isomorphisms from $j_{w*}\star j_{w'*}\star j_{w''*}$ (resp. from $j_{w!}\star j_{w'!}\star j_{w''!}$)
to $j_{ww'w''*}$ (resp. to $j_{ww'w''!}$) are the same.

Let us recall the following fundamental result due to I.Mirkovic
(cf. \cite[Appendix]{AB}). The proof for ramified groups is exactly
the same as for split groups. In fact, the proof works in the
general affine Kac-Moody setting.

\begin{prop}\label{Waki}Let $w,v\in\widetilde{W}$. Then
both $j_{w*}\star j_{v!}$ and $j_{w!}\star j_{v*}$ are perverse
sheaves. In addition, both sheaves are supported on the Schubert
variety $\Fl_{wv}$ and $j_{wv}^*(j_{w*}\star j_{v!})\cong
j_{wv}^*(j_{w!}\star j_{v*})\cong\overline{\bbQ}_\ell[\ell(wv)]$.
\end{prop}

Fix $w\in{W_0}$ to be an element in the finite Weyl group of $G$. We
are going to define the $w$-Wakimoto sheaves on $\Fl$. Recall the
definition of $\xcoch(T)_\Ga^+$ in \eqref{plus}. For
$\mu\in\xcoch(T)_\Ga$, we write $\mu=\la-\nu$ with $\la,\nu\in
w(\xcoch(T)_\Ga^+)$. Define
\begin{equation}\label{wakimoto}
J^w_{\mu}=j_{t_\la !}\star j_{t_\nu*},
\end{equation}
which is well-defined up to a canonical isomorphism (by
\eqref{Iwahori-Hecke}). By Proposition \ref{Waki}, $J^w_\mu\in
\on{P}_I(\Fl)$ and is supported on $\Fl_\mu$ with $j_{t_\mu}^*
J^w_\mu\cong\bar{\bbQ}_\ell[\ell(t_\mu)]$. Let us remark that for
$G$ being split and $w=w_0$ being the longest element in $W_0$, they
are the Wakimoto sheaves considered in \cite{AB}. In addition, we
have
\begin{equation}\label{sum}
J^w_\mu\star J^w_\la\cong J^w_{\mu+\la}.
\end{equation}
In fact, by \eqref{Iwahori-Hecke} and Lemma \ref{length}, this is
true for $\mu,\la$ for $\mu,\la\in w(\xcoch(T)^+_\Ga)$. The
extension to all $\mu,\la$ is immediate.

One of the important applications of the Wakimoto sheaves is as
follow. An object $\calF\in\on{P}_I(\Fl)$ is called convolution
exact if $\calF'\star\calF$ is perverse for any
$\calF'\in\on{P}_I(\Fl)$, and is called central if in addition
$\calF\star\calF' \cong\calF'\star\calF$. For example, $\calZ_\mu$
is central. The following proposition generalizes \cite[Proposition
5]{AB}, where the case $w=e$ is considered. The proof is basically
the same.

\begin{prop}\label{filtration}
Fix $w\in{W_0}$. Any central object in $\on{P}_I(\Fl)$ has a
filtration whose associated graded pieces are $J^w_\la,
\la\in\xcoch(T)_\Ga$.
\end{prop}
\begin{proof}
We begin with some general notations and remarks following \cite{AB}. For a triangulated category $D$ and a
set of objects $S \subset Ob(D)$, let $\langle S\rangle$ be the set of all objects obtained from elements
of $S$ by extensions; i.e. $\langle S\rangle$ is the smallest subset of $Ob(D)$ containing $S\cup\{0\}$ and
such that:
\begin{enumerate}
\item if $A\cong B$ and $A\in\langle S\rangle$, then $B\in\langle
S\rangle$; and
\item for all $A, B\in \langle S\rangle$ and an
exact triangle $A\to C\to B\to A[1]$, we have $C \in \langle
S\rangle$.
\end{enumerate}
Let $\calF\in D_I(\Fl)$. The $*$-support of $\calF$ is defined to be
\[W_\calF^*:=\{w\in\widetilde{W}\mid j_w^*\calF\neq 0\},\] and the $!$-support
of $\calF$ is the set
\[W_\calF^!:=\{w\in\widetilde{W}\mid j_{w}^!\calF\neq 0\}.\]
By the induction on the dimension of the support of $\calF$, it is
easy to see that if $\calF\in D_I(\Fl)^{p,\leq 0}$ ($p$ stands for
the perverse t-structure), then $\calF$ is contained in $\langle
j_{v!}[n]\mid v\in W^*_\calF, n\geq 0\rangle$. On the other hand, if
$\calF\in D_I(\Fl)^{p,\geq 0}$, then $\calF\in\langle j_{v*}[n]\mid
v\in W^!_\calF, n\leq 0\rangle$.

For any $\calF\in D_I(\Fl)$, there exists a finite subset
$S_\calF\subset\widetilde{W}$, such that
\[W^!_{j_{w*}\star\calF}, W^*_{j_{w!}\star\calF}\subset w\cdot S_\calF;\quad\quad W^!_{\calF\star j_{w*}},W^*_{\calF\star j_{w!}}\subset S_\calF\cdot w.\]
Namely, let $\Fl_v$ be a Schubert variety such that $\calF$ is
supported in $\Fl_v$ (in both the $*$-sense and the $!$-sense). Then by the
proper base change theorem, the above assertions will follow if we
can show that there exists $S_v\subset\widetilde{W}$ such that
\[C(w)\tilde{\times}\Fl_v\subset \bigcup_{v'\in wS_v}C(v'),\quad\quad \Fl_v\tilde{\times}C(w)\subset\bigcup_{v'\in S_vw}C(v'),\]
This can be proved easily by induction of the length of $v$.

Now we prove the proposition. Let $\calF\in\on{P}_I(\Fl)$ be a
central object, and let $S_\calF\subset\widetilde{W}$ be the finite
set associated to $\calF$ as above. Recall that we have the special
vertex $v_0$ in the building of $G$, which determines an isomorphism
$\widetilde{W}=\xcoch(T)_\Ga\rtimes{W_0}$ determined by $v_0$. Let
$\mu\in w(\xcoch(T)_\Ga^+)$ such that
\[t_\mu S_\calF\subset w(\xcoch(T)_\Ga^{++}){W_0}, \quad\quad S_\calF t_\mu\subset{W_0}w(\xcoch(T)_\Ga^{++}),\]
where $\xcoch(T)_\Ga^{++}$ is the subset of regular elements in $\xcoch(T)_\Ga^+$. 
This is always possible since $S_\calF$ is a finite set. We have
$J^w_\mu=j_{\mu!}$ and from $J^w_\mu\star\calF\cong\calF\star
J^w_\mu$, we have
\[W^*_{J^w_\mu\star\calF}\subset t_\mu S_\calF\cap S_\calF t_\mu\subset w(\xcoch(T)_\Ga^{++}){W_0}\cap {W_0}w(\xcoch(T)_\Ga^{++})=w(\xcoch(T)_\Ga^{++}).\]
Therefore, $J^w_\mu\star\calF\in\langle j_{t_\la !}[n]\mid\la\in
w(\xcoch(T)_\Ga^+), n\geq 0\rangle$. Observe that $J^w_\la=j_{t_\la
!}$ for $\la\in w(\xcoch(T)_\Ga^+)$. Then by \eqref{sum}, we have
\[\calF\in\langle J^w_\la[n]\mid\la\in\xcoch(T)_\Ga, n\geq 0 \rangle.\]
By choosing $\mu\in w(-\xcoch(T)_\Ga^+)$ large enough and using
$J^w_\la=j_{t_\la *}$ for $\la\in w(-\xcoch(T)^+_\Ga)$,  we have
\begin{equation}\label{FF}\calF':=j_{t_\mu*}\star\calF=J^w_\mu\star\calF\in\langle
j_{t_\la*}[n]\mid\la\in w(-\xcoch(T)^+_\Ga), n\geq 0\rangle.
\end{equation}
We claim that this already implies that $\calF'$ has a filtration
(in the category of perverse sheaves) with associated graded by
$j_{t_\mu*},\mu\in w(-\xcoch(T)^+_\Ga)$, and therefore implies the
proposition. Indeed, since $\calF'$ is perverse, for any $\nu\in
w(-\xcoch(T)_\Ga^+)$, the $!$-stalk of $\calF'$ at $t_\nu$ has
homological degree $\geq -\ell(t_\nu)$. On the other hand,
\eqref{FF} implies that the $!$-stalk of $\calF'$ at $t_\nu$ has
homological degree $\leq -\ell(t_\nu)$. The claim follows.
\end{proof}

To proceed, we now study the category of perverse sheaves on $\Fl$
that are generated by those $J^w_\la$.

\begin{lem}\label{Rhom}
For $\la,\mu\in\xcoch(T)_\Ga$, $R\Hom(J^w_{\la},J^w_{\mu})=0$ unless
$w^{-1}(\la)\preceq w^{-1}(\mu)$. Furthermore,
$R\Hom(J^w_{\mu},J^w_{\mu})\cong\overline{\bbQ}_\ell$.
\end{lem}
\begin{proof}
$R\Hom(J^w_{\la},J^w_{\mu})=R\Hom(J^w_{\la+\nu},J^w_{\mu+\nu})=R\Hom(j_{t_{\la+\nu}!},j_{t_{\mu+\nu}!})$
for $\nu\in\xcoch(T)_\Ga$ such that $\la+\nu,\mu+\nu\in
w(\xcoch(T)_\Ga^+)$. The above complex of
$\ell$-adic vector spaces is non-zero only if $\Fl_{t_{\la+\nu}}\subset \Fl_{t_{\mu_\nu}}$, i.e. $t_{\la+\nu}\leq
t_{\mu+\nu}$ in the Bruhat order. This is equivalent to
$t_{w^{-1}(\la+\nu)}\leq t_{w^{-1}(\mu+\nu)}$ by Lemma \ref{other chamber}, which is in turn
equivalent to $w^{-1}(\la+\nu)\preceq w^{-1}(\mu+\nu)$ by Lemma
\ref{order}, which is equivalent to $w^{-1}(\la)\preceq
w^{-1}(\mu)$. The second statement follows from
$R\Hom(J^w_\mu,J^w_\mu)\cong
R\Hom(J^w_0,J^w_0)\cong\overline{\bbQ}_\ell$.
\end{proof}

\begin{lem}\label{coh of waki}
Let $\calF\in D_I(\Fl)$. Then for any $\mu\in\xcoch(T)_\Ga$,
\[H^*(\Fl,J^w_\mu\star\calF)\cong
H^{*-(w^{-1}(\mu),2\rho)}(\Fl,\calF).\] In particular,
$H^*(\Fl,J^w_\mu)=H^{(w^{-1}(\mu),2\rho)}(\Fl,J^w_\mu)\cong\bar{\bbQ}_\ell$.
\end{lem}
\begin{proof}
For any $v\in\widetilde{W}$,
let $C(v)$ be the Schubert cell in $\Fl$ corresponding to $v$. Then
we have $m:C(v)\tilde{\times}\Fl\to\Fl$, which is an affine bundle over $\Fl$.  Then the isomorphism $j_{v*}\star\calF\cong
m_*(\overline{\bbQ}_\ell[\ell(v)]\tilde{\times}\calF)$ induces $H^*(\Fl,j_{v*}\star\calF)\cong H^*(\Fl,\calF)[\ell(v)]$. Therefore, for $\mu\in w(-\xcoch(T)_\Ga^+)$, the lemma
holds by the above fact and Lemma \ref{length}. If the lemma holds for $\la,\mu$, then
\[H^*(\Fl,\calF)\cong H^*(\Fl,J^w_\la\star J^w_{-\la}\star\calF)\cong H^{*-(w^{-1}(\la),2\rho)}(\Fl, J^w_{-\la}\star\calF),\]
\[H^*(\Fl,J^w_{\la+\mu}\star\calF)\cong H^{*-(w^{-1}(\la),2\rho)}(\Fl,J^w_\mu\star\calF)\cong H^{*-(w^{-1}(\la+\mu),2\rho)}(\Fl,\calF).\]
Therefore, the lemma holds for $-\la$ and $\la+\mu$. Now any element
in $\xcoch(T)_\Ga$ can be written as $\la-\mu$ with $\la,\mu\in
w(\xcoch(T)_\Ga^+)$. We are done.
\end{proof}

Let $\on{W}^w(\Fl)$ be the full abelian subcategory of
$\on{P}_I(\Fl)$ generated by those $J^w_{\mu}, \mu\in\xcoch(T)_\Ga$.
Let $\on{W}^w(\Fl)_{\succeq \mu}$ be the category of $\on{W}^w(\Fl)$
generated by $J^w_{\la}, w^{-1}(\la)\succeq w^{-1}(\mu)$. For each
object $\calF\in W^w(\Fl)$, we define a filtration
\[\calF=\bigcup_{\mu}\calF^w_{\succeq\mu},\]
where $\calF^w_{\succeq\mu}\in\on{W}^w(\Fl)_{\succeq\mu}$ is the
maximal subobject of $\calF$ belonging to
$\on{W}^w(\Fl)_{\succeq\mu}$. Then by Lemma \ref{Rhom},
\[\calF^w_{\succeq\mu}/\bigcup_{w^{-1}(\mu')\succ w^{-1}(\mu)}\calF^w_{\succeq\mu'}\cong
J^w_\mu\otimes {^wW_\calF^\mu},\] where ${^wW_\calF^\mu}$ is a
finite dimensional $\overline{\bbQ}_\ell$ vector space. A direct consequence of
Lemma \ref{coh of waki} is

\begin{cor}Suppose that the notations are as above. Then for any $\calF\in\on{W}^w(\Fl)$, we have
\[H^*(\Fl,\calF)\cong\bigoplus_{\mu\in\xcoch(T)_\Ga}H^*(\Fl,J^w_\mu)\otimes {^wW_\calF^\mu}.\]
\end{cor}

\subsection{Proof of Theorem \ref{top fiber}}\label{proof of thm top fiber} 

Finally, let us prove Theorem \ref{top fiber}. Let
$\mu\in\xcoch(T)_\Ga^+$. Let $\on{Supp}(\mu)$ denote the subset of
$\widetilde{W}$ consisting of those $w$ such that $\Fl_w\subset
(\bGr_{\calG,\mu})_{\tilde{0}}$. We need to show that
$\on{Supp}(\mu)=\Adm(\mu)$. We already know that
$\Adm(\mu)\subset\on{Supp}(\mu)$ (Lemma \ref{easy}).
By Proposition
\ref{filtration} and \ref{Waki}, we also know that the maximal
elements in $\on{Supp}(\mu)$ (under the Bruhat order) belong to
$\xcoch(T)_\Ga\subset\widetilde{W}$. Let $t_{\mu'}\in\on{Supp}(\mu)$
be a maximal element. Then there exists some $w\in{W_0}$ such that
$\mu'\in w(\xcoch(T)_\Ga^+)$. By Proposition \ref{filtration},
$\calZ_\mu\in\on{W}^w(\Fl)$. Write
$\calZ_\mu=\cup_{\la}(\calZ_\mu)^w_{\succeq\la}$ so that the
associated graded pieces are $J^w_\la\otimes {^wW_\mu^\la}$ as above
(we write ${^wW^\la_\mu}$ instead of ${^wW^\la_{\calZ_\mu}}$ for
brevity). By Lemma \ref{support}, ${^wW^{\mu'}_\mu}\neq 0$. In addition,
being a maximal element in $\on{Supp}(\mu)$, $t_{\mu'}$ must have
length $(\mu,2\rho)$. Therefore, $(w^{-1}(\mu'),2\rho)=(\mu,2\rho)$.
On the other hand, $t_{w(\mu)}\in\Adm(\mu)\subset\on{Supp}(\mu)$ is also a maximal element in
$\on{Supp}(\mu)$ since $\ell(t_{w(\mu)})=(2\rho,\mu)$ by Lemma \ref{length}. Therefore, ${^wW^{w(\mu)}_\mu}\neq 0$. We claim that
$\mu'=w(\mu)$. Otherwise, we would have
\[H^{(\mu,2\rho)}(\Fl,\calZ_\mu)\supset {^wW_\mu^{\mu'}}\oplus{^wW_\mu^{w(\mu)}}\]
whose dimension would be at least two.

On the other hand, the map $f:\bGr_{\calG,\mu}\to\tilde{C}$ is
proper, and therefore
$H^*(\Fl,\calZ_\mu)\cong\Psi_{\tilde{C}}(f_*\IC_\mu)$. Since
$\bGr_{\calG,\mu}|_{\tilde{C}^{\circ}}\cong\bGr_\mu\times\tilde{C}^{\circ}$,
we have $H^*(\Fl,\calZ_\mu)\cong IH^*(\bGr_\mu)$ where $IH^*$
denotes the intersection cohomology of $\bGr_\mu$. It is well-known
(for example see \cite{MV}) that
$IH^{(\mu,2\rho)}(\bGr_\mu)\cong\bbQ_\ell$, which contradicts the above unless $\mu'=w(\mu)$. In other words, all the maximal elements
on $\on{Supp}(\mu)$ are contained in $\Adm(\mu)$, which proves the theorem.

\begin{rmk}\label{strengthen}
One should be able to generalize \cite[Theorem 4]{AB} to
the ramified case, which will imply Theorem \ref{top fiber}
directly. We sketch here a possible approach. First, $\on{W}^w(\Fl)$ is indeed a monoidal abelian
subcategory of $P_I(\Fl)$ because $J^w_{\la}\star J^w_{\mu}\cong
J^w_{\la+\mu}$. Let $\on{GrW}^w(\Fl)$ be the submonoidal category
whose objects are direct sums of $J^w_\la$. One can see that this category
is equivalent to $\on{Rep}(\hat{T}^\Ga)$, where $\hat{T}$ is the
dual torus of $T$ defined over $\bar{\bbQ}_\ell$, and $\hat{T}^\Ga$
is the Galois fixed subgroup. By taking the associated graded of the
filtration of $\calF\in\on{W}^w(\Fl)$ defined before,  one obtains a
well-defined functor $\Gr:\on{W}^w(\Fl)\to\on{GrW}^w(\Fl)$. As
explained in \cite[Lemma 16]{AB}, this is a monoidal functor.

Since $\Gr_{\calG}|_{\tilde{C}^{\circ}}\cong
\Gr_H\times\tilde{C}^{\circ}$, the nearby cycle functor indeed
gives a monoidal functor from
$\calZ:\on{P}_{L^+H}(\Gr_H)\to\on{W}^w(\Fl)$, where
$\on{P}_{L^+H}(\Gr_H)$ is the category of $L^+H$-equivariant
perverse sheaves on $\Gr_H$, which is well-known to be equivalent to
the category of representations of the Langlands dual group
$\hat{H}$. One can use the similar argument proved by Gaitsgory in
\cite[Appendix]{B} to show that this functor is in fact central (see
Section 2 of \emph{loc. cit.} for the definition). Then by the same
argument as \cite{AB}, one can show that
$\Gr\circ\calZ:\on{P}_{L^+H}(\Gr_H)\to\on{GrW}^w(\Fl)$ is in fact a
tensor functor, which indeed equivalent to the restriction functor
from the representations of $\hat{H}$ to the representations of
$\hat{T}^\Ga$.
\end{rmk}

\begin{rmk}
This remark is independent of the paper.
As being a nearby cycle, $\calZ_\mu$ carries on the monodromy action of the Galois group of $F_0$. One can show that this action is purely unipotent (see \cite{G} for the case when $G$ is split and \cite[Theorem 10.9]{PZ} in general).
\end{rmk}

\section{Appendix I: line bundles on the local models for ramified unitary groups}\label{lines on local models}
Since Theorem \ref{MainI} is not quite identical to the original
coherence conjecture given by Pappas and Rapoport, we explain here
how to apply it to the local models. First, if the group $G$ is split of
type $A$ or $C$, we find that all the $a_i^\vee=1$ in this case, and
the formulation of Theorem \ref{MainI} coincides with the original
conjecture of Pappas and Rapoport. Namely, the central charge of
$\calL(\sum_{i\in Y}\epsilon_i)$ is $\sharp Y$. In fact, in these
cases, it is proven in \emph{loc. cit.} (using the result of
\cite{Go1,Go2,PR0}) that the coherence conjecture holds for $\mu$
being sum of minuscule coweights. In what follows, we mainly discuss
the ramified unitary groups. As the main application of the coherence conjecture, general cases are treated in \cite{PZ}. 
 
Let us change the notation in the main body of the paper to the
following. Let $\calO_{F_0}$ be a completed discrete valuation ring
with algebraically closed residue field $k$ with $\on{char}k\neq 2$ and
fractional field $F_0$. Let $\pi_0$ be the uniformizer. For example,
$\calO=k[[t]]$ with $\pi_0=t$ as in the main body of the paper, or
$\calO=\bbZ_p^{ur}$, the completion of the maximal unramified
extension of $\bbZ_p$ and $\pi_0=p$.

We will follow \cite{PR2} (see also \cite{PR1}). Let $F/F_0$ be a
quadratic extension. Let $(V,\phi)$ be a split hermitian vector
space over $F$ of dimension $\geq 4$. That is, $V$ is a vector space over $F$ and $\phi$
is a hermitian form such that there is a basis $e_1,\ldots,e_n$ of
$V$ satisfying
\[\phi(e_i,e_{n+1-j})=\delta_{ij},\quad i,j=1,\ldots,n.\]
Let $G=\GU(V,\phi)$ be the group of unitary similitudes for
$(V,\phi)$, i.e. for any $F_0$-algebra $R$,
\[G(R)=\{g\in \GL(V\otimes_{F_0}R)|\phi(gv,gw)=c(g)\phi(v,w) \mbox{ for some } c(g)\in R^\times\}.\]
Then $G\otimes_{F_0}F\cong\GL_n\times\bbG_m$. The derived group
$G_\der$ is the ramified special unitary group $\SU(V,\phi)$
consisting of those $g\in G(R)$ such that $\det(g)=c(g)=1$.

We fix a square root $\pi$ of $\pi_0$. There are two associated
$F_0$-bilinear forms,
\[(v,w)=\Tr_{F/F_0}(\phi(v,w)),\quad\langle v,w \rangle=\Tr_{F/F_0}(\pi^{-1}\phi(v,w)).\]
Then $(-,-)$ is symmetric bilinear and $\langle-,-\rangle$ is alternating.  For $i=0,\ldots,n-1$, set
\[\Lambda_i=\on{span}_{\calO_{\tilde{F}}}\{\pi^{-1}e_1,\ldots,\pi^{-1}e_i,e_{i+1},\ldots,e_n\},\]
and complete this into a selfdual periodic lattice chain by setting
$\Lambda_{i+kn}=\pi^{-k}\Lambda_i$. Then
$\langle-,-\rangle:\Lambda_{-j}\times\Lambda_j\to\calO_{F_0}$ is a
perfect pairing. In particular,  $\Lambda_0$ is self-dual for the
alternating form $\langle-,-\rangle$.

Let us fix a minuscule coweight $\mu_{r,s}$ of $G_{F}$ of signature $(r,s)$ with $r\leq s, r+s=n$. That is
\[\mu_{r,s}(a)=(\on{diag}\{a^{(s)},1^{(r)}\},a)\]
where $a^{(s)}$ denotes $s$-copies of $a$. Let $E=F$ if $r\neq s$
and $E=F_0$ if $r=s$. Let $m=[\frac{n}{2}]$. Let
$I\subset\{0,\ldots,m\}$ be a non-empty subset with the requirement
that if $n$ is even and $m-1\in I$, then $m\in I$ as well (see
\cite[\S 1.b]{PR2} or \cite[Remark 4.2.C]{PR1} for the reason why we
make this assumption).

Let us define the following moduli scheme $\calM^{\on{naive}}_I$ over $\calO_{E}$.  A point of $\calM^{\on{naive}}_I$ with values in an $\calO_{E}$-scheme $S$ is given by an $\calO_{F}\otimes_{\calO_{F_0}}\calO_S$-submodule
$\calF_j\subset\Lambda_j\otimes_{\calO_{F_0}}\calO_S$ for each $j\in \pm I+n\bbZ$ satisfying the following conditions:
\begin{enumerate}
\item as an $\calO_S$-module, $\calF_j$ is locally on $S$ a direct summand of rank $n$;
\item for each $j<j', j,j'\in \pm I+n\bbZ$, the natural inclusion \[\Lambda_j\otimes_{\calO_{F_0}}\calO_S\to \Lambda_{j'}\otimes_{\calO_{F_0}}\calO_S\] induces a morphism $\calF_j\to \calF_{j'}$, and the isomorphism $\pi:\Lambda_j\to\Lambda_{j-n}$ induces an isomorphism of $\calF_j$ with $\calF_{j-n}$;
\item under the perfect pairing induced by $\langle-,-\rangle:\Lambda_{-j}\times\Lambda_j\to\calO_{F_0}$, $\calF_{-j}=\calF_j^\perp$, where $\calF_j^\perp$ is the orthogonal complement of $\calF_j$;
\item the determinant condition as in \cite[\S 1.e.1, d)]{PR2}.
\end{enumerate}
As explained in \emph{loc. cit.}, for any $I$,
$\calM_I^{\on{naive}}\otimes_{\calO_E}E$ is isomorphic to the
Grassmannian $\bbG(s,n)$ of $s$-planes in $n$-space. In addition,
for $i\in I$, there is a natural projection
$\calM_I^{\on{naive}}\to\calM_{\{i\}}^{\on{naive}}$ (if $n$ is even
and $i=m-1$, $\{i\}$ will mean $\{m-1,m\}$). Now the local model
$\calM_I^{\on{loc}}$ is defined as the flat closure of the generic
fiber $\calM_I^{\on{naive}}\otimes E$ inside $\calM_I^{\on{naive}}$.

The special fiber $\calM_I^{\on{naive}}\otimes k$ (and therefore
$\calM_I^{\on{loc}}\otimes k$) embeds into the (partial) affine flag
variety of the unitary group over $k\ppart$. Namely, let
$(V',\phi')$ be a split hermitian space over $k\pparu$ ($u^2=t$)
with a standard basis $e_1,\ldots,e_n$, such that
$\phi'(e_i,e_{n+1-j})=\delta_{ij}$. Let $\la_j,
j\in\{0,1,\ldots,n-1\}$, be the standard lattices in $V'$ defined
similarly to $\Lambda_j$ (replacing $\pi$ by $u$ and $\calO_{F}$ by
$k[[u]]$ in the definition of $\Lambda_j$). For
$I\subset\{0,\ldots,m\}$ as before, write $I={i_0<i_1<\cdots<i_k}$
and let $P_I$ be the group scheme over $k[[t]]$ which is the
stabilizer of the lattice chain \[\la_{i_0}\subset\cdots\subset
\la_{i_k}\subset u^{-1}\la_{i_0}\] in $\GU(V',\phi')$. As explained
in \emph{loc. cit.}, this is not always a connected group scheme
over $k[[t]]$. But if it is, then it is a parahoric group scheme of
$\GU(V',\phi')$. In any case, the neutral connected component
$P_I^0$ of $P_I$ is a parahoric group scheme.

Consider the ind-scheme $\calF_I$ which to a $k$-algebra $R$
associates the set of sequences of $R[[u]]$-lattice chains
\[L_{i_0}\subset\cdots\subset L_{i_k}\subset u^{-1}L_{i_0}\]
in $V'\otimes_{k\pparu}R\pparu$ together with an $R[[t]]$-lattice
$L\subset R\ppart$ satisfying conditions a) and b) as in \cite[\S
3.b]{PR2} (observe that we replace $\al\in
R\ppart^\times/R[[t]]^\times$ in \emph{loc. cit.} by a lattice
$L\subset R\ppart$, which seems more natural). Then \[\calF_I\cong
L\GU(V',\phi')/L^+P_I\] and $L\GU(V',\phi')/L^+P_I^0$ is either
isomorphic to $L\GU(V',\phi')/L^+P_I$ or to the disjoint union of
two copies of $L\GU(V',\phi')/L^+P_I$. In addition, for such $I$,
one can canonically associate to it a subset $Y\subset \bold S$
($\bold S$ are the set of vertices in the local Dynkin diagram of
$\GU(V',\phi')$) such that $\Fl^Y=L\GU(V',\phi')/L^+P^0_I$. Indeed,
by \cite[Remark 10.3]{PR1} (see also \cite[\S 1.2.3]{PR2}), one can
identify $\bold S$ with $\{0,1,\ldots,m\}$, if $n=2m+1$, resp.
$\{0,1,\ldots,m-2,m,m'\}$, if $n=2m$, where $m'$ is a formal symbol
as defined in \cite[\S 4]{PR1}, to which a lattice of $V'$
\[\la_{m'}=\on{span}_{k[[u]]}\{u^{-1}e_1,\ldots,u^{-1}e_{m-1},e_m,u^{-1}e_{m+1},e_{m+2},\ldots,e_{2m}\}\]
is associated. Then $Y=I$ in all cases except when
$n=2m,\{m-1,m\}\subset I$, in which case
$Y=(I\setminus\{m-1\})\cup\{m'\}$.

\begin{rmk}\label{special parahoric for unitary}
(i) Observe that if $n=2m+1$, under our identification of
$\{0,1,\ldots,m\}$ with $\bold S$ (the set of vertices of the local
Dynkin diagram), $i$ goes to the label $m-i$ in Kac's book
(\cite[p.p. 55]{Kac}), and if $n=2m$, under the identification of
$\{0,\ldots,m-2,m,m'\}$ with $\bold S$, $i$ goes to $m-i$ for $i\leq
m-2$ and $\{m,m'\}$ go to $\{0,1\}$.

(ii) As pointed out in \cite{PR1,PR2}, if $n=2m+1$, then $P_{\{0\}}$
and $P_{\{m\}}$ are the special parahoric group schemes, and if
$n=2m$, then $P_{\{m\}}, P_{\{m'\}}$ are the special parahoric group
schemes. We further point out: (1) let $n=2m+1$. Then (a)
$P_{\{0\}}$ is the special parahoric determined by a pinning of
$\GL_{2n+1}\times\bbG_m$, i.e. the group scheme $\calG_{v_0}$ as in
\eqref{constr of special}, and its reductive quotient is
$\on{GO}_{2m+1}$; and (b) the special parahoric $P_{\{m\}}$ has
reductive quotient $\on{GSp}_{2m}$, but it is not of the form
\eqref{constr of special}. (2) Let $n=2m$. Then both $P_{\{m\}},
P_{\{m'\}}$ are of the form \eqref{constr of special}, and their
reductive quotients are both isomorphic to $\on{GSp}_{2m}$.
\end{rmk}

Fix the isomorphisms $\Lambda_j\otimes_{k[[t]]} k\cong\la_j\otimes
k$, compatible with the actions of $\pi$ and $u$, by sending $e_i\to
e_i$. Now we embed the special fiber $\calM_I^{\on{naive}}\otimes k$
into $\calF_I$ as follows: for every $k$-algebra $R$,
\[\calF_j\subset(\Lambda_j\otimes k)\otimes_k R\cong (\lambda_j\otimes k)\otimes_k R,\]
and let $L_j\subset \lambda_j\otimes R[[t]]$ be the inverse image of
$\calF_j$ under $\lambda_j\otimes R[[t]]\to \lambda_j\otimes
\calO_S$. In addition, let $L=t^{-1}R[[t]]\subset R\ppart$. This
gives the embedding
\[\iota_I:\calM_I^{\on{naive}}\otimes k\to\calF_I.\]
It is proved in \cite[Proposition 3.1]{PR2} that
$\calA^{I}(\mu_{r,s})$ is contained in
$\calM_I^{\on{loc}}\otimes_{\calO_E}k$ under $\iota_I$, where $\calA^I(\mu_{r,s})$ is as defined in \eqref{Amu}.
Here we show the following result, which was shown in \cite[Theorem 0.1]{PR2} to follow from
a slightly different version of the coherence conjecture.
\begin{thm}\label{mainPR}
One has the equality $\calA^{I}(\mu_{r,s})=\calM_I^{\on{loc}}\otimes_{\calO_E}k$.
Therefore, the special fiber of $\calM_I^{\on{loc}}$ is reduced and
each irreducible component is normal, Cohen-Macaulay and
Frobenius-split.
\end{thm}
To prove it, one needs to construct a natural line bundle on
$\calM_I^{\on{naive}}$ and apply the coherence conjecture to compare
the dimensions of the space of global sections of this line bundle over the generic and the
special fibers. There are several choices of natural line
bundles. One of them will be given in \cite{PZ}, after we give a
group theoretical description of $\calM_I^{\on{naive}}$. Here, we
follow the original approach of \cite{PR1,PR2} to construct
another line bundle $\calL_I$, which is more down to earth.

First, if $I=\{j\}$, we define the line bundle $\calL_{\{j\}}$ over
$\calM_{\{j\}}^{\on{naive}}$ whose value at the $\calO_S$-point
given by $\calF_j\subset\Lambda_j\otimes_{\calO_E}\calO_S$ is
$\det(\calF_j)^{-1}$. If $n=2m$, we also define $\calL_{\{m-1,m\}}$
over $\calM_{\{m-1,m\}}^{\on{naive}}$ whose value at the
$\calO_S$-point of given by $\calF_{m-1}\subset\calF_m$ is
$\det(\calF_{m-1})^{-1}\otimes\det(\calF_m)^{-1}$.  For general $I$,
the line bundle $\calL_I$ is defined as the tensor product of those
$\calL_{\{j\}}$ or $\calL_{\{m-1,m\}}$ along all possible
projections $\calM^{\on{naive}}_I\to\calM^{\on{naive}}_{\{j\}}$ or
$\calM^{\on{naive}}_I\to\calM^{\on{naive}}_{\{m-1,m\}}$.

The restriction of $\calL_{\{j\}}$ to the generic fiber
$\calM_{\{j\}}^{\on{naive}}\otimes_EF\cong \Gr(s,n)$ is isomorphic
to $\calL_{\det}^{\otimes 2}$, where $\calL_{\det}$ is the
determinant line bundle on $\Gr(s,n)$, which is the positive
generator of the Picard group of $\Gr(s,n)$. On the other hand,
the restriction of $\calL(\{m-1,m\})$ to the generic fiber of
$\calM_{\{m-1,m\}}^{\on{naive}}$ is isomorphic to
$\calL_{\det}^{\otimes 4}$. Recall the $\calA^I(\mu_{r,s})^\circ$ defined in \eqref{nAmu}, and recall that there is the canonical isomorphism 
$\calA^I(\mu_{r,s})^\circ\cong
\calA^I(\mu_{r,s})$ as $G_\der=\SU_n$ is simply-connected.

\begin{prop} Under the canonical isomorphism $\calA^I(\mu_{r,s})^\circ\cong
\calA^I(\mu_{r,s})$, the line bundle $\calL_I$, when restricted to
$\calA^I(\mu_{r,s})$ is isomorphic to the restriction of
$\calL(\sum_{j\in Y}\kappa(j)\epsilon_j)$ to $
\calA^I(\mu_{r,s})^\circ$, where
\begin{enumerate}
\item if $n=2m+1$, then $\kappa(j)=1$ for
$j=0,1,\ldots,m-1$ and $\kappa(m)=2$;
\item if $n=2m$, then $\kappa(j)=1$ for $j=0,\ldots,m-2$ and
$\kappa(m)=\kappa(m')=2$.
\end{enumerate}
\end{prop}
\begin{proof}
Let us first introduce a convention. In what follows, when we write
$\la_j$, we consider it as a $k[[u]]$-lattice. If we just remember
its $k[[t]]$-lattice structure, we denote it by $\la_j/k[[t]]$.

Clearly, we can assume that $I=\{j\}$ or when $n=2m$ we shall also
consider $I=\{m-1,m\}$. The latter case will be treated at the end
of the proof. So we first assume that $j\neq m-1$.

Observe that we have a natural closed embedding of ind-schemes
\[L\GU(V',\phi')/L^+P_{\{j\}}\cong\calF_{\{j\}}\to\Gr_{\GL(\la_j)}\times\Gr_{\bbG_m}\]
by just remembering the lattices
$L_j\subset\la_j\otimes_{k\pparu}R\pparu$ and $L\subset R\ppart$. By
definition, the line bundle $\iota_{\{j\}}^*\calL_{\{j\}}$ on
$\calM_{\{j\}}^{\on{naive}}\otimes_{\calO_E} k$ is the pullback of
the determinant line bundle on $\Gr_{\GL(\la_j)}$ along the above
map.

Let $\SU(V',\phi')$ be the special unitary group. As explained in
\cite[\S 4]{PR1}, $P'_I=P_I\cap\SU(V',\phi')$ is a parahoric group
scheme of $\SU(V',\phi')$. By \cite[\S 6]{PR1}, we have
\[\begin{CD}L\SU(V',\phi')/L^+P'_{\{j\}}@>>>L\GU(V',\phi')/L^+P_{\{j\}}\\
@VVV@VVV\\
\Gr_{\SL(\la_j)}@>>>\Gr_{\GL(\la_j)}
\end{CD}\]
where the ind-schemes in the left column are identified with the
reduced part of neutral connected components of the ind-schemes in
the right column. Since the isomorphism
$\calA^I(\mu_{r,s})^\circ\cong \calA^I(\mu_{r,s})$ is obtained from
the translation by some $g\in \GU(V',\phi')(F)$, it is enough to
prove
\begin{lem}The pullback of $\calL_{\det}$ by
$L\SU(V',\phi')/L^+P'_{\{j\}}\to\Gr_{\SL(\la_j)}$ is
$\calL(\kappa(j)\epsilon_j)$.
\end{lem}
\begin{proof}
Assume that $j\neq 0,m$, and in the case $n=2m$, $j\neq m-1$. By \eqref{rev:pic}, the pullback of $\calL_{\det}$ is of the form $\calL(m\epsilon_j)$ for some $m$. Consider the rational line
$\bbP^1_j\subset L\SU(V',\phi')/L^+P'_{\{j\}}$ given by the
$\bbA^1=\Spec k[s]$-family of lattices
\[L_{j,s}=u^{-1}k[[u]]e_1+\cdots+ u^{-1}k[[u]]e_{j-1}+u^{-1}k[[u]](e_j+s e_{j+1}) +k[[u]]e_{j+1}+\cdots+k[[u]]e_n.\]
It is easy to see that the restriction of $\calL_{\det}$ to this
rational line is $\calO(1)$. In fact, by the map
\[L_{j,s}\to L_{j,s}/(\sum_{r\leq j-1} u^{-1}k[[u]]e_r+\sum_{r\geq
j} k[[u]]e_r),\] this rational curve $\bbP^1_j$ is identified with
the $\Gr(1,2)$ classifying lines in the 2-dimensional $k$-vector
space generated by $\{u^{-1}e_j, u^{-1}e_{j+1}\}$ and clearly the
restriction of the determinant line bundle of $\Gr_{\SL(\la_j)}$ is
the determinant line bundle on $\Gr(1,2)$. Therefore, $\kappa(j)=1$
if $j\neq 0,m$ (and $j\neq m-1$ if $n=2m$).

If $j=0$, consider the rational line $\bbP^1_0\subset
L\SU(V',\phi')/L^+P'_{\{0\}}$ given by the $\bbA^1=\Spec
k[s]$-family of lattices
\begin{equation}\label{rational
line}L_s=k[[u]]e_1+\cdots+k[[u]]e_{n-1}+k[[u]](e_n+s
u^{-1}e_1).
\end{equation}
By the same reasoning as above, the restriction of $\calL_{\det}$ to
this rational line is $\calO(1)$. Therefore, $\kappa(0)=1$.

Now, if $n=2m+1$ and $j=m$ or $n=2m$ and  $j=m$ or $m'$, we will
prove that $2\mid \kappa(j)$. Assuming this, to prove the lemma it
is enough to find some rational line $\bbP^1_j\subset
L\SU(V',\phi')/L^+P'_{\{j\}}$ such that the restriction of
$\calL_{\det}$ to it is $\calO(2)$. If $n=2m+1$, we can take the
rational line $\bbP^1_m$ given by the $\bbA^1=\Spec k[s]$-family of
lattices
\begin{multline*}L_s=u^{-1}k[[u]]e_1+\cdots+u^{-1}k[[u]]e_{m-1}+\\
u^{-1}k[[u]](e_m+s
e_{m+1}-\frac{s^2}{2}e_{m+2})+k[[u]]e_{m+1}+\cdots+k[[u]]e_n.\end{multline*}
To see that $\calL_{\det}$ restricts to $\calO(2)$, consider the map
\[L_s\to L_s/(\sum_{r\leq m-1} u^{-1}k[[u]]e_r+\sum_{r\geq
m} k[[u]]e_r),\] which gives rise to embeddings,
$\bbP^1_m\subset\Gr(1,3)\subset\Gr_{\SL(\la_m)}$. Here $\Gr(1,3)$
classifies lines in the 3-dimensional $k$-vector space generated by
$\{u^{-1}e_m,u^{-1}e_{m+1},u^{-1}e_{m+2}\}$. Clearly, the pullback
of $\calL_{\det}$ along $\Gr(1,3)\to\Gr_{\SL(\la_m)}$ is the
determinant line bundle and the embedding $\bbP^1_m\to\Gr(1,3)$ is
quadratic, the claim follows.

If $n=2m$ and $j=m$ (the case $j=m'$ is similar), we can take the
rational line $\bbP^1_{m}$ given by the $\bbA^1=\Spec k[s]$-family
of lattices
\begin{multline*}L_s=u^{-1}k[[u]]e_1+\cdots+u^{-1}k[[u]]e_{m-2}+u^{-1}k[[u]](e_{m-1}+s e_{m+1})\\
+u^{-1}k[[u]](e_m-se_{m+2})+k[[u]]e_{m+1}+\cdots+k[[u]]e_n.\end{multline*}
To see that $\calL_{\det}$ restricts to $\calO(2)$, consider the map
\[L_s\to L_s/(\sum_{r\leq m-2} u^{-1}k[[u]]e_r+\sum_{r\geq
m-1} k[[u]]e_r),\] which gives rise to embeddings,
$\bbP^1_{m}\subset\Gr(2,4)\subset\Gr_{\SL(\la_m)}$. Here $\Gr(2,4)$
classifies planes in the 4-dimensional $k$ vector space generated by
$\{u^{-1}e_{m-1},\ldots,u^{-1}e_{m+2}\}$. The restriction of
$\calL_{\det}$ to $\Gr(2,4)$ is the determinant line bundle, and
therefore it is enough to see that the restriction of the
determinant line bundle on $\Gr(2,4)$ along $\bbP^1_m\to\Gr(2,4)$ is
$\calO(2)$. If we use the determinant line bundle on $\Gr(2,4)$ to
embed $\Gr(2,4)$ into $\bbP(V)$, where $V$ is generated by
$\{\{u^{-1}e_i\wedge u^{-1}e_j\mid m-1\leq i<j\leq m+2\}$, then the
composition $\Spec k[s]\subset
\bbP^1_m\to\Gr(2,4)\to\bbP(V)\setminus\{u^{-1}e_{m-1}\wedge
u^{-1}e_m\}$ is given by
\[s\mapsto su^{-1}e_{m-1}\wedge u^{-1}e_{m+2}+su^{-1}e_m\wedge u^{-1}e_{m+1}-s^2u^{-1}e_{m+1}\wedge u^{-1}e_{m+2}.\]
The claim is clear from this description.

So it remains to prove $2\mid \kappa(j)$ for $n=2m+1, j=m$, or
$n=2m, j=m$ or $m'$. Recall that when regarding $V'$ as a vector
space over $k\ppart$, it has a split symmetric bilinear form
\[(v,w)=\Tr_{k\pparu/k\ppart}(\phi'(v,w)).\]
Observe that when $n=2m+1, j=m$, or $n=2m, j=m$ or $m'$,
$\la_j/k[[t]]$ is maximal isotropic, i.e. $\la_j\subset\widehat{\la_j}^s$
and $\dim_k(\widehat{\la_j}^s/\la_j)=0$ or $1$, where
\[\widehat{\la_j}^s=\{v\in V'\mid (v,\la_j)\subset\calO\}.\]
Let $\on{Iso}(V')\subset\Gr_{\SL(\la_j/k[[t]])}$ denote the subspace
of maximal isotropic lattices in $V'$. Then the morphism
\[L\SU(V',\phi')/L^+P'_{\{j\}}\to \Gr_{\SL(\la_j)}\to \Gr_{\SL(\la_j/k[[t]])}\]
factors through
\[L\SU(V',\phi')/L^+P'_{\{j\}}\to \on{Lag}(V')\to \Gr_{\SL(\la_j/k[[t]])}.\]
It is by definition that the pullback of $\calL_{\det}$ along
$\Gr_{\SL(\la_j)}\to \Gr_{\SL(\la_j/k[[t]])}$ is $\calL_{\det}$, and
it is well-known (for example, see \cite[\S 4]{BD}) that the
pullback of $\calL_{\det}$ along
$\on{Iso}(V')\to\Gr_{\SL(\la_j/k[[t]])}$ admits a square root (the
Pffafian line bundle). The lemma follows.
\end{proof}

To deal with the case $n=2m$ and $I=\{m-1,m\}$, observe there is a
map
\[LGU(V',\phi')/L^+P_{I}\to\Gr_{\GL(\la_{m})}\times\Gr_{\GL(\la_{m'})}\]
by sending $L_{m-1}\subset L_m$ to $L_m, gL_m$, where $g$ is the
unitary transformation $e_m\mapsto e_{m+1}, e_{m+1}\mapsto e_m$ and
$e_i\mapsto e_i$ for $i\neq m,m+1$. One observes that
$\iota_{I}^*\calL_I$ on $\calM_{I}^{\on{naive}}\otimes_{\calO_E} k$
is the pullback along the above map of the tensor product of the
determinant line bundles (on each factor).
\end{proof}

Finally, let us see why this proposition can be used to deduce
Theorem 0.1 of \cite{PR2}. First let $a_i^\vee$ be the Kac labeling
as in \cite[\S 6.1]{Kac}. Using Remark \ref{special parahoric for
unitary} (i), by checking all the cases, we find that
$a_i^\vee\kappa(i)=2$. Let $\calL_I$ be the line bundle on
$\calM_I^{\on{loc}}$. Then for $a\gg 0$,
\begin{equation}\label{ap:eq}
\dim\Ga(\calM_I^{\on{loc}}\otimes_{\calO_E}k, \calL_I^a)=\dim\Ga(\calM_I^{\on{loc}}\otimes_{\calO_E}E,\calL_I^a).
\end{equation}
By the above proposition and \cite[Proposition 3.1]{PR2}, the left
hand side
\[\dim\Ga(\calM_I^{\on{loc}}\otimes_{\calO_E}k, \calL_I^a)\geq \dim\Ga(\calA^I(\mu_{r,s})^\circ,\calL(a\sum_{i\in
Y}\kappa(i)\epsilon_i)),\] and the central charge of
$\calL(a\sum_{i\in Y}\kappa(i)\epsilon_i)$ is \[\sum_{i\in
Y}aa_i^\vee\kappa(i)=2a\sharp I.\] The line bundle on right hand
side of \eqref{ap:eq} is just the $2a\sharp I$-power of the ample generator of the
Picard group of $\bbG(s,n)$. Then since Theorem \ref{MainI} holds, Theorem
\ref{mainPR} follows by the argument in \cite{PR1}.

\section{Appendix II: Some recollections and proofs}\label{recol}
We collect  and strengthen various results, which exist in literature, in the form needed in the main body
of the paper.

\subsection{Combinatorics of Iwahori-Weyl group}\label{cb of IW}
We recall a few facts about the translation elements in the Iwahori-Weyl group which
are used in the paper. We keep the notation as in \S \ref{grp data}. In particular, we identify the apartment $A(G,S)$ with $\xcoch(S)_{\bbR}$ via the special vertex $v$. We choose the alcove $\mathbf{a}$, whose closure $\overline{a}$ contains $v$, and whch is contained in the finite Weyl chanber of $G$ determined by the chosen Borel subgroup. We write $\widetilde{W}=\xcoch(T)_\Ga\rtimes W_0$ using the vertex $v$.

Let $2\rho$ be the sum of all positive roots (for $H$). Observe that
given $\mu\in\xcoch(T)_\Ga$,  the integer $(\tilde{\mu},2\rho)$ is independent of  its lifting
$\tilde{\mu}\in\xcoch(T)$. By abuse of notation, we denote this
number by $(\mu,2\rho)$.

\begin{lem}\label{length}Let $\mu\in\xcoch(T)_\Ga^+$, the set defined in \eqref{plus}.
Let $\Lambda=\subset\xcoch(T)_\Ga$ be the $W_0$-orbit associated to $\mu$ as in \S \ref{grp data}. Then for all $\nu\in \Lambda$,
$\ell(t_\nu)=(2\rho,\mu)$.
\end{lem}
\begin{proof}
(Let $x\in\mathbf{a}$ be a point in the interior of the alcove $\mathbf{a}$. Then for any $w\in\widetilde{W}$, 
\begin{equation}\label{length formula}
\ell(w)=\{\alpha \mbox{ is an affine root }\mid \alpha(x)>0, \alpha(w(x))<0\}.
\end{equation}
If $w=t_\nu$ is a translation element, then this is the number of affine roots $\alpha$ such that $0< \alpha(x)<(\dot\alpha,\nu)$, where $\dot\alpha$ is the vector part of $\alpha$ (so $\dot\alpha$ is a finite root of $G$). This number can be rewritten  as
\[\sum_{a\in\Phi, (a,\nu)>0}\sharp\{\alpha\mid  \dot{\alpha}=a, 0<\alpha(x)<(a,\nu) \}.\]
Let $j:\Phi(H,T_H)\to \Phi(G,S)$ be the restriction of the root system of $H$ (the absolute root system) to the root system of $G$ (the relative root system) . Then the lemma will follow from the equality 
\[\sharp\{\alpha\mid  \dot{\alpha}=a, 0<\alpha(x)<(a,\nu) \}=\sum_{\tilde{a}\in j^{-1}(a)}(\tilde{a},\tilde{\nu}).\]
This statement involves only one root of $G$. By checking the semi-simple subgroup of $G$ of semi-simple $F$-rank one (which are Weil restrictions of either $\SL_2$ or $\SU_3$), this equality holds.
\end{proof}

One can easily deduce the following lemma from \eqref{length formula}.
\begin{lem}\label{length criterion}
Let $w,w'\in\tilde{W}$.
Then $\ell(ww')=\ell(w)+\ell(w')$ if and only if the following two
statements holds: for $a+m$ an affine root,
\begin{enumerate}
\item if $\alpha(x)>0$ and $\alpha(w'(x))<0$, then $\alpha(ww'(x))<0$;
\item if $\alpha(x)>0$ and $\alpha(w'^{-1}(x))<0$, then $\alpha(w(x))>0$.
\end{enumerate}
\end{lem}

\begin{lem}\label{length2}
Let $\mu\in\xcoch(T)_\Ga^+$. Then
$\ell(t_\mu w_f)=\ell(t_\mu)+\ell(w_f)$. 
\end{lem}
\begin{proof}
Let $w=t_\mu$ with $\mu\in\xcoch(T)_\Ga^{+}$ and $w_f\in{W_0}$.
Assume that $\alpha(x)>0$ and $\alpha(w_f(x))<0$. As $v=\overline{\mathbf{a}}\cap\overline{w_f(\mathbf{a})}$, $\alpha(v)=0$. Let
$a=\dot{\alpha}$, then $a(x-v)=\alpha(x)>0$, i.e. $a$ is a positive root of $G$. Therefore $\alpha(t_\mu w_f(x))=\alpha(w_f(x))-(\mu,a)<0$. On the other hand,
assume that $\alpha(x)>0$ and $\alpha(w_f^{-1}(x))<0$. Then $\alpha(v)=0$,
and $w_f(a)$ is negative. Therefore $(\mu,a)\leq 0$. Then $\alpha(t_\mu(x))=\alpha(x)-(\mu,a)>0$. This proves
that $\ell(t_\mu w_f)=\ell(t_\mu)+\ell(w_f)$.
\end{proof}

On the finitely generated abelian group $\xcoch(T)_\Ga$, there are
two partial orders. One is the restriction of the Bruhat order on
$\widetilde{W}$, denoted by ``$\leq$". The other, denoted by
``$\preceq$", is defined as follows. Recall that the lattice
$\xcoch(T_\s)$ is the coroot lattice of $H$. The Galois group $\Ga$
acts on $\xcoch(T_\s)$ which sends the positive coroots of $H$
(determined by the pinning) to positive coroots. Therefore, it makes
sense to talk about positive elements in $\xcoch(T_\s)_\Ga$. Namely,
$\la\in\xcoch(T_\s)_\Ga$ is positive if its preimage in
$\xcoch(T_\s)$ is a sum of positive coroots (of $(H,T_H)$). Since
$\xcoch(T_\s)_\Ga\subset\xcoch(T)_\Ga$, we can define
$\la\preceq\mu$ if $\mu-\la$ is positive in $\xcoch(T_\s)_\Ga$.

\begin{lem}\label{order}
Let $\la,\mu\in\xcoch(T)^+_\Ga$. Then $\la\preceq\mu$ if and only if
$t_\la\leq t_\mu$ in the Bruhat order.
\end{lem}
\begin{proof}In the case that $G$ is split, the proof is contained
in \cite[Proposition 3.2, 3.5]{R}. The ramified case can be reduced
to the same proof as shown in \cite[Corollary 1.8]{Ri}. See
\cite[Remark 4.2.7]{PRS}.
\end{proof}

Recall the following lemma.
\begin{lem}
Let $x,y\in \widetilde{W}$ and $w\in W_{\on{aff}}$. Assume that $\ell(xw)=\ell(x)+\ell(w)$ and $\ell(yw)=\ell(y)+\ell(w)$. Then $x\leq y$ if and only if $xw\leq yw$.
\end{lem}
\begin{proof}By induction of the length of $w$, we can assume that $w$ is a simple reflection. Then the lemma is clear.
\end{proof}

\begin{lem}\label{other chamber}
Let $\la,\mu\in w(\xcoch(T)_\Ga^{+})$. Then $t_\la\leq t_\mu$ if and only if $t_{w^{-1}\la}\leq t_{w^{-1}\mu}$.
\end{lem}
\begin{proof}
Observe that $w^{-1}\la$ and $w^{-1}\mu$ are dominant. Combining Lemma \ref{length} and \ref{length2}, 
\[\ell(w^{-1}t_\la)=\ell(t_{w^{-1}\la})+\ell(w^{-1})=\ell(w^{-1})+\ell(t_\la).\]
Therefore by the above lemma, $t_{w^{-1}\la}\leq t_{w^{-1}\mu}$ if and only if $w^{-1}t_{\la}\leq w^{-1}t_{\mu}$ if and only if $t_{\la}\leq t_{\mu}$.
\end{proof}

\subsection{Deformation to the normal cone}\label{deform}
Let $C$ be a smooth curve over an algebraically closed field $k$.
Let $\calX$ be a scheme faithfully flat and affine over $C$. Let
$x\in C(k)$ be a point and let $\calX_x$ denote the fiber of $\calX$
over $x$. Let $Z\subset\calX_x$ be a closed subscheme. Consider the
following functor $\calX_Z$ on the category of flat $C$-schemes: for
each $V\to C$,
\[\calX_Z(V)=\{f\in\Hom_C(V,\calX)\mid f_x:V_x\to\calX_x \mbox{ factors through } V_x\to Z\subset \calX_x\}.\]
It is well known that this functor is represented by a scheme
affine and flat over $C$, usually called the deformation to the
normal cone (or called the dilatation of $Z$ on $\calX$, see
\cite[\S 3.2]{BLR}). Indeed, the construction is easy if $\calX$ is
affine over $C$. Namely, we can assume that $C$ is affine and $x$ is
defined by a local parameter $t$. Assume that $\calA$ be the
$\calO_C$-algebra defining $\calX$ over $C$, and let
$\calI\subset\calA$ be the ideal defining $Z\subset \calX$. Then
$t\calA\subset \calI$ and let $\calB=\calA[\frac{i}{t}, i\in
\calI]\subset\calA[t^{-1}]$. It is easy to see that $\calB$ is flat
over $\calO_C$ and $\Spec\calB$ represents $\calX_Z$.

There is a natural morphism $\calX_Z\to\calX$ which induces an
isomorphism over $C-\{x\}$ and over $x$ it factors as
$(\calX_Z)_x\to Z\to \calX_x$. If $\calX$ is smooth over $C$, and
$Z$ is a smooth closed subscheme of $\calX_x$, then $\calX_Z$ is
also smooth over $C$. Indeed, \'{e}tale locally on $\calX_x$, the
map $(\calX_Z)_x\to Z$ can be identified with the map from the
normal bundle of $Z$ inside $\calX_x$ to $Z$, which justifies the
name of the construction.

Now let $\calG_1$ be a connected affine smooth group scheme over the
curve $C$. Let $x\in C(k)$ and let $(\calG_1)_x$ be the fiber of
$\calG_1$ at $x$. Let $P\subset(\calG_1)_x$ be a smooth closed
subgroup. Let $\calG_2=(\calG_1)_P$. This is indeed a smooth
connected affine group scheme over $C$. By restriction to $x$, we
have $r:\Bun_{\calG_2}\to\bbB (\calG_2)_x$ and
$r:\Bun_{\calG_1}\to\bbB (\calG_1)_x$ (here we assume that
$C$ is a complete curve).
\begin{prop}\label{cart}
We have the following Cartesian diagram
\[\xymatrix{
\Bun_{\calG_2}\ar^r[r]\ar[d]&\bbB(\calG_2)_x\ar[r]&\bbB P\ar[d]\\
\Bun_{\calG_1}\ar^r[rr]&&\bbB(\calG_1)_x }\]
\end{prop}
\begin{proof}Let $V=\Spec R$ be a noetherian\footnote{This suffices since all the stacks are locally of finite presentation.} affine scheme. Let $\calE$ be a $\calG_1$-torsor on $C_R$ and
$\calE_P$ be a $P$-torsor on $V$ together with an isomorphism
$\calE_P\times^P(\calG_1)_x\cong\calE|_{\{x\}\times \Spec R}$. We
need to construct a $\calG_2$-torsor $\calE'$ satisfying the appropriate
conditions. This construction will provide the inverse to the morphism $\Bun_{\calG_2}\to \Bun_{\calG}\times_{\bbB(\calG_1)_x}\bbB P$.

As a scheme over $C$, $\calE$ is faithfully flat. Its fiber over $x$
is $\calE|_{\{x\}\times \Spec R}$. Let $Z$ be the closed subscheme
of $\calE_x$ given by the closed embedding
\[\calE_P\subset\calE_P\times^P(\calG_1)_x\cong\calE|_{\{x\}\times \Spec R}.\]
Then $\calE_Z$ is a scheme affine and flat over $C$, together with a
morphism $\calE_Z\to \calE$. Therefore, $\calE_Z$ is a scheme over
$C_R$. We claim that $\calE_Z$ is a $\calG_2$-torsor over $C_R$.
First, $\calE_Z$ is faithfully flat over $C_R$. Indeed, by the local
criterion of flatness, it is enough to prove that
$\calE_Z|_{\{x\}\times \Spec R}$ is faithfully flat over $\Spec R$.
But this is clear, since \'{e}tale locally on $\calE|_{\{x\}\times
\Spec R}$, there is an isomorphism between $\calE_Z|_{\{x\}\times
\Spec R}$ and the normal bundle of
$\calE_P\subset\calE_P\times^P(\calG_1)_x$. Next, there is an action
of $\calG_2$ on $\calE_Z$. Indeed, the map
$\calE_Z\times_{C_R}\calG_2\mapsto
\calE\times_{C_R}\calG_1\to\calE$, when restricted to the fiber over
$x$, factors through $Z$. Therefore, by the definition of $\calE_Z$,
it gives rise to a map \[\calE_Z\times_{C_R}\calG_2\mapsto
\calE_Z.\] Finally, it is easy to see that
\[\calE_Z\times_{C_R}\calE_Z\cong\calE_Z\times_{C_R}\calG_2.\] Indeed,
the left hand side represents the scheme
$(\calE\times_{C_R}\calE)_{Z\times_{\Spec R}Z}$ and the right hand side
represents the scheme $(\calE\times_{C_R}\calG_1)_{Z\times _{\Spec
R} P}$. Then the desired isomorphism follows from
\[(\calE\times_{C_R}\calE)_{Z\times_{\Spec R} Z}\cong
(\calE\times_{C_R}\calG_1)_{Z\times_{\Spec R} P}.\]
\end{proof}

\subsection{Frobenius morphisms}\label{Frob}
Let us review some basic facts about the Frobenius morphisms of a
variety $X$ over an algebraically closed field of characteristic
$p>0$. The book \cite[Chapter 1]{BK} provides a detailed account of the general theory.

First assume that $X$ is smooth and let $\omega_X$ be its
canonical sheaf. Then there is the following isomorphism (\cite[\S 1.3.7-1.3.8]{BK})
\begin{equation}\label{dual}
\frakD:F_*\omega_X^{1-p}\stackrel{\sim}{\to}\calH om_{\calO_X}(F_*\calO_X,\calO_X),
\end{equation}
where $F:X\to X$ is the absolute Frobenius map of $X$. The existence of this isomorphism follows from the Grothendieck duality theorem for finite morphisms (see \cite{BK}, the discussion before \S 1.3.1).  Explicitly, the isomorphism is given as follows. Let $x\in X$ be a closed point and let
$x_1,\ldots, x_n$ be a sequence of regular parameters of the complete local ring $\hat{\calO}_{X,x}$. Then in an \'{e}tale neighborhood of $x$ in $X$, the above isomorphism is given by
\begin{equation}\label{Cartier isom}
\begin{split}&\frakD(x_1^{m_1}\cdots x_n^{m_n}(dx_1\cdots dx_n)^{1-p})(x_1^{\ell_1}\cdots x_n^{\ell_n})\\
=&\left\{\begin{array}{ll}0 &\mbox{if } p\nmid m_i+\ell_i+1 \mbox{
for some } i\\ x_1^{(m_1+\ell_1-p+1)/p}\cdots
x_n^{(m_n+\ell_n-p+1)/p}\end{array}\right.\end{split}
\end{equation}

Next, assume that $X$ is normal (\cite[\S 1.3.12]{BK}). It is still make sense to talk
about the canonical sheaf $\omega_X$ and its any $n$th power
$\omega_X^{[n]}$. Namely, let $j:X^{sm}\to X$ be the open immersion
of the smooth locus into $X$. Then by definition
$\omega_X^{[n]}:=j_*\omega_{X^{sm}}^n$. The isomorphism \eqref{dual}
still holds in this situation. Observe that there is a natural map
$(\omega_X^{[\pm 1]})^{\otimes n}\to \omega_X^{[\pm n]} (n>0)$ which
is not necessarily an isomorphism. In what follows, we use
$\omega_X^n$ to denote $\omega_X^{[n]}$ if no confusion will rise.
Let us recall that if in addition $X$ is Cohen-Macaulay, $\omega_X$
is the dualizing sheaf.

Next, we consider a flat family $f:X\to V$ of varieties which is
fiberwise normal and Cohen-Macaulay. In addition, let us assume that
$V$ is smooth, so that the total space $X$ is also normal and
Cohen-Macaulay. In this case, the relative dualizing sheaf
$\omega_{X/V}$ commutes with base change and is flat over $V$. We
have $\omega_X\cong f^*\omega_V\otimes\omega_{X/V}$. Let $X^{(p)}$
be the Frobenius twist of $X$ over $V$, i.e. the pullback of $X$
along the absolute Frobenius endomorphism $F:V\to V$. Let
$F_{X/V}:X\to X^{(p)}$ be the relative Frobenius morphism, and let
$\varphi:X^{(p)}\to X$ be the map such that the composition
$\varphi\circ F_{X/V}$ is the absolute Frobenius morphism $F$ for
$X$. Then
\begin{equation}\label{rel dual}
\frakR\frakD:(F_{X/V})_*\omega_{X/V}^{1-p}\stackrel{\sim}{\to}\calH
om_{\calO_{X^{(p)}}}((F_{X/V})_*\calO_X,\calO_{X^{(p)}}).
\end{equation}
Here $\omega_{X/V}^{n}$, as in the absolute case, is the pushout of
the $n$th tensor power of the relative canonical sheaf on
$X^{rel,sm}$, the maximal open part of $X$ such that
$f|_{X^{rel,sm}}$ is smooth. In addition, we have the following
homorphisms
\begin{equation}\label{comp1}
f^*F_*\omega_V^{1-p}\stackrel{\sim}{\to} f^*\calH om(F_*\calO_V,\calO_V)\cong\calH
om(\varphi_*\calO_{X^{(p)}},\calO_X),
\end{equation}
\begin{equation}\label{comp2}
 F_*\omega_{X/V}^{1-p}\otimes
f^*F_*\omega_V^{1-p}\cong  F_*\omega_{X/V}^{1-p}\otimes
F_*f^*\omega_V^{1-p}\to F_*\omega_X^{1-p}.
\end{equation}
The homorphisms \eqref{dual}, \eqref{rel dual}-\eqref{comp2} fit into
the following commutative diagram
\begin{equation}\label{comm diag}\begin{CD}\varphi_*\calH om((F_{X/V})_*\calO_X,\calO_{X^{(p)}})@.\ \otimes_{\calO_X}\ @.\calH om(\varphi_*\calO_{X^{(p)}},\calO_X)@>>>\calH om(F_*\calO_X,\calO_X)\\
@AAA@.@AAA@AAA\\
F_*\omega_{X/V}^{1-p}@.\ \otimes_{\calO_X}\
@.f^*F_*\omega_V^{1-p}@>>>F_*\omega_X^{1-p}.
\end{CD}\end{equation}

Finally, let $W$ be another smooth variety over $k$ and let $g:W\to
V$ be a $k$-morphism (not necessarily flat). By abuse of 
notation, we still use $g$ to denote the base change maps $X_W\to X$
and $(X_W)^{(p)}\cong X_W^{(p)}\to X^{(p)}$. Then the following
diagram is commutative.
\begin{equation}\label{base change for F-split}
\begin{CD}
g^*(F_{X/V})_*\omega_{X/V}^{1-p}@>\cong>>g^*\calH
om((F_{X/V})_*\calO_X,\calO_{X^{(p)}})\\
@V\cong VV@VV\cong V \\
(F_{X_W/W})_*\omega_{X_W/W}^{1-p}@>\cong>>\calH
om((F_{X_W/W})_*\calO_{X_W},\calO_{X_W^{(p)}}).
\end{CD}
\end{equation}

To prove the isomorphism \eqref{rel dual}, and that \eqref{comm diag} and \eqref{base change
for F-split} are commutative, one can first assume that $X$ is smooth over $V$. In this case, the proof of \eqref{dual} (as in \cite[\S 1.3]{BK}) with obvious modifications applies to \eqref{rel dual}. In particular, \'{e}tale locally on $X$, \eqref{rel dual} can be described by the explicit formula as in \eqref{Cartier isom}, with $x_1,\ldots,x_n$ replaced by a system of local coordinates of $X$ relative to $V$. Then \eqref{comm diag} and \eqref{base change for F-split} follows from the direct calculation.
Then one can easily extends these to the case that $X$ is flat over $V$ with with normal and
Cohen-Macaulay fibers. Indeed, under our assumptions, all the sheaves appearing in \eqref{rel dual}, \eqref{comm diag} and \eqref{base change
for F-split} have the following property: Let $\calF$ be such a sheaf on $X$ and $j:X^{rel,sm}\to X$ be the open embedding as before, then $\calF\cong j_*(\calF|_{X^{rel,sm}})$.

\end{document}